\documentclass[11pt]{article}

\usepackage[T1]{fontenc}
\usepackage[latin1]{inputenc}
\usepackage{amsmath,amsthm,amssymb}

\usepackage{authblk}

\usepackage[margin=1in]{geometry}

\usepackage{hyperref}

\newtheorem{prop}{Proposition}[section]
\newtheorem{cor}[prop]{Corollary}
\newtheorem{defi}[prop]{Definition}
\newtheorem{lem}[prop]{Lemma}
\newtheorem{theo}[prop]{Theorem}

\newtheorem{theomain}{Theorem}

\def\ricc{\mbox{\rm Ricc}}
\def\tr{\mbox{\rm tr}}
\def\lip{\mbox{\rm lip}} 
    
\def\proj{\mbox{\rm proj}}
\newcommand{\point}{\mbox{\LARGE .}}
\DeclareMathOperator*{\argmin}{\arg\!\min}
\DeclareMathOperator*{\argmax}{\arg\!\max}

\newcommand{\CC}{\mathbb{C}}

\newcommand{\EE}{\mathbb{E}}

\newcommand{\LL}{\mathbb{L}}
\newcommand{\MM}{\mathbb{M}}

\newcommand{\RR}{\mathbb{R}}
\newcommand{\SB}{\mathbb{S}}
\def\SS{\mathbb{S}}

\newcommand{\WW}{\mathbb{W}}

\newcommand{\Ba}{ {\cal B }}
\newcommand{\Ca}{ {\cal C }}
\newcommand{\Da}{ {\cal D }}
\newcommand{\La}{ {\cal L }}

\newcommand{\Ea}{ {\cal E }}
\newcommand{\Sa}{ {\cal S }}
\newcommand{\Ra}{ {\cal R }}
\newcommand{\Va}{ {\cal V }}
\newcommand{\Ua}{ {\cal U }}
\newcommand{\Fa}{ {\cal F }}
\newcommand{\Ga}{ {\cal G }}

\newcommand{\Oa}{ {\cal O }}
\newcommand{\Ia}{ {\cal I }}

\newcommand{\Ma}{ {\cal M }}
\newcommand{\Ta}{ {\cal T}}

\newcommand{\Ja}{ {\cal J }}
\newcommand{\Pa}{ {\cal P }}

\newcommand{\Wa}{ {\cal W }}

\numberwithin{equation}{section}

\begin{document}

\title{Perturbations and projections of Kalman-Bucy semigroups}

\author[$1$]{Adrian N. Bishop}
\author[$2$,$3$]{Pierre Del Moral}
\author[$3$]{Sahani D. Pathiraja}
\affil[$1$]{{\small University of Technology Sydney (UTS) and Data61 (CSIRO)}}
\affil[$2$]{{\small INRIA, Bordeaux Research Center, France}}
\affil[$3$]{{\small University of New South Wales (UNSW), Australia}}

\date{}
\maketitle

\begin{abstract}
We analyse various perturbations and projections of Kalman-Bucy semigroups and Riccati equations. For example, covariance inflation-type perturbations and localisation methods (projections) are common in the ensemble Kalman filtering literature. In the limit of these ensemble methods, the regularised sample covariance tends toward a solution of a perturbed/projected Riccati equation. With this motivation, results are given characterising the error between the nominal and regularised Riccati flows and Kalman-Bucy filtering distributions. New projection-type models are also discussed; e.g. Bose-Mesner projections. These regularisation models are also of interest on their own, and in, e.g., differential games, control of stochastic/jump processes, and robust control.
\end{abstract}

\section{Introduction}

The purpose of this work is to analyse a number of perturbations and projections of Kalman-Bucy \cite{kalman61,ap-2016} semigroups and of the associated (matrix differential) Riccati flow.

The results in this work are of interest in their own mathematical right. However, a prime motivating application for this work is the ensemble Kalman filter ({\tt EnKF}) \cite{evensen03} and the various ``regularisation'' methods used to ensure well-posedness of the sample covariance (e.g. sufficient rank) and to ``move'' the sample covariance closer (in some sense) to the Riccati flow of the true Kalman filter \cite{kalman61,ap-2016}. For example, two common forms of regularisation are covariance inflation-type methods (perturbations) and so-called covariance localisation methods (projections). Covariance inflation is a simple idea that involves adding some positive definite matrix to the sample covariance in order to increase its rank \cite{anderson99}; i.e. more specifically to account for an under-representation of the true variance due to a potentially inferior sample size. Separately, the idea of covariance localization involves multiplying (element-wise) the {\tt EnKF} sample covariance matrix via Schur (or Hadamard) products with certain sparse ``masking''  matrices with the intent of reducing spurious long-range correlations and increasing the sample covariance rank \cite{Houtekamer2001,Mitchell2002}. See \cite{hamill01} for an empirical examination of both types of regularisation. In these two cases, choosing the right inflation or localization is non-trivial and numerous ideas exist; e.g. \cite{Gaspari1999,Gaspari2006,anderson07,li2009,anderson09,sakov11,anderson2012}. Other related, and/or more subtle, regularisation methods exist and we will cover more general models in more detail in later sections; see also \cite{Heemink2001,pham2001,anderson2003,tippett03,sakov2008a,johns2008,Saetrom11,hou16} for related {\tt EnKF} methodology. 

Note that the total literature on {\tt EnKF} methodology is too broad to cover adequately here. Results on {\tt EnKF} convergence are recent (relative to this work) and concern, e.g., weak convergence with sample size \cite{legland09,mandel2011,law2016}, and stability \cite{kelly14,tong2016a,tong16b,delmoral16a,delmoral16b,2018arXiv180800235}, etc. The articles \cite{tong16b,Majda16} concern stability and robustness of the {\tt EnKF} in the presence of specific inflation and localisation methods. The articles \cite{apa-2017,2018arXiv180800235} study the behaviour of matrix-valued Riccati diffusions that capture the flow of the sample covariance various {\tt EnKF} implementations; i.e. its moment behaviour (non-asymptotic bias and variance), convergence, central-limit-type behaviour, stability etc. 

From a purely mathematical vantage, regularisation amounts to studying various projections and perturbations of the ``standard'' Riccati flow (viz \cite{kalman61,ap-2016}). The analytical behaviour of general projections and perturbations are a major focus of this study. We consider a broad class of perturbation model. We consider a particular projection model, and a certain class of localizable/diagonalizable systems adapted to these projections; the details are specified later. New ideas concerning projections relevant to the {\tt EnKF} are also introduced within this class. Given this analysis, we then study the (nonlinear) Kalman-Bucy diffusion \cite{ap-2016} and provide a number of contraction-type convergence results between the corresponding perturbed/projected diffusion and the optimal Kalman-Bucy diffusion. We study convergence in the mean-square sense and also in terms of the law of the diffusion. 

While methods in data assimilation and ensemble Kalman filtering are the main drivers of this work, the types of perturbations considered herein are more widely relevant: For example, our analysis captures well those perturbations of the ``standard'' Riccati flow that arise in, e.g., linear quadratic differential games \cite{Bernhard79,McAsey06,delfour07}, in the control of linear stochastic jump systems \cite{desouza90,abou-kandil94}, in certain robust and H$^\infty$ control settings \cite{freiling96,basar2008}, etc; see also the early work of Wonham \cite{wonham68} in linear-quadratic stochastic control. We also highlight the text \cite[e.g. Chap. 6]{abou-kandil03} and the references therein. Separately, a specific projected Riccati flow is studied in \cite{callier96}. Other relevant and related literature considering similar-type projections in estimation theory is given in \cite{Rebeschini2015,Majda1606.09087}. Going forward, we primarily rely on {\tt EnKF} motivators, but we emphasise here that the mathematical development is more broadly applicable.

Further introduction, discussion, and background is given in later subsections with a more technical focus. The organisation of this article is as follows:

{\small
\setcounter{tocdepth}{2}
\tableofcontents
}

\subsection{Kalman-Bucy diffusions}\label{desc-sec-intro}

The notation used throughout this article is introduced later in Section \ref{intro-notation}. However, the set-up in this section is relatively standard. Consider a time homogeneous linear-Gaussian filtering model of the following form
\begin{equation}\label{lin-Gaussian-diffusion-filtering}
\left\{
\begin{array}{rcl}
	dX_t&=&A\,X_t~dt~+~R^{1/2}\,dW_t\\
	dY_t&=&C\,X_t~dt~+~\Sigma^{1/2}\,dV_{t}
\end{array}
\right.
\end{equation}
where $(W_t,V_t)$ is an $(r+r^{\prime})$-dimensional standard Brownian motion, $X_0$ is a $r$-valued
Gaussian random vector (independent of $(W_t,V_t)$) with mean $\EE(X_0)$ and covariance matrix $P_0$, the symmetric positive definite matrices $R^{1/2}$ and $\Sigma^{1/2}$ are invertible, $A$ is an arbitrary square $(r\times r)$-matrix, $C$ is an arbitrary $(r^{\prime}\times r)$-matrix, and $Y_0=0$. Let $\Fa_t=\sigma\left(Y_s,~s\leq t\right)$ be the $\sigma$-algebra filtration generated by the observations.

It is well-known \cite{ap-2016} that the conditional distribution $\eta_t$ of the signal state $X_t$ given $\Fa_t$ is a $r$-dimensional Gaussian distribution with a mean and covariance matrix given by
$$
\widehat{X}_t:=\EE(X_t~|~\Fa_t)\quad\mbox{\rm and}\quad
P_t:=\EE\left(\left(X_t-\EE(X_t~|~\Fa_t)\right)\left(X_t-\EE(X_t~|~\Fa_t)\right)^{\prime}\right)
$$ 
given by the Kalman-Bucy and the Riccati equations 
\begin{equation}
d\widehat{X}_t=A~\widehat{X}_t~dt+P_{t}~C^{\prime}\Sigma^{-1}~\left(dY_t-C\widehat{X}_tdt\right)\quad\mbox{\rm with}\quad
\partial_tP_t=\ricc(P_t).\label{nonlinear-KB-Riccati}
\end{equation}
In the above display, $\ricc$ stands for the Riccati drift function from $\SS^+_{r}$ into $\SS_{r}$ defined for any $Q\in \SS^+_{r}$ by
\begin{equation}\label{def-ricc}
\ricc(Q)=AQ+QA^{\prime}-QSQ+R\quad\mbox{\rm with}\quad  S:=C^{\prime}\Sigma^{-1}C. 
\end{equation}
We now consider the conditional nonlinear McKean-Vlasov type diffusion process
\begin{equation}\label{Kalman-Bucy-filter-nonlinear-ref}
d\overline{X}_t=A~\overline{X}_t~dt~+~R^{1/2}~d\overline{W}_t+\Pa_{\eta_t}C^{\prime}\Sigma^{-1}~\left[dY_t-\left(C\overline{X}_tdt+\Sigma^{1/2}~d\overline{V}_{t}\right)\right]
\end{equation}
where  $(\overline{W}_t,\overline{V}_t,\overline{X}_0)$ are independent copies of $(W_t,V_t,X_0)$ (thus independent of
 the signal and the observation path). The notation $\Pa_{\eta_t}$ stands for the covariance matrix
\begin{equation}\label{def-nl-cov}
\Pa_{\eta_t}=\eta_t\left[(e-\eta_t(e))(e-\eta_t(e))^{\prime}\right]
\quad\mbox{\rm with}\quad \eta_t:=\mbox{\rm Law}(\overline{X}_t~|~\Fa_t)\quad\mbox{\rm and}\quad
e(x):=x.
\end{equation}
We shall call this probabilistic model (\ref{Kalman-Bucy-filter-nonlinear-ref}) the Kalman-Bucy (nonlinear) diffusion process. 

The ensemble Kalman-Bucy filter ({\tt EnKF}) coincides with the mean-field particle approximation of the nonlinear diffusion process \eqref{Kalman-Bucy-filter-nonlinear-ref}. To be more precise we let
 $(\overline{W}^i_t,\overline{V}^i_t,\xi^i_0)_{1\leq i\leq N}$  be 
 $N$ independent copies of $(\overline{W}_t,\overline{V}_t,\overline{X}_0)$. 
In this notation, a naive {\tt EnKF} is given by the Mckean-Vlasov type interacting diffusion process
\begin{equation}\label{fv1-3}
\left\{
\begin{array}{rcl}
d\xi^i_t&=&A~\xi^i_tdt+R^{1/2}d\overline{W}_t^i+p_tC^{\prime}\Sigma^{-1}\left[dY_t-\left(C\xi^i_tdt+\Sigma^{1/2}~d\overline{V}^i_{t}\right)\right]\\
i&=&1,\ldots,N
\end{array}
\right.
\end{equation}
with the rescaled particle covariance $p_t:=\left(1-N^{-1}\right)^{-1}\Pa_{\eta^{N}_t}$ and where the covariance matrix $\Pa_{\eta^{N}_t}$ is defined similarly to (\ref{def-nl-cov}) but in terms of the empirical measures $\eta^{N}_t:=N^{-1}\sum_{1\leq i\leq N}\delta_{\xi^i_t}$.

We define the following semigroup notation. 

\begin{defi}
We let $\theta_{s,t}(x)$ be the stochastic flow associated with the underlying signal process (\ref{lin-Gaussian-diffusion-filtering}). We let $\phi_{s,t}(Q)$ be the semigroup associated with the matrix Riccati equation in (\ref{nonlinear-KB-Riccati}) with (\ref{def-ricc}). And we let $\psi_{s,t}(x,Q)$ and $\overline{\psi}_{s,t}(x,Q)$ be the vector stochastic flows associated with the Kalman-Bucy filter and the nonlinear diffusion defined in (\ref{nonlinear-KB-Riccati}) and (\ref{Kalman-Bucy-filter-nonlinear-ref}), with $s\leq t$ and $(x,Q)\in\RR^r\times\SS_{r}^+$. 
\end{defi}

We also make the following standing assumption: Throughout this work we take the standard controllability and observability conditions as holding; see Section \ref{cont-observ-conds} for a statement of these conditions, and \cite{kalman61,Antsaklis,ap-2016} for a broader discussion and details on controllability and observability in control and filtering theory.

A key feature of any {\tt EnKF} method, is the sample-based estimation of the solution to the Riccati equation using a collection of interacting Kalman-Bucy filters. Contrary to conventional covariance estimates based on independent random samples, the {\tt EnKF} is based on interacting samples. These samples are sequentially updated by a noisy observation process through a gain matrix that itself depends on the sample covariance. The corresponding process is highly nonlinear (even when the true signal and observation model is linear). In high dimensions, the interacting particle estimation of the Riccati solution experiences the same difficulties as any conventional sample covariance estimator. For example:
\begin{itemize}
  \item The sample covariance $p_t$ is the sample mean of $N-1$ independent unit-rank matrices and has null eigenvalues when $N-1<r$. Thus, in some principal directions, the {\tt EnKF} is driven solely by the signal diffusion. With unstable signals, the {\tt EnKF} will exhibit divergence as it is not corrected by the innovation process. In this setting, one cannot design a stable particle sampler of the nonlinear diffusion (\ref{Kalman-Bucy-filter-nonlinear-ref}) without some kind of regularization. 
  \item The estimation of sparse high-dimensional covariance matrices using a small number of independent samples cannot readily be achieved without incorporating some information on the sparsity structure of the desired limit. Several regularization techniques have been developed in the statistics literature; see e.g. \cite{Haff,ledoit2,dahl,chaudhuri,bickel,elkaroui,lam-fan,bien,levina,khare,chen,arcolano}. One key common feature is to eliminate (typically long-range) noisy-type empirical correlations when its known that the limiting correlation is null or very small. 
\end{itemize}

\subsection{Perturbations and projections}

From a pure mathematical position, our model of perturbation or projection is motivated by methodology that replaces the sample covariance $p_t$ in (\ref{fv1-3}) by some matrix $\pi(p_t)$, where $\pi:\SS^+_{r}\mapsto\SS^+_{r}$ is some judiciously chosen mapping. These methods coincide with the mean field particle approximation of the nonlinear diffusion $\overline{X}_t^{\pi}$ defined by (\ref{Kalman-Bucy-filter-nonlinear-ref}) with $\Pa_{\eta_t}$ replaced by $\pi(\Pa_{\eta^{\pi}_t})$, i.e.,
\begin{equation}\label{Kalman-Bucy-filter-nonlinear-ref-reg}
d\overline{X}_t^{\pi}=A~\overline{X}_t^{\pi}~dt~+~R^{1/2}~d\overline{W}_t+\pi(\Pa_{\eta^{\pi}_t})~C^{\prime}\Sigma^{-1}~\left[dY_t-\left(C\overline{X}^{\pi}_tdt+\Sigma^{1/2}~d\overline{V}_{t}\right)\right]
\end{equation}
where $\eta^{\pi}_t=\mbox{\rm Law}(\overline{X}_t^{\pi}~\vert~\Fa_t)$. The initial state $\overline{X}_0^{\pi}$ is a Gaussian random variable with some covariance matrix $\Pa_{\eta^{\pi}_0}$. We expect the empirical average of the {\tt EnKF} system associated with (\ref{Kalman-Bucy-filter-nonlinear-ref-reg}) to converge to the Kalman-Bucy filter defined by (\ref{nonlinear-KB-Riccati}) except with $P_t$ replaced by the matrix $\pi(P_t)$. \textit{From the statistical viewpoint, the Kalman-Bucy filter {\small$\widehat{X}^{\pi}:= \EE(\overline{X}_t^{\pi}\,|\,\Fa_t)$} defined by (\ref{Kalman-Bucy-filter-nonlinear-ref-reg}) captures the limiting bias of the {\tt EnKF} empirical mean, introduced by some perturbation and/or projection operator $\pi$.} The nonlinear diffusion (\ref{Kalman-Bucy-filter-nonlinear-ref-reg}) is well posed, and the flow of covariance matrices $P^{\pi}_t=\Pa_{\eta^{\pi}_t}$ satisfies
\begin{eqnarray}\label{pi-Riccati-def}
\partial_tP^{\pi}_t&=&\ricc^{\pi}(P^{\pi}_t)\nonumber\\
&:=&\left[A-\pi(P^{\pi}_t)S\right]P^{\pi}_t+P^{\pi}_t\left[A-\pi(P^{\pi}_t)S\right]^{\prime}+R+\pi(P^{\pi}_t)S\pi(P^{\pi}_t)
\end{eqnarray}
when $\pi$ is chosen so that (\ref{pi-Riccati-def}) has a unique positive definite solution; a proof of this assertion is provided in the Appendix. This equation captures the covariance flow of the limiting perturbed/projected Kalman-Bucy filter $\widehat{X}^{\pi}:= \EE\left(\overline{X}_t^{\pi}\,|\,\Fa_t\right)$ associated with (\ref{Kalman-Bucy-filter-nonlinear-ref-reg}). {\em Consequently, (\ref{pi-Riccati-def}) captures the bias in the limiting {\tt EnKF} sample covariance as $N\rightarrow\infty$}. \textit{This perturbed or projected Riccati equation (\ref{pi-Riccati-def}) is the main object of study in this work}.

Note that we focus on the limiting object (\ref{pi-Riccati-def}). Our analysis holds if one replaces (\ref{fv1-3}), or the regularised limiting object (\ref{Kalman-Bucy-filter-nonlinear-ref-reg}), with regularised versions of the (unregularised) ``deterministic'' {\tt EnKF} in \cite{sakov2008a} (or those in \cite{Reich2013,Taghvaei2016ACC}). The deterministic {\tt EnKF} in \cite{sakov2008a} swaps {\small$d\overline{V}^i_{t}$} in (\ref{fv1-3}) with a deterministic adjustment factor. Regularised versions of other {\tt EnKF} variants may also be considered, if they lead to the same regularised limiting object of interest, i.e. (\ref{pi-Riccati-def}); e.g. any {\tt EnKF} ``flavour'' leaving $\eta^{\pi=id}_t$ unchanged is covered. Of course, when studying the perturbed Riccati equation (\ref{pi-Riccati-def}) alone, our analysis is not even limited to {\tt EnKF}-type motivation, as noted in the introduction (and again later).

We define the following semigroup notation. 

\begin{defi}
 Given some mapping $\pi$ from $\SS_r^+$ into itself, we let $\phi_{s,t}^{\pi}(Q)$, resp. $\psi_{s,t}^{\pi}(x,Q)$ and $\overline{\psi}^{\pi}_{s,t}(x,Q)$ be the semigroup, respectively the stochastic flows associated with the Riccati equation (\ref{pi-Riccati-def}), respectively the Kalman-Bucy filter and the Kalman-Bucy diffusion associated with the nonlinear model (\ref{Kalman-Bucy-filter-nonlinear-ref-reg}), with $s\leq t$ and $(x,Q)\in\RR^r\times\SS_{r}^+$.  
 \end{defi}

In the further development we shall distinguish and analyze the two different cases: 
\begin{equation}\label{two-classes}
1)\quad \pi= id +\Delta \quad\mbox{\rm with}\quad\Delta\approx 0 \qquad\mbox{\rm or}\qquad 2)\quad   \pi\circ\pi=\pi
\end{equation}
where $id$ stands for the identity mapping.

The first class of model can be thought of as a local perturbation mapping. These mappings are associated to some parameter that describe the level of perturbation. This model includes the variance inflation techniques discussed in Section~\ref{var-inflation-sec} and Stein-Shrinkage models presented in Section~\ref{stein-sec}, among others. 

The second class of model corresponds to projection-type mappings such as masked projections (or localization methods) discussed in Section~\ref{block-diag-sec} and projection mappings on Bose-Mesner algebras discussed in Section~\ref{Bose-Mesner-sec}.

Later in Section~\ref{mean-repulsion-sec} we consider mean-repulsion type perturbations, and we highlight how the main results presented in this work can be applied more broadly than implied by (\ref{two-classes}) alone. 

We also show later that the first class of model can actually capture most projections considered herein, or more general classes of test-type driving estimators; see the discussion in Section \ref{examples-and-applications}.

\subsubsection{Discussion: Perturbation-type regularization}

Consider the first class of perturbation model in (\ref{two-classes}). Under this model, several variance inflation methods have been proposed in the data assimilation literature as a simple means to address some of these numerical issues \cite{anderson99,hamill01,anderson07,li2009,anderson09}. By far the simplest technique is to add an artificial diagonal (positive definite) matrix to the sample covariance matrix $p_t$ in (\ref{fv1-3}). Another strategy is to consider a general class of Stein-Shrinkage-type perturbations models. These two strategies are discussed in Section~\ref{var-inflation-sec} and Section~\ref{stein-sec}. 

As an example, in view of (\ref{Kalman-Bucy-filter-nonlinear-ref-reg}), (\ref{pi-Riccati-def}), a simple variance inflation method $\pi(Q):=Q+ \Delta(Q)$, yields the following Riccati evolution
\begin{eqnarray}\label{pi-Riccati-var-inflation}
\partial_tP^{\pi}_t&=&\ricc^{\pi}(P^{\pi}_t)\nonumber\\
&:=&\left[A-\pi(P^{\pi}_t)S\right]P^{\pi}_t+P^{\pi}_t\left[A-\pi(P^{\pi}_t)S\right]^{\prime}+R+\pi(P^{\pi}_t)S\pi(P^{\pi}_t)\nonumber\\
&=&\ricc(P^{\pi}_t)+\Delta(Q)S\Delta(Q)
\end{eqnarray}
Obviously, such artificial inflations introduce an extra bias in the particle estimates delivered by the {\tt EnKF} (beyond the bias caused by a finite sample size and (nonlinear) interacting particles). In this example, a non-vanishing inflation term would generally be the sole cause of bias in the limiting {\tt EnKF} empirical mean and covariance as $N\rightarrow\infty$.

Later, we consider more general perturbation mappings that may arise in scenarios outside (ensemble) Kalman filtering such as in differential games, or in the control of linear stochastic jump systems, etc. These applications were briefly referenced in the introduction. These models will capture the preceding perturbation map (\ref{pi-Riccati-var-inflation}) as a special case.
 
Analysis of any bias-variance relationship trade-off requires one to quantify somewhat these two terms. This work focuses on the bias, in particular as it follows from the mapping $\pi$. For example, with the {\tt EnKF}, the $\LL_2$-error estimate at the origin with respect to the Frobenius norm is
$$
\EE\left[\Vert \pi(p_0)-P_0\Vert^2_F\right]~=~\Vert\pi(P_0)-P_0\Vert^2_F+\EE\left[\Vert \pi(p_0)-\pi(P_0)\Vert^2_F\right]
$$
whenever $\EE(p_0)=P_0$ and $\EE[\pi(p_0)]=\pi(P_0)$. Unfortunately, this unbiasedness property is not preserved in time $t>0$, due to the mean-field interactions; i.e. the {\tt EnKF} estimate $p_t$ of $P_t$ is biased in any case (e.g. even with $\pi=id$) due to the particle approximation/interaction. We don't study the bias arising from the mean field approximation here, and our analysis is mostly deterministic and focused on the relevant regularisation mappings. See \cite{apa-2017,2018arXiv180800235} for a detailed study of the bias (and variance, etc) of a stochastic matrix Riccati diffusion that captures the flow of the (finite $N$) sample covariance in a naive (non-regularised $\pi=id$) {\tt EnKF} implementation.

The general class of all perturbation-type mappings considered in this work is discussed in Section~\ref{perturbation-models-sec} and Section~\ref{first-second-sec-KF} (see also Sections~\ref{var-inflation-sec} and \ref{stein-sec}).

\subsubsection{Discussion: Projection-type regularization}

Consider now the second class of projection models in (\ref{two-classes}). Under the {\tt EnKF} framework, these projections are often defined in terms of the Hadamard product (a.k.a. Schur product) of the sample covariance matrix with some mask \cite{Gaspari1999,Houtekamer2001}. Here we may approximate such masks with a matrix $L$ of $\{0,1\}$-valued entries. The null entries capture the desired sparsity of the estimate. In the signal processing and data assimilation literature, these projections are often referred to as localization techniques. The study of $\{0,1\}$-valued mask matrices $L$ allows us to make rigorous convergence statements, and these results may act as a proxy for qualitatively understanding the behaviour in more general cases such as those considered in \cite{Gaspari1999,Houtekamer2001}. In the statistics literature, a random matrix given by the Hadamard product $L\odot p_0$ associated with some sample covariance $p_0$ is called a masked (or banded) sample covariance estimator of some limiting matrix $P_0$, see~\cite{bickel,elkaroui,levina,chen}. 

These projection techniques require the solution of the true unperturbed Riccati equation (the desired limit of the sample covariance) to lie within some class of (at least ``approximately'') ``band-able'' covariance matrices. To avoid the introduction of a huge bias \cite{Mitchell2002}, some prior knowledge of the sparsity/correlation structure of the solution to the Riccati equation is typically needed. However, the sparsity structure of a prescribed filtering problem is generally difficult to extract from the signal and sensor models etc. In some cases, the sparsity structure of the matrices $P_t$ can be estimated online from the particle model; e.g. see the {\tt Isomap} algorithm described in~\cite{Tenenbaum2000,wagaman}.

As with the first class of perturbation models, the choice of mapping $\pi$ under the second class of projection model introduces a deterministic bias. For example, in the filtering problem discussed in Section~\ref{block-diag-sec}, $P_0$ is a block-diagonal covariance matrix associated with $n$-independent filtering problems. In this case, we have $\pi(P_0)=L\odot P_0=P_0$ for some judicious block-diagonal matrix $L$ with $\{0,1\}$-valued entries. With this choice, it also follows that $L\odot P_t=P_t$. However, as noted before, the {\tt EnKF} derived (finite) sample covariance matrices are always (randomly) biased due to the (random) particle approximations/interactions, so that $L\odot p_t\not=p_t$ for any $t> 0$. Hence the effect of this projection in practice is to ``enforce'' some structure on the sample covariance at each time. In the limit $N\rightarrow\infty$ one hopes to recover the property $L\odot p_t\rightarrow L\odot P_t=P_t$.

In the general case, the fluctuations of $L\odot p_0$ around its limiting average value $L\odot P_0$ depend only on the non-zero entries. More precisely, for any symmetric mask-matrix $L$ with $\{0,1\}$-entries and at most $l$-zeros in each row we have the Levina-Vershynin's inequality,
$$
\EE\left[\Vert L\odot (p_0-P_0)\Vert_2\right]~\leq~ c\,\log^3{(2r)}\,\left[\tfrac{l}{N}+\sqrt{\tfrac{l}{N}}\right]~\Vert P_0\Vert_2
$$
for some finite universal constant $c<\infty$; see~\cite{levina,chen}. Of course, as before, this relationship is not so nicely preserved in time $t>0$ when comparing $p_t$ and $P_t$, due to the random particle approximation/interaction which introduces its own bias and fluctuations. Again we point to \cite{apa-2017,2018arXiv180800235} for a discussion on these random (particle) induced fluctuations. 

Another example class of projections discussed in detail in Section \ref{Bose-Mesner-sec} are orthogonal projections on Bose-Mesner-type cellular algebras w.r.t. the Frobenius norm \cite{bose}. These more sophisticated projections are more interesting than those examples in Section \ref{block-diag-sec} and can be used to project sample covariance matrices based on the topological/graph structure of the matrices $(A,R,S)$. 

The general class of all projection-type mappings considered in this work is discussed in Section~\ref{projection-models-sec} and Section~\ref{projection-models-sec-2}; see also Section~\ref{block-diag-sec} and~\ref{Bose-Mesner-sec} for those examples discussed above.

\subsection{Some background notation}\label{intro-notation}

This section details some basic notation and terms used throughout the article.

Let $\left\Vert\point\right\Vert_2$ be the Euclidean norm on $\RR^{r}$, $r\geq 1$. We denote by $\MM_r$ the set of $(r\times r)$-square matrices with real entries, $\SS_r\subset \MM_r$ the set of $(r\times r)$ real symmetric matrices, and by $\SS_r^+\subset\SS_r$ the subset of symmetric positive (semi)-definite matrices. With a slight abuse of notation, we denote by $Id$ the $(r\times r)$ standard identity matrix (with the size obvious from the context). Given some subsets $\Ia,\Ja\subset \{1,\ldots,r\}$ we set $A_{\Ia,\Ja}=\left(A_{i,j}\right)_{(i,j)\in( \Ia\times\Ja)}$ and $A_{\Ia}=A_{\Ia,\Ia}$. 

Denote by $\lambda_i(A)$, with $1\leq i\leq r$, the non-increasing sequence of eigenvalues of a $(r\times r)$-matrix $A$ and let $\mbox{\rm Spec}(A)$ be the set of all eigenvalues. We often denote by $\lambda_{min}(A)=\lambda_{r}(A)$ and $\lambda_{max}(A)=\lambda_{1}(A)$ the minimal and the maximal eigenvalue. We set $A_{sym}:=(A+A^{\prime})/2$ for any $(r\times r)$-square matrix $A$.
We define the logarithmic norm $\mu(A)$ of an $(r_1\times r_1)$-square matrix $A$ by
\begin{equation}\label{def-log-norm}
\begin{array}{rcl}
\mu(A)&:=&\inf{\{\alpha:\forall x,~\langle x,Ax\rangle\leq \alpha \left\Vert x\right\Vert_2^2\}} \\
&=&\lambda_{max}\left(A_{sym}\right)\\
&=&\inf{\{\alpha:\forall t\geq 0,~\Vert \exp{(At)}\Vert_2\leq \exp{(\alpha t)}\}}.
\end{array}
\end{equation}
The above equivalent formulations show that 
$$
\mu(A)\,\geq\, \varsigma(A):=\max{\left\{\mbox{\rm Re}(\lambda)~:~\lambda\in \mbox{\rm Spec}(A)\right\}}
$$
where $\mbox{\rm Re}(\lambda)$ stands for the real part of the eigenvalues $\lambda$.
The parameter $\varsigma(A)$ is often called the spectral abscissa of $A$. Also note that $A_{sym}$ is negative definite as soon as $\mu(A)<0$. The Frobenius matrix norm of a $(r_1\times r_2)$ matrix $A$ is defined by
$$
\left\Vert A\right\Vert_{F}^2=\mbox{\rm tr}(A^{\prime}A)
\qquad\mbox{\rm
with the trace operator $\mbox{\rm tr}(\point)$.}
$$
If $A$ is a matrix $(r\times r)$, we have $\left\Vert A\right\Vert_{F}^2=\sum_{1\leq i,j\leq r}A(i,j)^2$. For any $(r\times r)$-matrix $A$, we recall norm equivalence formulae
$$\Vert A\Vert_2^2=\lambda_{max}(A^{\prime}A)\leq \mbox{\rm tr}(A^{\prime}A)=\Vert A\Vert_F^2\leq r~\Vert A\Vert_2^2.$$
For any matrices $A$ and $B$ we also have the estimate
$$
\lambda_{min}(AA^{\prime})^{1/2}~~\Vert B\Vert_F\leq \Vert AB\Vert_F\leq \lambda_{max}(AA^{\prime})^{1/2}~\Vert B\Vert_F.
$$

We also quote a Lipschitz property of the square root function on (symmetric) positive definite matrices. For any $Q_1,Q_2\in\SS_{r}^+$
\begin{equation}\label{square-root-key-estimate}
\Vert Q_1^{1/2}- Q_2^{1/2}\Vert \leq \left[\lambda^{1/2}_{min}(Q_1)+\lambda^{1/2}_{min}(Q_2)\right]^{-1}~\Vert Q_1- Q_2\Vert
\end{equation}
for any unitary invariant matrix norm (such as the $\LL_2$-norm or the Frobenius norm). See for instance Theorem 6.2 on page 135 in~\cite{higham}, as well as Proposition 3.2 on page 591 in~\cite{hemmen}.

The Hadamard-Schur product of two $(r\times r^{\prime})$-matrices $A$ and $B$ of the same size is defined by the matrix $A\odot B$
with entries $(A\odot B)_{i_1,i_2}=A_{i_1,i_2}B_{i_1,i_2}$ for any $1\leq i_1\leq r$ and $1\leq i_2\leq r^{\prime}$. With a slight abuse of notation, we denote by $J$ the $(r\times r^{\prime})$ Hadamard-Schur identity matrix with all unit entries. By Theorem 17 in~\cite{johnson}, we recall that for any symmetric positive semi-definite matrices
$(A,B,P,Q)$ we have
\begin{equation}\label{ref-johnson-17}
P\geq Q\geq 0\quad\mbox{\rm and}\quad A\geq B\geq 0 \qquad \Longrightarrow \quad P\odot A\geq Q\odot B.
\end{equation}

Now, given some random variable $Z$ with some probability measure or distribution $\eta$ and some measurable function
$f$ on some product space $\RR^r$, we let $$\eta(f)=\EE(f(Z))=\int~f(x)~\eta(dx)$$ be the integral of $f$ w.r.t. $\eta$ or the expectation of $f(Z)$.  As a rule any multivariate variable, say $Z$, is represented by a column vector and we use the transposition operator $Z^{\prime}$ to denote the row vector (similarly for matrices; already seen above).

We also need to consider the $n$-th Wasserstein distance between two probability measures $\nu_1$ and $\nu_2$ on $\RR^r$ defined by
$$
	\WW_n(\nu_1,\nu_2)=\inf{\left\{\EE\left(\Vert Z_1-Z_2\Vert^n_2\right)^{\frac{1}{n}} \right\}}.
$$
The infimum in the above formula is taken over all pairs of random variable $(Z_1,Z_2)$ such that $\mbox{\rm Law}(Z_i)=\nu_i$, with $i=1,2$. We denote by $ \mbox{\rm Ent}\left(\nu_1~|~\nu_2\right)$ the Boltzmann-relative entropy
$$
 \mbox{\rm Ent}\left(\nu_1~|~\nu_2\right):=\int~\log\left(\frac{d\nu_1}{d\nu_2}\right)~d\nu_1\quad\mbox{if $\nu_1\ll\nu_2$, and $+\infty$ otherwise.}
$$

\subsection{Statement of the main results}\label{statement-sec}

In Section~\ref{perturbation-models-sec} and Section~\ref{projection-models-sec} (cf. Theorem~\ref{bplem} and formula (\ref{domination})) we will check that 
$$
	\phi^{\pi}_t(Q)\geq \phi_t(Q).
$$ 
This property shows that any $\pi$-perturbation or $\pi$-projection of the Kalman-Bucy diffusion induces a larger covariance matrix w.r.t. the Loewner order. This property is one key driving motivation for regularisation in the {\tt EnKF} literature.

Our first contribution concerns the continuity properties of the first class of perturbation models presented in (\ref{two-classes}) and introduced more formally in Section~\ref{perturbation-models-sec}. We consider a compact subset $\Pi$ of continuous mappings $\pi~:~\SS_r^+\mapsto~\SS_r^+$ equipped with the uniform norm induced by the $\LL_2$-norm on $\SS_{r}^+$. We let $B(\delta)$ be a $\delta$-ball around the identity mapping. For example, consider (\ref{Kalman-Bucy-filter-nonlinear-ref-reg}), (\ref{pi-Riccati-def}) and suppose further that
\begin{eqnarray}\label{pi-Riccati-def-gamma-term}
	\partial_tP^{\pi}_t&=&\ricc^{\pi}(P^{\pi}_t)\nonumber\\
		&:=&\left[A-\pi(P^{\pi}_t)S\right]P^{\pi}_t+P^{\pi}_t\left[A-\pi(P^{\pi}_t)S\right]^{\prime}+R+\pi(P^{\pi}_t)S\pi(P^{\pi}_t)\nonumber\\
		&=&\ricc(P^{\pi}_t)+\Gamma_{\pi}(P^{\pi}_t) 
\end{eqnarray}
with the quadratic positive mapping $\Gamma_{\pi}$ defined by
$$
	\Gamma_{\pi}(Q)=B_0+B_1Q+QB_1^{\prime}+QB_2Q+\Ra(Q)
$$
for some matrices $(B_0,B_1,B_2)\in \SS_{r}^3$ with $B_2\leq S$ and $\varpi:=\sup_{Q\in \SS^+_{r}}\Vert\Ra(Q)\Vert_2<\infty ~\Rightarrow~ \Ra(Q)\leq \varpi~Id.$ This model captures, e.g., simple inflation models like (\ref{pi-Riccati-var-inflation}), and Stein-Shrinkage methods like those discussed in Section~\ref{stein-sec}. This model also captures those perturbations relevant in, e.g., linear-quadratic differential games, control of stochastic jump processes, robust control theory, etc.

This mapping $\Gamma_{\pi}$ already hints that the analysis of the semigroups $\phi^{\pi}_t$ is a delicate mathematical problem, since it cannot be deduced directly from that of the Riccati flow $\phi_t$. By the Cauchy-Lipschitz theorem, the existence and the uniqueness of the flow of matrices $\phi^{\pi}_t(Q)$ {\em for any starting covariance matrix $Q$} is ensured by the local Lipschitz property of the drift function $\ricc^{\pi}$, on some open interval that may depend on $Q$. The existence of global solutions on the real line is not ensured as the quadratic term may induce a blow up on some finite time horizon.

In this setting, our first main result concerns the first class of perturbation models presented in (\ref{two-classes}), and takes the following {\em{mildly informal}} form.

 \begin{theomain}\label{theo-main-1}
 Assume that the filtering problem is observable and controllable. In this situation, under some regularity conditions,
there exists some $\delta>0$ such that for any $\epsilon<\delta$, any $\pi\in B(\epsilon)$, and any $n\geq 1$ we have
the uniform estimates
\begin{equation}\label{theo-main-1-estimates}
\sup_{t\geq 0}{\Vert\phi^{\pi}_{t}(Q)-\phi_{t}(Q)\Vert_2}\leq c(\delta)~\epsilon\quad\mbox{and}\quad
\sup_{t\geq 0}{\EE\left[\Vert \psi_{0,t}^{\pi}(x,Q)-\psi_{0,t}(x,Q) \Vert_2^{2n}\right]^{\frac{1}{2n}}}\leq c(\delta)\sqrt{n}~\epsilon
\end{equation}
for some finite constant $c(\delta)$ whose values only depend on the parameter $\delta$. 
 \end{theomain}
 
A {\em{precise statement of this result}} is given in Theorem~\ref{theo-continuity} and Theorem~\ref{theo-KB-stoch-flows}; e.g. with clarification of the required regularity conditions. The proof of the Riccati estimates in the l.h.s. of (\ref{theo-main-1-estimates}) is provided in Section~\ref{robustness-sec}, dedicated to the boundedness and the robustness properties of Riccati semigroups (cf. Theorem~\ref{theo-continuity}; see also \cite{ap-2018-franklin} for further discussion on these robustness and related results). The proof of the r.h.s. estimates in (\ref{theo-main-1-estimates}) is provided in Section~\ref{first-second-sec-KF} dedicated to the continuity properties of Kalman-Bucy stochastic flows (cf. Theorem~\ref{theo-KB-stoch-flows}).
 
 The preceding theorem concerns time-uniform bounds on the mean and the covariance of the Kalman-Bucy flows. Our second objective, given the first class of perturbations, is to quantify the difference between the conditional distributions of the nonlinear Kalman-Bucy diffusion,
$$
	\eta_{s,t}(x,Q):=\mbox{\rm Law}\left(\overline{\psi}_{s,t}(x,Q)~|~\Fa_{s,t}\right) \quad\mbox{\rm and}\quad \eta_{s,t}^{\pi}(x,Q):=\mbox{\rm Law}\left(\overline{\psi}_{s,t}^{\,\pi}(x,Q)~|~\Fa_{s,t}\right)
$$
where $\Fa_{s,t}=\sigma(Y_u, s\leq u\leq t)$ denotes the $\sigma$-field generated by the observations from time $s$ to the 
time horizon $t$. By construction $\overline{\psi}_{s,t}$ and $\overline{\psi}_{s,t}^{\,\pi}$ are time-varying Ornstein-Uhlenbeck-type diffusions (linear stochastic differential equations) \cite{ap-2016} and consequently $\eta_{s,t}(x,Q)$ and $\eta_{s,t}^{\pi}(x,Q)$ are both Gaussian distributions. Our next main result {\em{informally}} takes the following form.

\begin{theomain}\label{theo-main-2}
Under the assumptions of Theorem~\ref{theo-main-1}, for any $n\geq 1$, we have the almost sure relative entropy and Wasserstein distance estimates
\begin{eqnarray*}
 \mbox{\rm Ent}\left(\eta_{s,t}^{\pi}(x,Q)~|~\eta_{s,t}(x,Q)\right)
 &\leq&  \displaystyle c
 \left[\left\Vert{\psi}^{\pi}_{s,t}(x,Q)-{\psi}_{s,t}(x,Q)\right\Vert^2_2+\Vert \phi_{s,t}(Q)-\phi^{\pi}_{s,t}(Q)\Vert_2\right]\\
\WW_{2n}\left[\eta_{s,t}^{\pi}(x,Q),\eta_{s,t}(x,Q)\right]&\leq &\Vert{\psi}^{\pi}_{s,t}(x,Q)-{\psi}_{s,t}(x,Q)\Vert_2+
c ~\sqrt{n}~\Vert 
\phi_{s,t}^{\pi}(Q)-\phi_{s,t}(Q)
\Vert_2
\end{eqnarray*}
for some constant $c<\infty$ that depends on the system and observation matrices.
\end{theomain}
The proof of these estimates, with a more precise description of the constant $c$, is provided in Section~\ref{first-second-sec-KF}; e.g. see {\em{the precise statement of these results}} in Theorem~\ref{theo-entropy} and Theorem~\ref{theo-wass}.

The impact of these two theorems is illustrated in Section~\ref{var-inflation-sec} and Section~\ref{stein-sec} in terms of the variance inflation and the Stein-Shrinkage methods common in the data assimilation literature. 

Our second contribution concerns the continuity properties of the second class of projection mappings presented in (\ref{two-classes}) and discussed further in Section~\ref{projection-models-sec}. We assume that $\pi$ is some positive map from $\MM_r$ into itself, of the form
$$
	\pi(Q)=\argmin_{B\in \Ba}{\pi\left[(Q-B)(Q-B)^{\prime}\right]}\quad\mbox{\rm for some matrix ring  $\Ba\subset\MM_r$}.
$$
From the geometrical viewpoint, these orthogonal projections map the set $\SS_r^+$ into the set of matrices with the same sparsity structure as the matrices in the ring $\Ba$. These projection techniques are unbiased when the covariance graph of the filtering model (reflecting the sparsity structure of the matrices $P_t$) is defined in terms of the same association scheme. Thus, the {\em optimal} use of these projections requires some prior knowledge on the sparsity structure of the solution to the Riccati equation. This is a special class of projection model differing somewhat from the typical localization used in the {\tt EnKF} literature; e.g. see \cite{Gaspari1999,Houtekamer2001}. However, under the particular chosen class of projection, explicit and rigorous convergence results are possible when the correlation structure is well-enough adapted to the projection. Heuristically, these results may act as a proxy to gain intuitive, or qualitative, insight into the behaviour of more practical localization implementations \cite{Houtekamer2001}; e.g. and can be taken in combination with the first class of perturbation model for this purpose.

A prototype model satisfying these conditions are orthogonal projections onto the set of block-diagonal matrices $\Ba=\Ma_{r[1]}\oplus\ldots\oplus\Ma_{r[n]}\subset \MM_r$, with $r=\sum_{1\leq q\leq n}r[q]$. Another important class of models satisfying the above conditions are orthogonal projections on Bose-Mesner-type cellular algebras w.r.t. the Frobenius norm \cite{bose}. These more sophisticated projections are interesting and can be used to project sample covariance matrices based on the topological/graph structure of the matrices $(A,R,S)$. 

See Section~\ref{block-diag-sec} for applications to block-diagonal masking matrices and Section~\ref{Bose-Mesner-sec} for further discussion on Bose-Mesner projections; e.g. Section~\ref{riccati-solver} provides an explicit solution of the Riccati equation as soon as the matrices $(A,R,S)$ and the initial condition belong to some Bose-Mesner algebra.

In this context, our third main result takes the following {\em{mildly informal}} form.

\begin{theomain}\label{theo-projection-intro}
 Assume that the filtering problem is observable and controllable and assume that $(A,A^{\prime},S,R)\in \Ba$. In this situation we have
\begin{equation}\label{theo-main-2-estimates}
\phi_t^{\pi}\circ\pi=\phi_t\circ\pi \quad\mbox{and}\quad \psi_{s,t}^{\pi}(x,Q)= 
\psi_{s,t}(x,\pi(Q))
\end{equation}
for any $(x,Q)\in (\RR^r\times\SS_r^+)$ and $t\geq 0$.
In addition, there exists some $\rho>0$ such that for any $Q\in \SS_r^+$ and any time horizon 
$t\geq 0$ we have the local exponential-Lipschitz  inequality
\begin{equation}\label{theo-main-2-estimates-2}
\Vert \phi_t^{\pi}(Q)-\phi_t(Q)\Vert_2\leq c_Q~e^{-\rho t}~\Vert Q-\pi(Q)\Vert_2
\end{equation}
for some finite constant $c_Q$ whose values only depend on $\Vert Q\Vert_2$.
\end{theomain}

The result in (\ref{theo-main-2-estimates}) is {\em{stated precisely}} as Theorem \ref{theo-pi-1} and is covered also in Section~\ref{projection-models-sec-2}. The estimate (\ref{theo-main-2-estimates-2}) is {\em{stated precisely}} in Theorem \ref{exp-contraction-proj}; see also the corollaries in Section~\ref{expo-concentration}.

The relationship (\ref{theo-main-2-estimates}) shows that the set $\Ba$ is stable w.r.t. the $\pi$-projected Riccati flow. The exponential estimate (\ref{theo-main-2-estimates-2}) shows that, for any initial condition, the Kalman-Bucy stochastic flow as well as the $\pi$-projected Riccati flow converges to the set $\Ba$ as the time horizon $t$ tends to $\infty$.

Last, but not least, Theorem~\ref{theo-projection-intro} allows one to transfer, without further work, all the exponential contraction inequalities developed in~\cite{ap-2016}, dedicated to the stability properties of Kalman-Bucy diffusions.

\subsection{Some background results}

\subsubsection{Observability, controllability and the steady-state Riccati equation} \label{cont-observ-conds}

We  assume that $(A,R^{1/2})$ is a controllable pair and $(A,C)$ is observable in the sense that
\begin{equation}\label{def-contr-obs}
\left[R^{1/2},A(R^{1/2})\ldots, A^{r-1}R^{1/2}\right]\quad
\mbox{\rm and}\quad
\left[\begin{array}{c}
C\\
CA\\
\vdots\\
CA^{r-1}
\end{array}
\right]
\end{equation}
have rank $r$. With $R$ positive definite as assumed here (and common in filtering problems), the controllability condition follows immediately. We consider the observability and controllability Gramians $(\Oa_{t},\Ca_{t}(\Oa))$ and $(\Ca_{t},\Oa_{t}(\Ca))$ 
associated with the triplet $(A,R,S)$ and defined by
\begin{eqnarray*}
\Oa_{t}:= \int_{0}^{t}~e^{-A^{\prime}s}~S~e^{-As}~ds \qquad\mbox{and}\qquad \Ca_{t}(\Oa)&:=& \Oa_{t}^{-1}\left[\int_0^t~e^{-(t-s)A^{\prime}}~\Oa_{s}~R~\Oa_{s}~e^{-(t-s)A}~ds\right]\Oa_{t}^{-1}
\\
\Ca_{t} := \int_{0}^{t}~e^{As}~R~e^{A^{\prime}s}~ds \qquad\mbox{and}\qquad \Oa_{t}(\Ca)&:=&\Ca_{t}^{-1}\left[\int_0^t~e^{(t-s)A}~\Ca_{s}~ S~\Ca_{s}~e^{(t-s)A^{\prime}}~ds\right]\Ca_{t}^{-1}.
\end{eqnarray*}
Given the rank assumptions on (\ref{def-contr-obs}), there exists some parameters $\upsilon,\varpi^{o,c}_{\pm},\varpi^{c}_{\pm}(\Oa),\varpi^{o}_{\pm}(\Ca)>0$ such that
\begin{equation}\label{steady-state-eq-2}
\varpi_-^{c}~Id\leq \Ca_{\upsilon}
\leq \varpi_+^{c}~Id
\quad\mbox{\rm
and}
\quad
\varpi_-^{o}~Id\leq \Oa_{\upsilon}\leq \varpi_+^{o}~Id
\end{equation}
as well as
$$
 \varpi_-^{c}(\Oa)~Id
\leq \Ca_{\upsilon}(\Oa)\leq \varpi_+^{c}(\Oa)~Id\quad
\mbox{\rm and}\quad
 \varpi_-^{o}(\Ca)~Id
\leq \Oa_{\upsilon}(\Ca)\leq \varpi_+^{o}(\Ca)~Id.
$$
The parameter $\upsilon$ is often called the interval of observability-controllability. By Theorem 4.4 in~\cite{ap-2016}, for any $t\geq \upsilon$ and any $Q\in \SB_{r}^+$ we have the uniform estimates
\begin{equation}\label{after-upsilon}
\left(\Oa_{\upsilon}(\Ca)+ \Ca_{\upsilon}^{-1}\right)^{-1}~\leq \phi_{t}(Q)~\leq 
\Oa_{\upsilon}^{-1}+\Ca_{\upsilon}(\Oa).
\end{equation}
When (\ref{steady-state-eq-2}) is satisfied, we say that a triplet $(A,R,S)$ satisfy the Gramian condition for some parameters $\upsilon,\varpi^{o,c}_{\pm}>0$. These conditions ensure the existence and the  uniqueness of a positive definite fixed-point matrix $P$ solving the so-called algebraic Riccati equation
\begin{equation}\label{steady-state-eq}
\ricc(P):=AP+PA^{\prime}-PSP+R=0.
\end{equation}
Importantly, in this case, the matrix difference $A-PS$ is asymptotically stable (Hurwitz stable) even when the signal matrix $A$ is unstable \cite[Theorems 9.12, 9.15]{Lancaster1995}. More relaxed conditions (i.e. detectability and stabilisability) for a stabilising solution (perhaps only positive semi-definite) to exist are discussed widely in the literature; see \cite{kucera72,Molinari77,Lancaster1995} and the convergence results in \cite{Kwakernaak72,callier81}.

\subsubsection{Exponential and Kalman-Bucy semigroup estimates}

The transition matrix associated with a smooth flow of $(r\times r)$-matrices $A:u\mapsto A_u$ is denoted by
$$
\Ea_{s,t}(A)=\exp{\left[\oint_s^t A_u~du\right]}\Longleftrightarrow \partial_t \Ea_{s,t}(A)=A_t~\Ea_{s,t}(A)\quad\mbox{\rm and}\quad
\partial_s \Ea_{s,t}(A)=-\Ea_{s,t}(A)~A_s
$$
for any $s\leq t$, with $\Ea_{s,s}=Id$, the identity matrix. Equivalently in terms of the fundamental solution matrices
$\Ea_t(A):=\Ea_{0,t}(A)$ we have
$
\Ea_{s,t}(A)=\Ea_t(A)\Ea_s(A)^{-1}
$.

The following technical lemma provides a pair of semigroup estimates  of the state transition matrices associated with a sum of drift-type matrices. 
\begin{lem}[\cite{ap-2016}]\label{perturbation-lemma-intro}
Let $A~:~u\mapsto A_u$ and $B~:~u\mapsto B_u$ be some smooth flows of $(r\times r)$-matrices. For any $s\leq t$  and any matrix norm $\Vert\cdot\Vert$ we have
$$
	\Vert  \Ea_{s,t}(A) \Vert\leq \alpha_A~\exp{\left(-\omega_A~(t-s)\right)}\Rightarrow
 		\left\Vert   \Ea_{s,t}(A+B)\right\Vert\leq \alpha_A \exp{\left[-\omega_A(t-s)+\alpha_A\int_s^t \Vert B_u\Vert~du\right]}.
$$
for some positive constant $\alpha_A$ and some parameter $\omega_A$.
\end{lem}

For any $s\leq t$ and $Q\in\SS_{r}^+$ we set
$$
E_{s,t}(Q):=\exp{\left[\oint_s^t\left(A-\phi_u(Q)S\right)~du\right]}.
$$
When $s=0$ sometimes we write $E_{t}(Q)$ instead of $E_{0,t}(Q)$. In this notation we have
$$
E_{s,t}(Q)=E_{t}(Q)E_{s}(Q)^{-1}.
$$
For any $s\leq u\leq t$ and $Q\in\SS_{r}^+$ we set
$$
	E_{t\vert s}(Q)=\exp{\left[\oint_s^t\left(A-\phi_{s,v}(Q)~S\right)~dv\right]} \quad\mbox{and}\quad E_{u,t\vert s}(Q):=E_{t\vert s}(Q)E_{u\vert s}(Q)^{-1}.
$$
Also observe that
\begin{eqnarray*}
	E_{s,t}(Q)&=&\exp{\left[\oint_{s}^{t}\left(A-\phi_{s,u}(\phi_s(Q))S\right)du\right]}=E_{t\vert s}(\phi_s(Q)).
\end{eqnarray*}
For any $s\leq u\leq t$ and any $Q\in\SS_{r}^+$ we have
\begin{equation}\label{ref-E-s-t}
	E_{t\vert s}(Q)=E_{t-s}(Q)\quad\mbox{and}\quad E_{u,t\vert s}(Q)=E_{(u-s),(t-s)}(Q).
\end{equation}

Observe that the Riccati equation is time-homogeneous so that 
$$
	\phi_{s,s+t}(Q)=\phi_{t}(Q):=\phi_{0,t}(Q).
$$
By Proposition 4.3 in~\cite{ap-2016} we have
\begin{equation}\label{upper-bound-Phi}
	0\leq \phi_{t}(Q)\leq P+e^{(A-PS)t}(Q-P)e^{(A-PS)t}\Longrightarrow \Vert \phi_{t}(Q)\Vert_2\leq \Vert P\Vert_2+\kappa \Vert Q-P\Vert_2
\end{equation}
for some constant $\kappa$ whose values doesn't depend on the time parameter nor on $Q$. We also have the following contraction result.

\begin{theo}[\cite{ap-2016}]
For any $Q_1,Q_2\in\SS^{+}_{r}$ and for any $t\geq 0$ we have the local contraction inequality
\begin{eqnarray}
\Vert E_t(Q_1)\Vert_2&\leq& \kappa_E(\Vert Q_1\Vert_2)~e^{-2\nu t}\label{expo-E}\\
\Vert \phi_t(Q_2)-\phi_t(Q_1)\Vert_2&\leq &\kappa_{\phi}(\Vert Q_1\Vert_2,\Vert Q_2\Vert_2)~e^{-2\nu t}~\Vert Q_2-Q_1\Vert_2\label{expo-phi}
\\
\Vert E_t(Q_2)-E_t(Q_1)\Vert_2&\leq &\kappa_{E}(\Vert Q_1\Vert_2,\Vert Q_2\Vert_2)~e^{-\nu t}~\Vert Q_2-Q_1\Vert_2\label{expo-Elip}
\end{eqnarray}
for some rate $\nu>0$, and some finite non-decreasing functions $\kappa_E(q_1),\kappa_{E}(q_1,q_2)$, $\kappa_{\phi}(q_1,q_2)<\infty$.
\end{theo}

\section{Riccati semigroups}

\subsection{Variational and backward semigroups}

We let $\La(\SS_r,\SS_{r})$ be the set of bounded linear functional from $\SS_r$ into itself, and equipped with the Frobenius norm. A mapping $\phi:\SS^+_{r}\mapsto \SS^+_r$ is Fr\'echet differentiable at some $Q_1\in \SS_r^+$ if there exists a continuous linear functional
$
\partial \phi(Q_1)\in  \La(\SS_r,\SS_{r})
$
such that
$$
\lim_{Q_2\rightarrow Q_1}\Vert Q_2-Q_1\Vert_F^{-1}
\Vert\phi(Q_2)-\phi(Q_1)-\partial \phi(Q_1)\cdot (Q_2-Q_1) \Vert_F=0.
$$

For instance the first-order Frechet-derivative of the Riccati quadratic drift function $$\ricc~:~Q \in\SS_{r}^+\mapsto \ricc(Q)\in \SS_{r}$$
defined in (\ref{def-ricc}) is given for any $ (Q_1,Q_2)\in (\SS_{r}^+\times \SS_r)$ by the formula
\begin{equation}\label{partial-ricc}
\partial \mbox{\rm Ricc}(Q_1)\cdot Q_2 =(A-Q_1S) Q_2+Q_2(A-Q_1S) ^{\prime}.
\end{equation}

\begin{lem}\label{lem-frechet-derivatives}
For any $t\geq 0$   the mapping  $Q\mapsto  \phi_t(Q)$ is Fr\'echet  differentiable and
 for any $(Q_1,Q_2)\in(\SS_{r}^+\times \SS^+_{r})$ we have
the formulae
\begin{eqnarray*}
\partial \phi_t(Q_1)\cdot Q_2&=&
   E_t(Q_1)~Q_2~   E_t(Q_1)^{\prime}.
\end{eqnarray*}

\end{lem}

\proof
Using the decomposition
\begin{eqnarray*}
\phi_t(Q_1)-\phi_t(Q_2)&
\displaystyle=&E_{s,t}(Q_2)\left[\phi_s(Q_1)-\phi_s(Q_2)\right]~E_{s,t}(Q_2)^{\prime}\\
&&\hskip-.5cm\displaystyle-\int_s^t~
E_{u,t}(Q_2)~\left[\phi_u(Q_1)-\phi_u(Q_2)\right]~S~\left[\phi_u(Q_1)-\phi_u(Q_2)\right]~
E_{u,t}(Q_2)^{\prime}~du
\end{eqnarray*}
we have
\begin{eqnarray*}
\phi_t(Q_2)-\phi_t(Q_1)&
\displaystyle=&E_{t}(Q_1)\left[Q_2-Q_1\right]~E_{t}(Q_1)^{\prime}\\
&&\hskip-.5cm\displaystyle-\int_0^t~
E_{u,t}(Q_1)~\left[\phi_u(Q_2)-\phi_u(Q_1)\right]~S~\left[\phi_u(Q_2)-\phi_u(Q_1)\right]~
E_{u,t}(Q_1)^{\prime}~du.
\end{eqnarray*} 
We end the proof of the first assertion using the Lipschitz property (\ref{expo-phi}). The proof of the  lemma is completed. \qed

We have the following backward flow and first-order variational result that will be used subsequently, but which is also of interest in its own right.

\begin{prop}\label{prop-key}
For any $Q\in\SS_{r}^+$ and any $0\leq s\leq t$ we have
$$
\partial_{s}\phi_{s,t}(Q)=-\ricc(\phi_{s,t}(Q))
\quad\mbox{and}\quad
\partial_{t}\phi_{s,t}(Q)=\ricc(\phi_{s,t}(Q))=\partial \phi_{s,t}(Q)\cdot \ricc(Q).
$$
In addition, the first-order variational equation associated with the Riccati equation is given by
the composition formula
\begin{equation}\label{first-order-var-riccati}
\partial_t\left(\partial \phi_t(Q)\right)=\partial\mbox{\rm Ricc}(\phi_t(Q))\circ \partial \phi_t(Q).
\end{equation}
\end{prop}
\proof

For any $Q\in\SS_{r}^+$ we have
$$
\partial_{s}\phi_{s,t}(Q)=\partial_{s}\phi_{0,t-s}(Q)=-\ricc(\phi_{0,t-s}(Q))=-\ricc(\phi_{s,t}(Q)).
$$
On the other hand, we have
$$
\begin{array}{l}
\displaystyle\Vert\ricc(\phi_{s-h,u}(Q))-\ricc(Q)\Vert_F\leq c_Q~h\\
\\
\Longrightarrow\displaystyle\Vert\int_{s-h}^{s}\left[\ricc(\phi_{s-h,u}(Q))-\ricc(Q)\right]~du\Vert_F\leq c_Q~h^2\\
\end{array}
$$
for some finite constant $c_Q$ whose values only depends on $\Vert Q\Vert_F$.
Using Lemma~\ref{lem-frechet-derivatives} this yields
\begin{eqnarray*}
\phi_{s-h,t}(Q)-\phi_{s,t}(Q)&=&\phi_{s,t}(\phi_{s-h,s}(Q))-\phi_{s,t}(Q)\\
&=&\phi_{s,t}\left(Q+\int_{s-h}^{s}\ricc(\phi_{s-h,u}(Q))du\right)-\phi_{s,t}(Q)\\
&=& \partial \phi_{s,t}(Q)\cdot \left[\int_{s-h}^{s}\ricc(\phi_{s-h,u}(Q))du\right]+\mbox{\rm o}(h)\\
&=& \partial \phi_{s,t}(Q)\cdot \ricc(Q)~h\\
&&\hskip2cm+\partial \phi_{s,t}(Q)\cdot \left[\int_{s-h}^{s}\left[\ricc(\phi_{s-h,u}(Q))-\ricc(Q)\right]du\right]+\mbox{\rm o}(h)\\
&=& \partial \phi_{s,t}(Q)\cdot \ricc(Q)~h+\mbox{\rm o}(h).\
\end{eqnarray*}
This implies that
$$
\partial_s\phi_{s,t}(Q)=\lim_{h\rightarrow 0}\frac{1}{-h}\left[\phi_{s-h,t}(Q)-\phi_{s,t}(Q)\right]=-\partial \phi_{s,t}(Q)\cdot \ricc(Q)
$$
from which we conclude that
\begin{eqnarray}\label{la-ref}
-\ricc(\phi_{s,t}(Q))+\partial \phi_{s,t}(Q)\cdot \ricc(Q) 
&=&\partial_s\phi_{s,t}(Q)+\partial_t\phi_{s,t}(Q)=0.
\end{eqnarray}
Finally, by Lemma~\ref{lem-frechet-derivatives} and (\ref{partial-ricc}) we have
\begin{eqnarray*}
\partial_t\left[\partial \phi_t(Q_1)\cdot Q_2\right]&=&
   \left[A-\phi_t(Q_1)S\right]~\left[\partial \phi_t(Q_1)\cdot Q_2\right]+\left[\partial \phi_t(Q_1)\cdot Q_2\right]~\left[A-\phi_t(Q_1)S\right]^{\prime}\\
   &=& \partial \mbox{\rm Ricc}(\phi_t(Q_1))\cdot \left[\partial \phi_t(Q_1)\cdot Q_2\right] \\
   &=&
   \left[\partial \mbox{\rm Ricc}(\phi_t(Q_1))\circ \partial \phi_t(Q_1)\right](Q_2).
\end{eqnarray*}
This ends the proof of the proposition.
\qed

\subsection{Perturbation-type models}\label{perturbation-models-sec}

\subsubsection{First and second order perturbations}\label{second-order-appr-sec}

We consider perturbation-type distortions in (\ref{pi-Riccati-def}) of the first type in (\ref{two-classes}). Formally, consider (\ref{pi-Riccati-def}) and (\ref{pi-Riccati-def-gamma-term}) and the class of perturbation mappings $\Gamma_{\pi}$ in (\ref{pi-Riccati-def-gamma-term}) with the hypothesis
$$
	\hskip-3cm\mbox{\rm (H)}_{0}\hskip3cm \Gamma_{\pi}(Q)=B_0+B_1Q+QB_1^{\prime}+QB_2Q+\Ra(Q)
$$
for some given matrices $(B_0,B_1,B_2)\in \SS_{r}^3$ such that $B_2\leq S$, and a uniformly bounded (symmetric) remainder term 
$$
	\varpi:=\sup_{Q\in \SS^+_{r}}\Vert\Ra(Q)\Vert_2<\infty ~\Longrightarrow~ \Ra(Q)\leq \varpi~Id.
$$

In this situation, the $\pi$-Riccati drift function $\ricc^{\pi}$ in (\ref{pi-Riccati-def}) takes the form
$$
\ricc^{\pi}(Q)~=~\ricc_{\pi}(Q)+\Ra_{\pi}(Q)~\leq~ \ricc_{\pi}(Q)
$$
with
\begin{equation}\label{def-Api}
\ricc_{\pi}(Q)
:=A_{\pi}Q+QA_{\pi}^{\prime}+R_{\pi}-QS_{\pi}Q,\qquad \Ra_{\pi}(Q)=\Ra(Q)-\varpi Id\leq 0
\end{equation}
and the matrices
$$
R_{\pi}:=R+B_0+\varpi Id\qquad  A_{\pi}:=A+B_1\quad\mbox{\rm and}\quad S_{\pi}:=S-B_2\geq 0.
$$

\begin{defi}
We let $\phi_{\pi,t}$, resp. $\phi^{\pi}_t$ be the Riccati flows associated with the drift function
$\ricc_{\pi}$ and resp. $\ricc^{\pi}$. We consider the observability and the controllability Gramians $(\Oa_{\pi,t},\Ca_{\pi,t}(\Oa))$ and $(\Ca_{\pi,t},\Oa_{\pi,t}(\Ca))$ associated with the triplet $(A_{\pi},R_{\pi},S_{\pi})$.
\end{defi}

We also let $ \Xi_{\pi}$ be the mapping from $\SS_{r}$ into itself defined by
$$
\Xi_{\pi}(Q):=\ricc_{\pi}(Q)-\ricc(Q)=
(A_{\pi}-A)Q+Q(A_{\pi}-A)^{\prime}+(R_{\pi}-R)-Q(S_{\pi}-S)Q.
$$
We also set
$$
\gamma(\pi):=\Vert A_{\pi}-A\Vert_2+\Vert R_{\pi}-R\Vert_2+\Vert S_{\pi}-S\Vert_2.
$$

We consider the following condition,
$$
\mbox{\rm (H)}_{1}\hskip1cm
\mbox{\rm $(A_{\pi},R_{\pi},S_{\pi})$ satisfies the Gramian condition (\ref{steady-state-eq-2}) for some
 $\upsilon_{\pi},\varpi^{o,c}_{\pm}(\pi)>0$}.
$$

We recall that this condition ensures the existence and the  uniqueness of
a positive definite fixed-point matrix $P_{\pi}$ solving the so-called algebraic Riccati equation
\begin{equation}\label{steady-state-eq-pi}
\ricc_{\pi}(P_{\pi}):=A_{\pi}P_{\pi}+P_{\pi}A_{\pi}^{\prime}-P_{\pi}S_{\pi}P_{\pi}+R_{\pi}=0.
\end{equation}
In addition, the matrix difference $A_{\pi}-P_{\pi}S_{\pi}$ is asymptotically stable.

Our first objective is to analyze the existence and the uniqueness of the flow $\phi^{\pi}_t$.

 \begin{theo}\label{bplem} 
Assume $\mbox{\rm (H)}_{0}$ and $\mbox{\rm (H)}_{1}$. For any $t\geq (\upsilon\vee \upsilon_{\pi})$ and $Q\in \SS^+_{r}$,
\begin{equation}\label{backward-analysis-cons} 
 \left(\Oa_{\upsilon}(\Ca)+ \Ca_{\upsilon}^{-1}\right)^{-1}\leq\phi_{t}(Q)\leq 
\phi^{\pi}_{t}(Q)\leq \phi_{\pi,t}(Q)\leq \Oa_{\pi,\upsilon_{\pi}}^{-1}+\Ca_{\pi,\upsilon_{\pi}}(\Oa).
\end{equation}
\end{theo}

\proof
By (\ref{after-upsilon}) and $\mbox{\rm (H)}_{1}$ we have the uniform
estimates
$$
\left(\Oa_{\pi,\upsilon_{\pi}}(\Ca)+ \Ca_{\pi,\upsilon_{\pi}}^{-1}\right)^{-1}\leq \phi_{\pi,t}(Q)\leq
\Oa_{\pi,\upsilon_{\pi}}^{-1}+\Ca_{\pi,\upsilon_{\pi}}(\Oa).
$$

We let $E_{\pi,t|s}(Q)$ be the transition semigroups  defined as $E_{t|s}(\phi^{\pi}_{s}(Q))$
by replacing $(A,\phi_{t})$ by $(A_{\pi},\phi_{\pi,t})$.
In this notation, the proof (\ref{backward-analysis-cons}) is a direct consequence of the backward perturbation formulae
\begin{equation}\label{backward-analysis-pi} 
\phi^{\pi}_{t}(Q)- \phi_{\pi,t}(Q)=\int_0^t~ E_{\pi,t|s}(\phi^{\pi}_{s}(Q))~ \Ra_{\pi}\left[\phi^{\pi}_{s}(Q)\right]~  E_{\pi,t\vert s}(\phi^{\pi}_{s}(Q))^{\prime}~ds\leq 0
\end{equation}
as well as
\begin{equation}\label{backward-analysis} 
\phi^{\pi}_{t}(Q)- \phi_{t}(Q)=\int_0^t~ E_{t\vert s}(\phi^{\pi}_{s}(Q))~ \Gamma_{\pi}\left[\phi^{\pi}_{s}(Q)\right]~  E_{t\vert s}(\phi^{\pi}_{s}(Q))^{\prime}~ds\geq 0.
\end{equation}

That is, the l.h.s. estimate in (\ref{backward-analysis-cons}) is a direct consequence of (\ref{after-upsilon}) and the relationship $\phi_{t}(Q)\leq \phi^{\pi}_{t}(Q)\leq \phi_{\pi,t}(Q)$ following from (\ref{backward-analysis-pi}) and (\ref{backward-analysis}). The r.h.s. estimate in (\ref{backward-analysis-cons}) follows obviously from the above.

To check (\ref{backward-analysis}) we use the interpolating path
$$
s\in [0,t]\mapsto \phi_{s,t}(\phi^{\pi}_{s}(Q))\quad\mbox{\rm from}\quad \phi_{t}(Q)\quad\mbox{\rm to}\quad \phi^{\pi}_{t}(Q).
$$
By Proposition~\ref{prop-key} we have
\begin{eqnarray*}
\partial_s\phi_{s,t}(\phi^{\pi}_{s}(Q))&=&-\ricc(\phi_{s,t}(\phi^{\pi}_{s}(Q)))+
 \partial\phi_{s,t}(\phi^{\pi}_{s}(Q))\cdot \partial_s\phi^{\pi}_{s}(Q)\\
 &=&
 \partial\phi_{s,t}(\phi^{\pi}_{s}(Q))\cdot\Gamma_{\pi}\left[\phi^{\pi}_{s}(Q)\right]=  E_{t\vert s}(\phi^{\pi}_{s}(Q))~\Gamma_{\pi}\left[\phi^{\pi}_{s}(Q)\right]~  E_{t\vert s}(\phi^{\pi}_{s}(Q)).
\end{eqnarray*}
This ends the proof of (\ref{backward-analysis}). The proof of (\ref{backward-analysis-pi}) follows the same arguments, thus it is skipped. This ends the proof of the theorem.
\qed

The next lemma compares the semigroups $\phi_{\pi,t}(Q)$ and $\phi_{t}(Q)$ when the matrices
$(A_{\pi},R_{\pi},S_{\pi})$ are close to $(A,R,S)$.

 \begin{lem}
Assume $\mbox{\rm (H)}_{0}$ and $\mbox{\rm (H)}_{1}$. For any $t\geq 0$ and $Q\in \SS^+_{r}$ we have
\begin{equation}\label{backward-analysis-pi-2} 
 \phi_{\pi,t}(Q)-\phi_t(Q)=\int_0^t~ E_{t|s}(\phi_{\pi,s}(Q))~  \Xi_{\pi}\left[\phi_{\pi,s}(Q)\right]~   E_{t|s}(\phi_{\pi,s}(Q))^{\prime}~ds
\end{equation}
as well as
\begin{equation}\label{backward-analysis-2}
\phi_{t}(Q)- \phi_{\pi,t}(Q)=\int_0^t~ E_{\pi,t\vert s}(\phi_{s}(Q))~  \Xi_{\pi}\left[\phi_{s}(Q)\right]~  E_{\pi,t\vert s}(\phi_{s}(Q))^{\prime}~ds.
\end{equation}

\end{lem}

The proof of this lemma follows the same arguments as the proof of Theorem~\ref{bplem}; thus it is skipped. Observe that
$$
(\ref{backward-analysis-pi-2})\Longrightarrow
P_{\pi}-P=\phi_t(P_{\pi})-P+\int_0^t~ E_{t|s}(P_{\pi})~  \Xi_{\pi}\left[P_{\pi}\right]~   E_{t|s}(P_{\pi})^{\prime}~ds
$$
and
$$
(\ref{backward-analysis-2})\Longrightarrow
P-P_{\pi}= \phi_{\pi,t}(P)-P_{\pi}+\int_0^t~ E_{\pi,t\vert s}(P)~  \Xi_{\pi}\left[P\right]~  E_{\pi,t\vert s}(P)^{\prime}~ds.
$$

\subsubsection{Robustness theorems}\label{robustness-sec}

We equip the set $\Ca(\SS_{r}^+,\SS_{r}^+)$ of continuous mappings $\pi~:~\SS_{r}^+\mapsto \SS_{r}^+$ with the uniform norm
$$
\Vert \pi_1-\pi_2\Vert=\sup_{Q\in \SS_{r}^+}\Vert \pi_1(Q)-\pi_2(Q)\Vert_2.
$$

Let $\Pi\subset \Ca(\SS_{r}^+,\SS_{r}^+)$ be a compact subset, and let $t>0$ be some fixed time horizon.
For any $\delta>0$, we let $B(\delta)$ be the $\delta$-ball around the identity mapping; that is
$$
B(\delta):=\{\pi\in \Pi~:~\Vert \pi-id\Vert\leq \delta\}.
$$
We consider the following continuity condition 
$$
\begin{array}{l}
\hskip-1cm\mbox{\rm (H)}_2\qquad 
\forall\epsilon,\alpha\in ]0,1]\quad  \exists\delta>0\quad \mbox{\rm such that}\quad\forall \pi\in B(\delta)\quad~\mbox{\rm we have}~\\
\\
\hskip1cm \alpha~S~\leq S_{\pi}\leq \alpha^{-1}~S
\qquad \alpha~R~\leq R_{\pi}\leq \alpha^{-1}~R\qquad\mbox{\rm and}\qquad \Vert A-A_{\pi}\Vert_2 \leq \epsilon.
\end{array}
$$

Assume that $\mbox{\rm (H)}_0$ and $\mbox{\rm (H)}_1$ are met. Importantly, in this situation we have
\begin{equation}\label{pi-above-below}
\Vert\phi^{\pi}_{t}(Q)-\phi_{t}(Q)\Vert_2\leq \Vert\phi_{\pi,t}(Q)-\phi_t(Q)\Vert_2
\end{equation}
We check this claim using the fact that
$$
(\ref{backward-analysis-cons}) \Longrightarrow~  0\leq \phi^{\pi}_{t}(Q)-\phi_{t}(Q)\leq \phi_{\pi,t}(Q)-\phi^{\pi}_{t}(Q) + \phi^{\pi}_{t}(Q)-\phi_{t}(Q) = \phi_{\pi,t}(Q)-\phi_{t}(Q).
$$
Of course, it is important to note that we can take $\mbox{\rm (H)}_2$ without $\mbox{\rm (H)}_0$, and consider directly just the flow $\phi_{\pi,t}$. The main objective of this section is to prove the following robustness theorem.
\begin{theo}\label{theo-continuity}
Let $(A_{\pi},R_{\pi}, S_{\pi})$ be a collection of matrices satisfying condition 
 $\mbox{\rm (H)}_2$. In this case, there exists some $\delta>0$ such that
for any $\pi\in B(\delta)$, any horizon $t\geq 0$ and any $Q\in \SS_{r}^+$ we have
\begin{eqnarray*}
\Vert\phi_{\pi,t}(Q)-\phi_t(Q)\Vert_2 ~&\leq&~ \left[ \chi_1(\delta)+ e^{-4t\nu}~\chi_2(\delta,\Vert Q\Vert_2)\right]~\gamma(\pi)
\end{eqnarray*}
for some finite constants $\chi_1(\delta)$ and $\chi_2(\delta,\Vert Q\Vert_2 )$, whose values only depend on the parameters $\delta$, and $(\delta,\Vert Q\Vert_2)$ respectively. In particular we have
$$
\forall \pi\in B(\delta)\qquad
\Vert P_{\pi}-P\Vert_2\leq 
\chi_1(\delta)~\gamma(\pi).
$$
\end{theo}

Note that whenever we take $\mbox{\rm (H)}_0$ as holding, then we take $\mbox{\rm (H)}_2$ as being compatible in the definition of $(A_{\pi},R_{\pi}, S_{\pi})$ as given in (\ref{def-Api}). This allows us to compare $\phi^\pi_{t}$ and $\phi_{\pi,t}$. In this case, (\ref{pi-above-below}) means this theorem guarantees the boundedness of those perturbation models satisfying (\ref{pi-Riccati-def}), (\ref{pi-Riccati-def-gamma-term}) and $\mbox{\rm (H)}_{0}$. See also \cite{ap-2018-franklin} for a refined/corrected discussion on this and related results.

Note also however, that $\mbox{\rm (H)}_{2}$, and Theorem~\ref{theo-continuity}, capture a broader class of perturbation model than (\ref{pi-Riccati-def}), (\ref{pi-Riccati-def-gamma-term}) and $\mbox{\rm (H)}_{0}$ alone. In particular, $\mbox{\rm (H)}_{2}$ is simply concerned with direct perturbations of the original $(A,R,S)$ system matrices in the Riccati operator (\ref{def-ricc}). 

The proof of the preceding theorem relies on the following proposition.

\begin{prop}\label{theo-comparison}
When $\mbox{\rm (H)}_2$ is met,
 for any $\alpha\in ]0,1]$ there exists some $\delta>0$ such that the matrices $(A_{\pi},R_{\pi},S_{\pi})$ indexed by mappings $\pi$ in the $\delta$-ball
$B(\delta)$  satisfy the  Gramian condition with a common 
interval of observability-controllability  $\upsilon_{\pi}=\upsilon$ and
 some parameters 
 $$
 \alpha~\varpi^{o,c}_{\pm}\leq 
 \varpi^{o,c}_{\pi,\pm}\leq \alpha^{-1}~\varpi^{o,c}_{\pm}$$
 and well as
 $$  \alpha~\varpi^{c}_{\pm}(\Oa)\leq 
 \varpi^{c}_{\pi,\pm}(\Oa)\leq \alpha^{-1}~\varpi^{c}_{\pm}(\Oa)\quad\mbox{and}\quad
  \alpha~\varpi^{o}_{\pm}(\Ca)\leq 
 \varpi^{o}_{\pi,\pm}(\Ca)\leq \alpha^{-1}~\varpi^{o}_{\pm}(\Ca).
$$
\end{prop}

We already quote a direct consequence of Theorem~\ref{bplem} and Proposition~\ref{theo-comparison}.

\begin{cor}\label{cor-H2}
Assume $\mbox{\rm (H)}_2$. In this situation, for any $\alpha\in ]0,1]$ there exist some $\delta>0$ such that for any $\pi\in B(\delta)$, any $Q\in\SS_{r}^+$ and any $t\geq 0$ we have the common uniform estimates
\begin{equation}\label{backward-analysis-cons-2} 
	\alpha\left(\varpi^{o}_{+}(\Ca)+ 1/\varpi^{c}_{-}\right)^{-1}Id~\leq~\phi_{t+\upsilon}(Q),\, \phi_{\pi,t+\upsilon}(Q)~\leq~ \alpha^{-1}\left[\varpi^{c}_{+}(\Oa)+1/\varpi^{o}_{-}\right]Id.
\end{equation}
If $\mbox{\rm (H)}_0$ also holds, then we know additionally that $\phi_{t}(Q)\leq \phi^{\pi}_{t}(Q)\leq \phi_{\pi,t}(Q)$.
\end{cor}

We have already studied the ordering between $\phi_{t}$, $\phi^{\pi}_{t}$ and $\phi_{\pi,t}$ in Theorem \ref{bplem} when $\mbox{\rm (H)}_{0}$ holds. But if we take $\mbox{\rm (H)}_2$ as holding alone (i.e. forget $\mbox{\rm (H)}_0$ and $\phi^{\pi}_{t}$), then we cannot compare $\phi_{\pi,t}(Q)$ with $\phi_{t}(Q)$ in the same way. To be more specific, $\phi_{\pi,t+\upsilon}(Q)\geq\phi_{t+\upsilon}(Q)$ holds when $\Xi_{\pi}(Q)\geq0$, and $\phi_{\pi,t+\upsilon}(Q)\leq\phi_{t+\upsilon}(Q)$ holds when $\Xi_{\pi}(Q)\leq0$. To check this simply see (\ref{backward-analysis-pi-2}). 

The proof of the above Proposition \ref{theo-comparison} relies on the next couple of comparison lemmas of interest on their own.

\begin{lem}\label{comp-lem}
Let $V_1,V_2\in\SS_{r}^+$ be a couple of definite positive matrices s.t. $V_1\geq V_2$. We set
$$
	Q_1:=U_1V_1U_1^{\prime}\quad\mbox{and}\quad Q_2:=U_2V_2U_2^{\prime}
$$
for some $(U_1,U_2)\in\MM_{r}^2$. Assume that 
$
Q_2\geq q_2~Id
$,
for some $q_2>0$. With $U_2$ invertible we have
$$
\Vert Q_2\Vert_2~\Vert U_1U_2^{-1}-Id\Vert_2~<\sqrt{1+q_2}-1\Longrightarrow
Q_1\geq q_{1,2}~
Q_2
$$
with 
$$
q_{1,2}=\left[1- q_2^{-1}\left\{\left(1+\Vert Q_2\Vert_2~\Vert  U_1U_2^{-1}-Id\Vert_2 \right)^2-1\right\}\right].
$$
\end{lem}
\proof
We set
$$
U_2U_1^{-1}=Id+U_{1,2}\quad\mbox{\rm and}\quad U_1U_2^{-1}=Id+U_{2,1}.
$$
Observe that
\begin{eqnarray*}
Q_2\geq q_2~Id&\Rightarrow &\lambda_{min}(Q_2)\geq q_2\Rightarrow \lambda_{min}(Q^{1/2}_2)\geq \sqrt{q_2}\\
&\Rightarrow& \lambda_{max}(Q^{-1/2}_2)\geq 1/\sqrt{q_2}\Rightarrow\Vert Q^{-1/2}_2\Vert_2^2\leq q_2^{-1}.
\end{eqnarray*}
In this situation, we have
$$
\begin{array}{l}
\displaystyle  Q_1\geq U_1U_2^{-1}~Q_2~  (U_1U_2^{-1})^{\prime}\\
\\
\qquad \Longrightarrow
\displaystyle  Q_1^{-1}\leq ~Q_2^{-1/2}\left[ Q_2^{-1/2}U_1U_2^{-1}~Q_2~  (U_1U_2^{-1})^{\prime}Q_2^{-1/2}\right]^{-1}Q_2^{-1/2}.
\end{array}
$$
Observe that
$$
\begin{array}{l}
\displaystyle ~\left[ Q_2^{-1/2}U_1U_2^{-1}~Q_2~  (U_1U_2^{-1})^{\prime}Q_2^{-1/2}\right]^{-1}\\
\\
\displaystyle \qquad\qquad\qquad\qquad =\left[ Id-Q_2^{-1/2}\left\{Q_2-U_1U_2^{-1}~Q_2~  (U_1U_2^{-1})^{\prime}\right\}Q_2^{-1/2}\right]^{-1}\\
\\
\displaystyle \qquad\qquad\qquad\qquad =\sum_{n\geq 0}\left[Q_2^{-1/2}\left\{Q_2-U_1U_2^{-1}~Q_2~  (U_1U_2^{-1})^{\prime}\right\}Q_2^{-1/2}\right]^{n}.
\end{array}
$$

On the other hand, we have
$$
\begin{array}{l}
\displaystyle \left\{Q_2-U_1U_2^{-1}~Q_2~  (U_1U_2^{-1})^{\prime}\right\}\\
\\
\qquad\qquad=\displaystyle \left\{Q_2-\left[Id+U_{2,1}\right]~Q_2~ \left[Id+U_{2,1}^{\prime}\right]\right\}
~~=~~\displaystyle -\left[U_{2,1}Q_2+Q_2U_{2,1}^{\prime}\right]-U_{2,1}Q_2U_{2,1}^{\prime}.
\end{array}
$$
This implies that
$$
\begin{array}{l}
\displaystyle \Vert Q_2^{-1/2}\left\{Q_2-U_1U_2^{-1}~Q_2~  (U_1U_2^{-1})^{\prime}\right\}Q_2^{-1/2}\Vert\\
\\
\qquad\qquad \leq \displaystyle \displaystyle q_2^{-1}~\Vert U_{2,1}\Vert_2\Vert Q_2\Vert_2 ~\left[2+\Vert U_{2,1}\Vert_2\Vert Q_2\Vert_2\right]
~~=~~ \displaystyle \displaystyle q_2^{-1}\left\{\left(1+\Vert U_{2,1}\Vert_2\Vert Q_2\Vert_2 \right)^2-1\right\}~<~1
\end{array}
$$
from which we conclude that
$$
\left[ Q_2^{-1/2}U_1U_2^{-1}~Q_2~  (U_1U_2^{-1})^{\prime}Q_2^{-1/2}\right]^{-1}~ \leq~ \left[1- q_2^{-1}\left\{\left(1+\Vert U_{2,1}\Vert_2\Vert Q_2\Vert_2 \right)^2-1\right\}\right]^{-1}~Id.
$$
This yields the estimate
$$
 Q_1^{-1}\leq q_{1,2}^{-1}~
Q_2^{-1}\quad\Longleftrightarrow\quad  Q_1\geq q_{1,2}~
Q_2
$$
with
$$
q_{1,2}^{-1}=~ \left[1- q_2^{-1}\left\{\left(1+\Vert U_{2,1}\Vert_2\Vert Q_2\Vert_2 \right)^2-1\right\}\right]^{-1}.
$$
This ends the proof of the lemma.\qed

\begin{lem}\label{lemm-comp-flows}
Let $\Ua,\Va$ be a pair of bounded functions from $[0,t]\times\Pi$ into $\SS_{r}^+$. We consider the integral mappings
$$
(s,\pi)\in ([0,t]\times\Pi)\mapsto
\Wa_s(\pi):=\int_0^s \Ua_u(\pi)  \Va_u(\pi)~\Ua^{\prime}_u(\pi)~du\in \SS_{r}^+.
$$
Let $\pi_1,\pi_2\in \Pi$  be such that
\begin{equation}\label{comparison-Wa-condition}
\forall s\in [0,t]\quad\Va_s(\pi_1)\geq \Va_s(\pi_2)\quad\mbox{and}\quad 
\Wa_t(\pi_2)\geq \varpi_{-,t}(\pi_2)~Id\quad\mbox{\rm for some}\quad \varpi_{-,t}(\pi_2)>0.
\end{equation}
Also assume that the flow of matrices $\Ua_s(\pi_2)$ are invertible
for any $s\in [0,t]$ and they satisfy the following Lipschitz inequality
\begin{equation}\label{hypothesis-Lip-Ua}
\sup_{s\in [0,t]}{\Vert~ \Ua_s(\pi_1)~\Ua_s(\pi_2)^{-1}-Id\Vert_2}\leq \lip_t(\Ua)~\Vert \pi_1-\pi_2\Vert
\end{equation}
for some finite constant $ \lip_t(\Ua)$.
In this situation, for any $\epsilon\in ]0,1]$ there exists some parameter $\delta=\delta(t,\epsilon,\pi_2)>0$
such that
$$
\Vert \pi_1-\pi_2\Vert\leq \delta\Longrightarrow
\Wa_t(\pi_1)\geq \epsilon~\Wa_t(\pi_2).
$$
\end{lem}
The proof of this lemma follows the same arguments as those in the proof of Lemma~\ref{comp-lem}. For the convenience of the reader a detailed proof of the lemma is given in the Appendix.

Now we come to the proof of Proposition~\ref{theo-comparison}.

\noindent {\bf Proof of Proposition~\ref{theo-comparison}:}

We assume that for any $\epsilon\in ]0,1]$ there exists some $\delta>0$ such that
$$
\Vert \pi-id\Vert\leq \delta\Longrightarrow
S_{\pi}\geq (1-\epsilon)~S\quad\mbox{\rm and}\quad R_{\pi}\geq (1-\epsilon)~R.
$$
We apply Lemma~\ref{lemm-comp-flows} to the functions
\begin{eqnarray*}
 \Wa_t^o(\pi):=(1-\epsilon)^{-1}~\Oa_{\pi,t}&\Longrightarrow &\Wa_t^o(id)=(1-\epsilon)^{-1}\Oa_t\\
 \Wa_t^c(\pi):=(1-\epsilon)^{-1}~\Ca_{\pi,t}&\Longrightarrow &\Wa_t^c(id)=(1-\epsilon)^{-1}\Ca_t
\end{eqnarray*}
with $(\pi_1,\pi_2)=(\pi,id)$ and the time horizon $t=\upsilon$. For any $\epsilon_1\in ]0,1]$ there exists
some parameter $\delta_1=\delta(\epsilon_1,\upsilon)$ such that
$$
\Vert \pi-id\Vert\leq \delta_1\Longrightarrow
\Oa_{\pi,\upsilon}\geq \epsilon_1~\Oa_{\upsilon}
\quad\mbox{\rm and}\quad
\Ca_{\pi,\upsilon}\geq \epsilon_1~\Ca_{\upsilon}.
$$

We assume that for any $\epsilon\in ]0,1]$ there exists some $\delta>0$ such that
$$
\Vert \pi-id\Vert\leq \delta\Longrightarrow
S\geq (1-\epsilon)~S_{\pi}\quad\mbox{\rm and}\quad R\geq (1-\epsilon)~R_{\pi}.
$$
We apply Lemma~\ref{lemm-comp-flows} to the functions 
$$
 \Wa_t^o(\pi)=\Oa_{\pi,t}\quad\mbox{\rm and}\quad  \Wa_t^c(\pi)=\Ca_{\pi,t}
$$
with $(\pi_1,\pi_2)=(id,\pi)$ and the time horizon $t=\upsilon$. From previous estimates we have
$$
\Vert \pi_2-id\Vert\leq \delta_1\Longrightarrow
 \Wa_{\upsilon}^o(\pi_2)\geq~\epsilon_1~\varpi_-^o~Id\quad
 \mbox{\rm and}\quad  \Wa_{\upsilon}^c(\pi_2)\geq~\epsilon_1~\varpi_-^c~Id.
$$
By Lemma~\ref{lemm-comp-flows} for any $\epsilon_2\in ]0,1]$ we can choose $\delta=\delta(\epsilon_1,\epsilon_2,\upsilon)$ such that
$$
\Vert \pi-id\Vert\leq \delta\Longrightarrow \Oa_{\upsilon}\geq \epsilon_2~\Oa_{\pi,\upsilon}\quad \mbox{\rm and}\quad
 \Ca_{\upsilon}\geq \epsilon_2~\Ca_{\pi,\upsilon}.
$$
This shows that
$$
\Vert \pi-id\Vert\leq \delta\Longrightarrow
 \epsilon_1~\Oa_{\upsilon}\leq 
~\Oa_{\pi,\upsilon}\leq  \epsilon_2^{-1}~\Oa_{\upsilon}\quad\mbox{\rm and}\quad
 \epsilon_1~\Ca_{\upsilon}\leq 
~\Ca_{\pi,\upsilon}\leq  \epsilon_2^{-1}~\Ca_{\upsilon}.
$$
In the same vein we prove the estimates of the Gramians $\Oa_{\pi,\upsilon}(\Ca)$ and  $\Ca_{\pi,\upsilon}(\Oa)$. This ends the proof of the proposition.\qed

We are now in position to prove Theorem~\ref{theo-continuity}

\noindent {\bf Proof of Theorem~\ref{theo-continuity}:}

By Corollary~\ref{cor-H2}, there exist some $\delta>0$ s.t. for any $\pi\in B(\delta)$ we have the uniform estimate
$$
\sup_{t\geq 0}{\sup_{Q\in\SS_{r}^+}{\Vert \Xi_{\pi}(\phi_{\pi,t+\upsilon}(Q))\Vert_2}} ~\leq~ \chi_1(\delta)~\left[\Vert A_{\pi}-A\Vert_2+\Vert R_{\pi}-R\Vert_2+\Vert S_{\pi}-S\Vert_2 \right]
$$
for some finite constant $\chi_1(\delta)$ whose values only depend on $\delta$. We have
$$
(\ref{upper-bound-Phi})~~\Longrightarrow~~\sup_{\pi\in B(\delta)}{\sup_{0\leq t\leq \upsilon}{\Vert \phi_{\pi,t}(Q)\Vert_2}} ~\leq~ \chi_2(\delta)~(1+\Vert Q\Vert_2)
$$
for some constant $\chi_2(\delta)$ whose values only depend on $\delta$.
Combining (\ref{expo-E}) with (\ref{backward-analysis-pi-2}) we have
$$
\begin{array}{l}
 \Vert\phi_{\pi,t}(Q)-\phi_t(Q)\Vert_2
 \displaystyle~\leq~ \left[\kappa_E(\chi_2(\delta)~(1+\Vert Q\Vert_2))\right]^2 \\
 \\
\hskip2cm\times ~\left[2\chi_2(\delta)~(1+\Vert Q\Vert_2)\Vert A_{\pi}-A\Vert_2+\Vert R_{\pi}-R\Vert_2+\chi_2(\delta)^2~(1+\Vert Q\Vert_2)^2\Vert S_{\pi}-S\Vert_2 \right]~\\
\\
\hskip4cm\times
[e^{-4\nu (t-\upsilon)}-e^{-4t\nu}]/(4\nu)\\
 \\
  \hskip2cm\displaystyle+\chi_1(\delta)~ \kappa_E(\chi_1(\delta))^2~\left[\Vert A_{\pi}-A\Vert_2+\Vert R_{\pi}-R\Vert_2+\Vert S_{\pi}-S\Vert_2 \right]
\left[1-e^{-4\nu (t-\upsilon)}\right]/(4\nu).
\end{array}
$$
We collect constants in the first term with $\chi_2(\delta,\Vert Q\Vert )$ and in the second term (via some notation abuse) with $\chi_1(\delta)$. This ends the proof of the first assertion. 

To check the last assertion we simply let $t\uparrow\infty$. More precisely,  observe that
$$
 \Vert P_{\pi}-P\Vert_2 ~\leq~ \Vert\phi_{\pi,t}(P_{\pi})-\phi_t(P_{\pi})\Vert_2+ \Vert\phi_{t}(P_{\pi})-\phi_t(P)\Vert_2.
$$
This implies, using (\ref{expo-phi}) and the first assertion of this theorem, that
$$
\begin{array}{l}
\Vert P_{\pi}-P\Vert_2 ~\leq~  \kappa_{\phi}(\Vert P_{\pi}\Vert_2,\Vert P\Vert_2)~e^{-2\nu t}~\Vert P_{\pi}-P\Vert_2\\
\\
\qquad\qquad\qquad\qquad +
\left[ \chi_1(\delta)+ e^{-4t\nu}~\chi_2(\delta,\Vert P_{\pi}\Vert_2)\right]
\left[\Vert A_{\pi}-A\Vert_2+\Vert R_{\pi}-R\Vert_2+\Vert S_{\pi}-S\Vert_2 \right].
\end{array}
$$
Letting $t\rightarrow\infty$ we end the proof of the desired estimate. This ends the proof of the theorem.\qed

See also \cite{ap-2018-franklin} for further discussion on these robustness and related results.

\subsection{Projection-type models}\label{projection-models-sec}

We consider projection-type mappings in (\ref{pi-Riccati-def}) of the second type in (\ref{two-classes}). Let $\pi$ be some positive map from $\MM_r$ into itself; that is $\pi(\SS_{r}^+)\subseteq \SS^{+}_{r}$. We first assume the matrices $(A,R,S)$ satisfy $$(\pi(A),\pi(A^{\prime}),\pi(S),\pi(R))=(A,A^{\prime},R,S)
$$
and we let $\Ba\subset\MM_r$ be a given {\em matrix ring}. Also assume that the pair $(\pi,\Ba)$ satisfies the following orthogonality property:
  \begin{equation}\label{commutation}
\hskip-1cm\mbox{\rm (H)}_3\hskip2cm \forall Q\in\MM_r\qquad\forall B\in\Ba
\qquad \pi(B[Q-\pi(Q)]+[Q-\pi(Q)]B)=0.
 \end{equation}
In this situation, we have
\begin{eqnarray*}
\pi\left[(Q-B)(Q-B)^{\prime}\right]&=&\pi\left[(Q-\pi(Q))(Q-\pi(Q))^{\prime}\right]+\pi(BB^{\prime})\\
&\geq&
\pi\left[(Q-\pi(Q))(Q-\pi(Q))^{\prime}\right] ~\geq~ 0.
\end{eqnarray*}
This shows that $\pi$ can be interpreted as a $\pi$-orthogonal projection
$$
\mbox{\rm (H)}_3~\Longleftrightarrow~
\pi(Q)=\argmin_{B\in \Ba}{\pi\left[(Q-B)(Q-B)^{\prime}\right]}.
$$
In addition, we have the Cauchy-Schwartz inequality
$$
\pi\left[(Q-\pi(Q))(Q-\pi(Q))^{\prime}\right]\geq 0 ~\Longrightarrow~ \pi(QQ^{\prime})\geq \pi(Q)\pi(Q^{\prime}).
$$
Whenever $\Ba$ is closed by transposition we have
$$
\pi(Q^{\prime})=\pi(Q)^{\prime}~ \Longrightarrow ~\pi(Q)\pi(Q)^{\prime}\leq \pi(QQ^{\prime}).
$$
The identity mapping $\pi=id$ and the set $\Ba=\MM_r$ clearly satisfies the above properties. 

The prototype of a non-trivial pair $(\Ba,\pi)$ satisfying $\mbox{\rm (H)}_3$ are orthogonal projections $\pi=\proj_{\Ba}$ (w.r.t. the Frobenius norm) onto cellular (a.k.a. coherent) algebras; that is a sub-algebra of matrices which are additionally closed under the Hadamard-Schur product and contains the identity elements $Id$ and $J$, where $J$ stands for the matrix with all ones entries. Up to a unitary change of basis, these projections can be reformulated in terms of block-diagonal matrices~\cite{barker}. By~\cite[page 57]{brouwer}, a sub-algebra of $\SS_r$ is a Bose-Mesner algebra \cite{bose} of some association scheme if and only if it contains $I$ and $J$, and it is closed under the Hadamard-Schur product. This shows that cellular sub-algebras of $\SS_r$ coincide with the Bose-Mesner algebras (of some association scheme). We refer to Section~\ref{Bose-Mesner-sec} for a detailed discussion on these models.

The set $\Ba=\Ma_{r[1]}\oplus\ldots\oplus\Ma_{r[n]}\subset \MM_r$  (with $r=\sum_{1\leq q\leq n}r[q]$) of block-diagonal matrices with null entries outside the blocks is also a matrix ring which is closed under the Hadamard-Schur product w.r.t. any matrix in $\MM_r$; but $\Ba$ is not a cellular algebra since $J\not\in\Ba$. The orthogonal projection from $\MM_r$ onto this $\Ba$ is given by
\begin{equation}\label{ex-Ba}
	\pi(Q):=L\odot Q\quad\mbox{\rm with the block-diagonal matrix}\quad
		 L:=\mbox{\rm diag}(J_1,\ldots,J_n)\geq 0.
\end{equation}
In the above display, $J_i$ stands for the $i$-th block unit matrix w.r.t. the Hadamard-Schur product; that is the $(r[i]\times r[i])$-square matrix with all unit entries. It is readily checked that $(\Ba,\pi)$ satisfies condition $\mbox{\rm (H)}_3$. We refer to Section~\ref{block-diag-sec} for a discussion on these models.

We let $\phi^{\pi}_t$ be the $\pi$-Riccati semigroup defined in Section~\ref{statement-sec}. By (\ref{backward-analysis}) we have the domination property
\begin{equation}\label{domination}
	\forall Q \in\SS_{r}, \qquad \phi^{\pi}_t(Q)\geq \phi_t(Q).
\end{equation}

In contrast with the second order approximation models discussed in Section~\ref{second-order-appr-sec} these projection techniques don't depend on some perturbation index that quantifies the distance between $\pi$ and the identity mapping. 

When $(\pi(A),\pi(S),\pi(R))\not=(A,R,S)$ we can replace $(A,R,S)$ by their projections $(A_{\pi},R_{\pi},S_{\pi})$. In this case,
$\phi_{\pi,t}$ is the Riccati semigroup associated with the drift function $\ricc_{\pi}(Q)$ defined simply by $\ricc(Q)$ with $(A,R,S)$ replaced by $(A_{\pi},R_{\pi},S_{\pi})$. The difference between $\phi_{\pi,t}$ and $\phi_t$ can be analyzed as in Theorem~\ref{theo-continuity}. It is not possible to ensure that $\phi_{\pi,t}$ is arbitrarily close to $\phi_t$ without some continuity conditions.

Section~\ref{stein-sec} discusses a way to combine these projection-type models with the perturbation-type models discussed in Section~\ref{second-order-appr-sec}

In the latter development of Section~\ref{expo-concentration} we will provide exponential concentration inequalities that ensure the $\pi$-projected Riccati flows converge exponentially fast to the solution of the (nominal, Kalman-Bucy) Riccati equation, viz \cite{kalman61,ap-2016}, as the time horizon tends to $\infty$, and as soon as condition $\mbox{\rm (H)}_3$ is met. Speaking somehow loosely we shall show that
$$
	 \phi_t\circ\pi=\phi^{\pi}_t\circ\pi\qquad  \mbox{\rm and} \qquad \phi^{\pi}_t~ \simeq ~ \phi_t.
$$

\begin{defi}
We let 
$\varphi^{\pi}_{t}(Q)$ be the flow of the projected Riccati equation
$$
\partial_t\varphi^{\pi}_{t}(Q)=\pi\left[\ricc\left(\varphi^{\pi}_{t}(Q)\right)\right]
$$
\end{defi}

The next theorem shows that the flows $\varphi^{\pi}_t$ and $\phi^{\pi}_t$ coincide with the $\pi$-projection of the Riccati flow $\pi\circ\phi_t$ as soon as we start from an initial state $Q$ that satisfies $\pi(Q)=Q$ and $(A,R,S)\in\Ba^3$. It also provides an explicit description of the flow $\phi^{\pi}_t$ in terms of $\phi_t$ and $\pi$ when $\mbox{\rm (H)}_3$ is satisfied.
\begin{theo}\label{theo-pi-1}
Assume $\mbox{\rm (H)}_{3}$ and recall $\ricc^{\pi}$ defined in (\ref{pi-Riccati-def}). For any time horizon $t\geq 0$  we have the formula
$$
\pi\circ\ricc\circ\pi=\ricc^{\pi}\circ\pi=\ricc\circ\pi=\pi\circ\ricc^{\pi}\circ\pi
$$
as well as the
semigroup commutation properties
\begin{equation}\label{commutation-prop-ref}
\pi\circ\phi_t\circ\pi=\pi\circ\varphi^{\pi}_{t}\circ\pi=\varphi^{\pi}_{t}\circ\pi=\phi_t\circ\pi =\phi^{\pi}_t\circ\pi.
\end{equation}
In addition, we have the formula
\begin{equation}\label{explicit-formula}
\phi_t(Q)~\leq~ \phi^{\pi}_t(Q)= \left[\phi_t\circ \pi\right](Q)+E_t(\pi(Q))(Q-\pi(Q))E_t(\pi(Q))^{\prime}.
\end{equation}
\end{theo}
\proof
Recall that $(A,A^{\prime},R,S)\in \Ba$. Since $\Ba$ is a matrix ring we have
\begin{eqnarray*}
\pi\left[A\pi(Q)+\pi(Q)A^{\prime}+R-\pi(Q)S\pi(Q)\right]=A\pi(Q)+\pi(Q)A^{\prime}+R-\pi(Q)S\pi(Q)
\end{eqnarray*}
or equivalently 
$$
\pi\circ\ricc\circ\pi=\ricc\circ\pi.
$$
Also observe that
\begin{eqnarray*}
\ricc^{\pi}(Q)&=&(A-\pi(Q)S)Q+Q(A-\pi(Q)S)^{\prime}+R+\pi(Q)S\pi(Q)\\
&=&\ricc(\pi(Q))+
(A-\pi(Q)S)(Q-\pi(Q))+(Q-\pi(Q))(A-\pi(Q)S)^{\prime}\\
&=&\left[\pi\circ\ricc\circ\pi\right](Q)+
(A-\pi(Q)S)(Q-\pi(Q))+(Q-\pi(Q))(A-\pi(Q)S)^{\prime}
\end{eqnarray*}
and thus
$$
\ricc^{\pi}\circ\pi=\pi\circ\ricc\circ\pi=\pi\circ\ricc^{\pi}\circ\pi.
$$

Now, we also have
$$
\pi^2:= \pi\circ\pi=\pi ~~\Rightarrow~~\partial_t\pi(\varphi^{\pi}_{t}(Q))=\partial_t\varphi^{\pi}_{t}(Q)~~\Rightarrow~~ \pi(\varphi^{\pi}_{t}(Q))=\varphi^{\pi}_{t}(Q)
+\pi(Q)-Q.
$$
This implies that
$$
\pi\circ\varphi^{\pi}_{t}\circ\pi=\varphi^{\pi}_{t}\circ\pi.
$$
This yields
\begin{eqnarray*}
\partial_t \left[\varphi^{\pi}_{t}\circ\pi\right](Q)&=&\partial_t \left[\pi\circ\varphi^{\pi}_{t}\circ\pi\right](Q)~=~
\pi\left(\left[\ricc\circ\varphi^{\pi}_{t}\circ\pi\right](Q)\right)\\
&=&
\pi\left(\left[\ricc\circ\pi\circ\varphi^{\pi}_{t}\circ\pi\right](Q)\right)~=~\left[\ricc\circ\pi\circ\varphi^{\pi}_{t}\circ\pi\right](Q)\\
&=&\left[\ricc\circ
\varphi^{\pi}_{t}\circ\pi\right](Q)
\end{eqnarray*}
and by the uniqueness of the solution of the Riccati equation we conclude that 
$$
\varphi^{\pi}_{t}\circ\pi=
\phi_t\circ\pi ~~\Rightarrow~~\pi\circ\phi_t\circ\pi=\pi\circ\varphi^{\pi}_{t}\circ\pi=\varphi^{\pi}_{t}\circ\pi=\phi_t\circ\pi.
$$
We also have
\begin{eqnarray*}
\partial_t\left[\varphi^{\pi}_{t}\circ\pi\right](Q)&=&\partial_t\left[\pi\circ\varphi^{\pi}_{t}\circ\pi\right](Q)
~=~\pi\left[\ricc^{\pi}\left(\left[\pi\circ\varphi^{\pi}_{t}\circ\pi\right](Q)\right)\right]\\
&=&
\ricc^{\pi}\left(\left[\pi\circ\varphi^{\pi}_{t}\circ\pi\right](Q)\right)~=~
\ricc^{\pi}\left(\left[\varphi^{\pi}_{t}\circ\pi\right](Q)\right)
\end{eqnarray*}
which implies that
$$
\phi^{\pi}_t\circ\pi=\varphi^{\pi}_{t}\circ\pi.
$$
This completes the proof of (\ref{commutation-prop-ref}).

Now, we have
 $$
 \begin{array}{l}
  \partial_t\left[\phi_t^{\pi}(Q)- \phi_t(\pi(Q))\right]\\
  \\
  \qquad=\ricc(\pi(\phi_t^{\pi}(Q)))-\ricc(\phi_t(\pi(Q)))\\
  \\
  \qquad\qquad\qquad  +
  (A-\phi_t(\pi(Q))S)(\phi_t^{\pi}(Q)-\pi(\phi_t^{\pi}(Q)))+(\phi_t^{\pi}(Q)-\pi(\phi_t^{\pi}(Q)))(A-\phi_t(\pi(Q))S)^{\prime}\\
  \\
  \qquad  = (A-\phi_t(\pi(Q))S)(\phi_t^{\pi}(Q)-\pi(\phi_t^{\pi}(Q)))+(\phi_t^{\pi}(Q)-\pi(\phi_t^{\pi}(Q)))(A-\phi_t(\pi(Q))S)^{\prime}.
\end{array} $$
This implies that
$$
\phi_t^{\pi}(Q)- \phi_t(\pi(Q))=E_t(\pi(Q))(Q-\pi(Q))E_t(\pi(Q))^{\prime}.
$$
The l.h.s. estimate in (\ref{explicit-formula}) is a consequence of the domination property (\ref{domination}).
  This ends the proof of the theorem. \qed

\subsubsection{Exponential contraction inequalities}\label{expo-concentration}
  
We continue with the projection-type models and $\mbox{\rm (H)}_{3}$ holding.   
\begin{theo}\label{exp-contraction-proj}
For any $Q_1,Q_2\in\SS^{+}_{r}$ and for any $t\geq 0$ we have
\begin{eqnarray}\label{expo-phi-pi}
\Vert \phi_t^{\pi}(Q_2)-\phi_t^{\pi}(Q_1)\Vert_2 ~\leq~ \kappa_{\phi^{\pi}}(Q_1,Q_2)~e^{-2\nu t}~\left[\Vert \pi(Q_2)-\pi(Q_1)\Vert_2+e^{-2\nu t}~\Vert Q_2-Q_1\Vert_2\right]
\end{eqnarray}
some finite constant  $\kappa_{\phi^{\pi}}(Q_1,Q_2)<\infty$ whose values only depend on $\left(\Vert Q_1\Vert_2,\Vert Q_2\Vert_2\right)$. This implies the existence of an unique fixed point $P^{\pi}=\phi^{\pi}_t(P^{\pi})$ with $\pi(P^{\pi})=P$. In addition, for any $Q\in\SS^{+}_{r}$ and for any $t\geq 0$ we have
\begin{equation}\label{expo-phi-pi-bis}
	\Vert \pi\left[\phi_t^{\pi}(Q)\right]-\phi_t^{\pi}(Q)\Vert_2
		~\leq~ \kappa_{\phi^{\pi}}(Q,\pi(Q))~e^{-4\nu t}~\Vert \pi(Q)-Q\Vert_2.
\end{equation}
 \end{theo}

\proof
We have
\begin{eqnarray*}
\phi^{\pi}_t(Q)&=&\phi^{\pi}_t(\pi(Q))+E_t(\pi(Q))(Q-\pi(Q))
E_t(\pi(Q))^{\prime}\\
&=&\phi_t(\pi(Q))+E_t(\pi(Q))(Q-\pi(Q))
E_t(\pi(Q))^{\prime}.
\end{eqnarray*}
This implies that
$$
\begin{array}{l}
\phi^{\pi}_t(Q_1)-\phi^{\pi}_t(Q_2) \,=\, \phi_t(\pi(Q_1))-\phi_t(\pi(Q_2))\\
\\
	\quad\qquad\qquad\qquad\quad+~\left[E_t(\pi(Q_2))(\pi(Q_2)-Q_2) E_t(\pi(Q_2))^{\prime}-E_t(\pi(Q_1))(\pi(Q_1)-Q_1) E_t(\pi(Q_1))^{\prime}\right].
\end{array}
$$
Using (\ref{expo-E}) we find that
$$
\Vert \phi_t(\pi(Q_1))-\phi_t(\pi(Q_2))\Vert_2\leq \kappa_{\phi}(\pi(Q_2),\pi(Q_1))~e^{-2\nu t}~\Vert \pi(Q_2)-\pi(Q_1)\Vert_2.
$$

To estimate the second term, we use the decomposition
$$
\begin{array}{l}
E_t(\pi(Q_2))~(\pi(Q_2)-Q_2)~
E_t(\pi(Q_2))^{\prime}-E_t(\pi(Q_1))~(\pi(Q_1)-Q_1)~
E_t(\pi(Q_1))^{\prime}\\
\\
\qquad\qquad\qquad=\left[E_t(\pi(Q_2))-E_t(\pi(Q_1))\right]~(\pi(Q_2)-Q_2)~
E_t(\pi(Q_2))^{\prime}\\
\\
\qquad\qquad\qquad\qquad+
E_t(\pi(Q_1))\left[\left\{(\pi(Q_2)-\pi(Q_1))-(Q_2-Q_1)\right\}
E_t(\pi(Q_2))^{\prime}\right]\\
\\
\qquad\qquad\qquad\qquad\qquad+E_t(\pi(Q_1))\left[
(\pi(Q_1)-Q_1)~\left[E_t(\pi(Q_2))^{\prime}-
E_t(\pi(Q_1))^{\prime}\right]\right].
\end{array}
$$
Combining  (\ref{expo-E}) with (\ref{expo-Elip}) we find that
$$
\begin{array}{l}
\Vert E_t(\pi(Q_2))~(\pi(Q_2)-Q_2)~
E_t(\pi(Q_2))^{\prime}-E_t(\pi(Q_1))~(\pi(Q_1)-Q_1)~
E_t(\pi(Q_1))^{\prime}\Vert_2\\
\\
\qquad\qquad\qquad\leq\kappa_{E}(\pi(Q_2),\pi(Q_1))~e^{-3\nu t}~\\
\\\hskip3cm\times
~\Vert \pi(Q_2)-\pi(Q_1)\Vert_2~\left[\kappa_E(\pi(Q_2))\Vert \pi(Q_2)-Q_2\Vert_2+\kappa_E(\pi(Q_1))\Vert \pi(Q_1)-Q_1
\Vert_2\right]\\
\\
\qquad\qquad\qquad\qquad+~\kappa_E(Q_1)\kappa_E(Q_2)~e^{-4\nu t}\left[\Vert \pi(Q_2)-\pi(Q_1)\Vert_2~+\Vert Q_2-Q_1\Vert_2~\right].
\end{array}
$$
To prove  (\ref{expo-phi-pi-bis}) we recall from Theorem~\ref{theo-pi-1} that $\pi\circ\phi^{\pi}_t=\phi^{\pi}_t\circ\pi$. This implies that
$$
\begin{array}{l}
\Vert \pi\left[\phi_t^{\pi}(Q_2)\right]-\phi_t^{\pi}(Q_1)\Vert_2\\
\\
\qquad\qquad\qquad \leq \kappa_{\phi^{\pi}}(Q_1,\pi(Q_2))~e^{-2\nu t}~\left[\Vert \pi(Q_2)-\pi(Q_1)\Vert_2+e^{-2\nu t}~\Vert \pi(Q_2)-Q_1\Vert_2\right].
\end{array}$$
If we set $Q_1=Q_2$ we obtain (\ref{expo-phi-pi-bis}). This ends the proof of the theorem.
\qed

Combining (\ref{expo-phi-pi}) with the fact that $\phi^{\pi}_t\circ\pi=\phi_t\circ\pi$  and $\pi(P^{\pi})=P$ we 
readily prove the 
following estimate.

\begin{cor}\label{expo-phi-pi-phi-1-cor}
For any $Q\in\SS^{+}_{r}$ and for any $t\geq 0$ we have
\begin{equation}\label{expo-phi-pi-phi}
\Vert \phi_t^{\pi}(Q)-\phi_t(\pi(Q))\Vert_2
\leq \kappa_{\phi^{\pi}}(Q,\pi(Q))~e^{-4\nu t}~\Vert Q-\pi(Q)\Vert_2~\Longrightarrow~P^{\pi}=P.
\end{equation}
 In addition, we have
\begin{equation}\label{expo-phi-pi-phi-1}
\Vert \phi_t^{\pi}(Q)-\phi_t(Q)\Vert_2\leq e^{-2\nu t}~\left[\kappa_{\phi}(Q,\pi(Q))+
 \kappa_{\phi^{\pi}}(Q,\pi(Q))~e^{-2\nu t}~\right]\Vert Q-\pi(Q)\Vert_2.
\end{equation}
 \end{cor}
 The estimate (\ref{expo-phi-pi-phi-1}) is a direct consequence of (\ref{expo-phi}) and (\ref{expo-phi-pi-phi}).
 Replacing $Q$ by $P=P^{\pi}$ in (\ref{expo-phi-pi-phi}) we obtain the following exponential decays to equilibrium.

\begin{cor}\label{cor-inter-1}
For any $Q\in\SS^{+}_{r}$ and for any $t\geq 0$ we have
 the local Lipschitz
estimate
 \begin{equation}\label{expo-phi-pi-2}
\Vert \phi_t^{\pi}(Q)-P\Vert_2
\leq \kappa_{\phi^{\pi}}(Q,P)~e^{-2\nu t}~\left[\Vert \pi(Q)-P\Vert_2+e^{-2\nu t}~\Vert Q-P\Vert_2\right].
\end{equation}
This yields the uniform estimate
$$
\Vert \phi^{\pi}(Q)\Vert_2:=
\sup_{t\geq 0}{\Vert \phi_t^{\pi}(Q)\Vert_2}\leq \Vert P\Vert_2+ \kappa_{\phi^{\pi}}(Q,P)~\left[\Vert \pi(Q)-P\Vert_2+\Vert Q-P\Vert_2\right].
$$
 \end{cor}

Combining (\ref{expo-Elip}) with Corollary~\ref{cor-inter-1} we prove the following local Lipschitz contraction.

\begin{cor}
For any $0\leq s\leq t$ and any $Q\in\SS_{r}^+$
$$
\Vert E_{t|s}(\phi_s^{\pi}(Q))-E_{t\vert s}(P)\Vert_2\leq \kappa_{E,\pi}(Q,P)~e^{-\nu (t+s)}~\left[\Vert \pi(Q)-P\Vert_2+e^{-2\nu s}~\Vert Q-P\Vert_2\right]
$$
 some finite constant  $\kappa_{E,\pi}(Q,P)<\infty$ whose values only depend respectively on 
$(\Vert P\Vert_2,\Vert Q\Vert_2)$. 
\end{cor}

\section{Kalman-Bucy stochastic flows}

\subsection{Perturbation-type models} \label{first-second-sec-KF}

We consider the perturbation models discussed in Section~\ref{perturbation-models-sec}. We set
$$
\sigma^2_{\delta}(Q) := 2\sqrt{2}~\kappa_{\delta,E}(Q)~\left[~r~
\left[\Vert R\Vert_2+\Vert S\Vert_2~\left(\delta+\Vert\phi(Q)\Vert_{\delta,2}\right)^2 \right]/((1-\delta)\nu)\right]^{1/2}
$$
with
$$
\kappa_{\delta,E}(Q):=\kappa_E(\Vert Q\Vert_2)~\exp{\left[\chi_2(\delta,\Vert Q\Vert_2)/(4\nu)\right]}\quad
\mbox{\rm and}\quad
\Vert\phi(Q)\Vert_{\delta,2}:=
\sup_{t\geq 0}{\sup_{\pi\in B(\delta)}{
\Vert\phi^{\pi}_{t}(Q)\Vert_2}}<\infty
$$
where $\chi_2(\cdot,\cdot)<\infty$ is introduced in Theorem \ref{theo-continuity} and $\nu>0$ and $\kappa_E(\cdot)<\infty$ are defined in (\ref{expo-E}). 

Recall the semigroup and stochastic flow notation defined in Section~\ref{statement-sec}. The first lemma in this section concerns the convergence of the perturbed Kalman-Bucy filter to the true underlying signal process, both in a mean-square sense and in terms of actual sample paths.

\begin{lem}\label{control-lemma}
Assume $\mbox{\rm (H)}_0$ and $\mbox{\rm (H)}_2$ are satisfied. For any $\epsilon>0$ there exists some parameter $0<\delta<\epsilon$ such that for any $\pi\in B(\delta)$, $0\leq s\leq t$, $Q\in \SS_{r}^+$ and any $n\geq 1$ we have
$$
\EE\left[\Vert \psi_{s,t}^{\pi}(x,Q)- \theta_{s,t}(X_s) \Vert_2^{2n}\vert X_s\right]^{\frac{1}{2n}}\leq \kappa_{\delta,E}(Q)~e^{-2(1-\epsilon)\nu(t-s)}~\Vert X_s-x\Vert_2+\sqrt{n}~\sigma_{\delta}(Q). 
$$
In addition, the conditional probability of the following event
$$
\Vert \psi_{s,t}^{\pi}(x,Q)- \theta_{s,t}(X_s) \Vert_2\leq 
  \kappa_{\delta,E}(Q)~e^{-2(1-\epsilon)\nu(t-s)}~\Vert X_s-x\Vert_2+~\sigma_{\delta}(Q)~\frac{e^2}{\sqrt{2}}~\left[\frac{1}{2}+
  \left(z+\sqrt{z}\right)\right]
$$
given the state variable $X_s$ is greater than $1-e^{-z}$.
\end{lem}
\proof
We have
$$
\begin{array}{l}
	d\left[\psi_{s,t}^{\pi}(x,Q)- \theta_{s,t}(X_s) \right] =\left[A-\pi\left(\phi^{\pi}_{s,t}(Q)\right)S\right]~\left[\psi^{\pi}_{s,t}(x,Q)- \theta_{s,t}(X_s) \right]~dt +dM_{s,t}^{\pi}
\end{array}
$$
In the above display, $t\in [s,\infty[\mapsto M_{s,t}^{\pi}$ stands for the $r$-dimensional martingale
$$
dM_{s,t}^{\pi}:=R^{1/2}dW_t+\pi\left(\phi_{s,t}^{\pi}(Q)\right)~C^{\prime}R^{-1/2}_{2}~dV_t
$$
with angle bracket
$$
\left(\partial_t\langle M_{s,t}^{\pi}(k),M_{s,t}^{\pi}(l) \rangle\right)_{1\leq k,l\leq r}:=R+\pi\left(\phi_{s,t}^{\pi}(Q)\right)S\pi\left(\phi_{s,t}^{\pi}(Q)\right).
$$
This yields the formula
$$
\psi_{s,t}^{\pi}(x,Q)- \theta_{s,t}(X_s) -E^{\pi}_{t|s}(Q)(x-X_s)=\int_s^tE^{\pi}_{u,t|s}(Q)~dM_{s,u}^{\pi}
$$
with the exponential semigroup $E^{\pi}_{u,t|s}(Q)$ defined for any $s\leq u\leq t$ by
$$
E^{\pi}_{u,t|s}(Q)=\exp{\left(\oint_u^t\left[A-\pi\left(\phi^{\pi}_{s,u}(Q)\right)S\right]~du\right)}.
$$
We have the decomposition
$$
A-\pi\left(\phi^{\pi}_{s,u}(Q)\right)S=A-\phi_{s,u}(Q)S+\left[\phi_{s,u}(Q)-\phi^{\pi}_{s,u}(Q)+\phi^{\pi}_{s,u}(Q)-\pi\left(\phi^{\pi}_{s,u}(Q)\right)\right]S.
$$
By (\ref{expo-E}) we have
$$
\Vert E_{t|s}(Q_1)\Vert_2\leq \kappa_E(\Vert Q_1\Vert_2)~e^{-2\nu (t-s)}.
$$
By Theorem~\ref{theo-continuity} there exists some $\delta>0$ such that for any $\pi\in B(\delta)$
$$
\Vert \pi\left[\phi^{\pi}_{t}(Q)\right]-\phi^{\pi}_{t}(Q)\Vert_2\leq \delta~
\Longrightarrow~\sup_{t\geq 0}{\sup_{\pi\in B(\delta)}{
\Vert \pi\left[\phi^{\pi}_{t}(Q)\right]\Vert_2}}\leq \delta+\Vert\phi(Q)\Vert_{\delta,2}.
$$
In addition, we have
$$
\Vert \phi_{s,t}(Q)-\phi^{\pi}_{s,u}(Q)\Vert_2\leq  \left[ \chi_1(\delta)+ e^{-4(t-s)\nu}~\chi_2(\delta,\Vert Q\Vert_2)\right]\gamma(\pi).
$$
Applying Lemma~\ref{perturbation-lemma-intro} we find that
$$
\Vert E^{\pi}_{u,t|s}(Q)\Vert_2\leq \kappa_{\delta,E}(Q)~\exp{\left[-\left\{2\nu-\delta-\gamma(\pi)~\kappa_E(\Vert Q\Vert_2)
  \chi_1(\delta)\right\}~ (t-s)\right]}.
$$
For any $\epsilon>0$ we choose $\epsilon>\delta>\delta^{\prime}>0$ and so that for any $\pi\in B(\delta^{\prime})$
$$
\delta+\gamma(\pi)~\kappa_E(\Vert Q\Vert_2)
  \chi_1(\delta)\leq 2\epsilon \nu\Rightarrow \Vert E^{\pi}_{u,t|s}(Q)\Vert_2\leq \kappa_{\delta,E}(Q)~\exp{\left[-2(1-\epsilon)\nu(t-s)\right]}.
$$ 

Following the proof of Lemma 5.3 in~\cite{ap-2016}, for any $n\geq 1$ we have
$$
\begin{array}{l}
\displaystyle\EE\left[\left(\Vert\int_s^tE^{\pi}_{u,t|s}(Q)~dM_{s,u}^{\pi}\Vert_2^{2n}\right)\right]^{\frac{1}{n}}\\
\\
\qquad\qquad\qquad \displaystyle\leq 4^2n~r\int_s^t
\Vert R+\pi\left(\phi_{s,u}^{\pi}(Q)\right)S\pi\left(\phi_{s,u}^{\pi}(Q)\right)\Vert_2~
\Vert E^{\pi}_{u,t|s}(Q)\Vert_2~du\\
\\
\qquad\qquad\qquad\displaystyle\leq 8n~r~\kappa_{\delta,E}(Q)
\left[\Vert R\Vert_2+\Vert S\Vert_2~\left(\delta+\Vert\phi(Q)\Vert_{\delta,2}\right)^2 \right]/((1-\epsilon)\nu)\leq n~\sigma^2_{\delta}(Q).
\end{array}
$$
The end of the proof of the first assertion is now easily completed. The proof of the exponential concentration inequality follows the same line of argument as the proof of Theorem 5.2 in~\cite{ap-2016}. This ends the proof of the lemma.
\qed

The next three theorems concern convergence of the $\pi$-perturbed Kalman-Bucy filter/diffusion to the true, optimal, Kalman-Bucy filter \cite{ap-2016}. The first concerns the stochastic flow of the two filters themselves, while the second two theorems concern the associated (conditional) distributions.

\begin{theo}\label{theo-KB-stoch-flows}
Assume $\mbox{\rm (H)}_0$ and $\mbox{\rm (H)}_2$ are satisfied. In this situation, there exists some parameter $\delta>0$ such that
for any $0 <\epsilon<\delta$, $\pi\in B(\epsilon)$, $0\leq s\leq t$, $Q\in\SS_{r}^+$, and any $n\geq 1$
we have
$$
\EE\left[
\Vert \psi_{s,t}^{\pi}(x,Q)-\psi_{s,t}(x,Q)\Vert_{2}^{2n}~|~X_s
\right]^{\frac{1}{2n}}
\leq \epsilon~{\chi}(\delta,Q)~\left[\sqrt{n}+e^{-(1-\epsilon)\nu(t-s)}~\Vert X_s-x\Vert_2\right]
$$
for some finite constants ${\chi}(\delta,Q)$ whose value only depends on the parameters $(\delta,\Vert Q\Vert_2)$.
\end{theo}
\proof
We have
$$
\begin{array}{l}
d\left[\psi_{s,t}^{\pi}(x,Q)-\psi_{s,t}(x,Q)\right]\\
\\
\qquad\quad=\left\{\left[A-\pi\left(\phi^{\pi}_{s,t}(Q)\right)S\right]-\left[A-\phi_{s,t}(Q)S\right]\right\}~\psi^{\pi}_{s,t}(x,Q)\\
\\
\qquad\qquad\quad-
\left[A-\phi_{s,t}(Q)S\right]~\left[\psi_{s,t}(x,Q)-\psi^{\pi}_{s,t}(x,Q)\right]~dt
+\left[\pi\left(\phi_{s,t}^{\pi}(Q)\right)-\phi_{s,t}(Q)\right]~C^{\prime}\Sigma^{-1}~dY_t.
\end{array}
$$
This implies that
$$
\begin{array}{l}
d\left[\psi_{s,t}^{\pi}(x,Q)-\psi_{s,t}(x,Q)\right]\\
\\
\qquad\qquad=\left[\phi_{s,t}(Q)-\pi\left(\phi^{\pi}_{s,t}(Q)\right)\right]S~\psi^{\pi}_{s,t}(x,Q)~dt\\
\\
\hskip4cm+\left[\pi\left(\phi_{s,t}^{\pi}(Q)\right)-\phi_{s,t}(Q)\right]~C^{\prime}\Sigma^{-1}~\left(C\theta_{s,t}(X_s)dt+R^{1/2}_2dV_t\right)\\
\\
\hskip7cm+
\left[A-\phi_{s,t}(Q)S\right]~\left[\psi^{\pi}_{s,t}(x,Q)-\psi_{s,t}(x,Q)\right]~dt\\
\\
\qquad\qquad=\left[A-\phi_{s,t}(Q)S\right]~\left[\psi^{\pi}_{s,t}(x,Q)-\psi_{s,t}(x,Q)\right]~dt\\
\\
\hskip5cm+\left[\phi_{s,t}(Q)-\pi\left(\phi^{\pi}_{s,t}(Q)\right)\right]S~\left[\psi^{\pi}_{s,t}(x,Q)-\theta_{s,t}(X_s)\right]~dt+dM^{\pi}_{s,t}
\end{array}
$$
with the $r$-dimensional martingale $t\in [s,\infty[\mapsto M^{\pi}_{s,t}$ defined by
$$
dM^{\pi}_{s,t}=\left[\pi\left(\phi_{s,t}^{\pi}(Q)\right)-\phi_{s,t}(Q)\right]~C^{\prime}R^{-1/2}_{2}dV_t.
$$
This implies that
$$
\begin{array}{l}
\psi_{s,t}^{\pi}(x,Q)-\psi_{s,t}(x,Q)\qquad \\
\\
\qquad\displaystyle=\int_s^tE_{u,t|s}(Q)\left[\phi_{s,u}(Q)-\pi\left(\phi^{\pi}_{s,u}(Q)\right)\right]S\left[\psi^{\pi}_{s,u}(x,Q)-\theta_{s,u}(X_s)\right]\,du\displaystyle+\int_s^tE_{u,t|s}(Q)\,dM^{\pi}_{s,u}.
\end{array}
$$
Arguing as in the proof of Lemma~\ref{control-lemma}, there exists some $0<\epsilon<\delta$ such that for any $\pi\in B(\epsilon)$
$$
\Vert \pi\left[\phi^{\pi}_{t}(Q)\right]-\phi^{\pi}_{t}(Q)\Vert_2\leq \epsilon~
$$
and
$$
\Vert \phi_{s,u}(Q)-\phi^{\pi}_{s,u}(Q)\Vert_2\leq  \epsilon~\left[ \chi_1(\delta)+ e^{-4(u-s)\nu}~\chi_2(\delta,\Vert Q\Vert_2)\right].
$$
By the generalized Minkoswki iequality, we have
$$
\begin{array}{l}
\displaystyle\Vert \int_s^t\,E_{u,t|s}(Q)~\left[\phi_{s,u}(Q)-\pi\left(\phi^{\pi}_{s,u}(Q)\right)\right]S~\left[\psi^{\pi}_{s,u}(x,Q)-\theta_{s,u}(X_s)\right]\,du\Vert_{2}
\displaystyle~\leq~\epsilon~\kappa_E(Q)~ \Vert S \Vert_{2} \\
 \\
\qquad\qquad\quad \displaystyle\times~ \int_s^t\left[ e^{-(t-u)\nu}(1+\chi_1(\delta))+ e^{-4(t-s)\nu}\chi_2(\delta,\Vert Q\Vert_2)\right]\Vert \psi^{\pi}_{s,u}(x,Q)-\theta_{s,u}(X_s) \Vert_{2}~du.
\end{array}
$$
This implies that
$$
\begin{array}{l}
\displaystyle\EE\left[\Vert \int_s^t~E_{u,t|s}(Q)~\left[\phi_{s,u}(Q)-\pi\left(\phi^{\pi}_{s,u}(Q)\right)\right]S~\left[\psi^{\pi}_{s,u}(x,Q)-\theta_{s,u}(X_s)\right]~du\Vert_{2}^{2n}~\vert~X_s\right]^{\frac{1}{2n} \qquad}\\
\\
 \qquad\qquad\displaystyle\leq\epsilon~\kappa_E(Q)~ \Vert S \Vert_{2}~ \int_s^t~ ~\left[ e^{-(t-u)\nu}(1+\chi_1(\delta))+ e^{-4(t-s)\nu}~\chi_2(\delta,\Vert Q\Vert_2)\right]~\\
 \\
\hskip7cm \displaystyle\times~
\EE\left[\Vert \psi^{\pi}_{s,u}(x,Q)-\theta_{s,u}(X_s) \Vert_{2}^{2n}~|~X_s\right]^{\frac{1}{2n}}
~du.
\end{array}
$$
By Lemma~\ref{control-lemma} we have
$$
\begin{array}{l}
\displaystyle\EE\left[\Vert \int_s^t~E_{u,t|s}(Q)~\left[\phi_{s,u}(Q)-\pi\left(\phi^{\pi}_{s,u}(Q)\right)\right]S~\left[\psi^{\pi}_{s,u}(x,Q)-\theta_{s,u}(X_s)\right]~du\Vert_{2}^{2n}~\vert~X_s\right]^{\frac{1}{2n}}\qquad\qquad\\
\\
 \qquad\qquad\qquad\qquad\displaystyle\leq~\epsilon~\kappa_E(Q)~ \Vert S \Vert_{2}~ \int_s^t~ ~\left[ e^{-(t-u)\nu}(1+\chi_1(\delta))+ e^{-4(t-s)\nu}~\chi_2(\delta,\Vert Q\Vert_2)\right]~\\
 \\
\qquad\qquad\qquad\qquad\qquad\qquad\qquad\qquad \displaystyle\times~\left[
\kappa_{\delta,E}(Q)~e^{-2(1-\epsilon)\nu(u-s)}~\Vert X_s-x\Vert_2+\sqrt{n}~\sigma_{\delta}(Q) \right]
~du.
\end{array}
$$
This yields the estimate
$$
\begin{array}{l}
\displaystyle\EE\left[\Vert \int_s^t~E_{u,t|s}(Q)~\left[\phi_{s,u}(Q)-\pi\left(\phi^{\pi}_{s,u}(Q)\right)\right]S~\left[\psi^{\pi}_{s,u}(x,Q)-\theta_{s,u}(X_s)\right]~du\Vert_{2}^{2n}~\vert~X_s\right]^{\frac{1}{2n}}\qquad\\
\\
 \qquad\qquad\qquad\qquad\displaystyle\leq~\epsilon~\kappa_E(Q)~ \Vert S \Vert_{2}~\left\{\sqrt{n}~ \left[ \overline{\chi}_1(\delta,Q)+ e^{-3(t-s)\nu}~\overline{\chi}_2(\delta,Q)\right]\right.\\
\\
\qquad\qquad\qquad\left. \displaystyle\hskip3cm+~\Vert X_s-x\Vert_2~e^{-(1-\epsilon)\nu(t-s)}~~\left[
\underline{\chi}_1(\delta,Q)+e^{-3(t-s)\nu}~\underline{\chi}_2(\delta,Q)\right]~
\right\}
 \end{array}
$$
with
$$
\overline{\chi}_1(\delta,Q):=\sigma_{\delta}(Q) ~(1+\chi_1(\delta))/\nu\quad
\mbox{\rm and}\quad
\overline{\chi}_2(\delta,Q):=\sigma_{\delta}(Q) \chi_2(\delta,\Vert Q\Vert_2)
$$
and
$$
\underline{\chi}_1(\delta,Q):=\kappa_{\delta,E}(Q)~(1+\chi_1(\delta))/((1-\delta)\nu)\quad
\mbox{\rm and}\quad
\underline{\chi}_2(\delta,Q):=\kappa_{\delta,E}(Q) \chi_2(\delta,\Vert Q\Vert_2)/(2(1-\delta)\nu).
$$
Following the proof of Lemma 5.3 in~\cite{ap-2016}, for any $n\geq 1$ we have
$$
\begin{array}{l}
\displaystyle\EE\left[\left(\Vert\int_s^tE_{u,t|s}(Q)~dM_{s,u}^{\pi}\Vert_2^{2n}\right)\right]^{\frac{1}{n}}\qquad\\
\\
\qquad\qquad\displaystyle\leq 4^2n~r\int_s^t
\Vert \left[\phi_{s,u}(Q)-\pi\left(\phi^{\pi}_{s,u}(Q)\right)\right]S\left[\phi_{s,u}(Q)-\pi\left(\phi^{\pi}_{s,u}(Q)\right)\right]\Vert_2~
\Vert E_{u,t|s}(Q)\Vert_2~du\\
\\
\qquad\qquad\displaystyle\leq 8~\epsilon^2~n~r~\kappa_{E}(Q)\Vert S\Vert_2/\nu~=~\epsilon^2~\overline{\sigma}^2(Q)
\end{array}
$$
with
$$
\overline{\sigma}^2(Q):=8\,r\,n\,\kappa_{E}(Q)\Vert S\Vert_2/\nu.
$$
This yields
$$
\begin{array}{l}
\epsilon^{-1}~\EE\left[
\Vert \psi_{s,t}^{\pi}(x,Q)-\psi_{s,t}(x,Q)\Vert_{2}^{2n}
\right]^{\frac{1}{2n}}\qquad\\
\\
\qquad\qquad\leq 
~\sqrt{n}~\overline{\sigma}(Q)
+\kappa_E(Q)~ \Vert S \Vert_{2}~
\left\{
\sqrt{n}~ \left[ \overline{\chi}_1(\delta,Q)+ e^{-3(t-s)\nu}~\overline{\chi}_2(\delta,Q)\right]\right.\\
\\
\left. \displaystyle\hskip3cm+~\Vert X_s-x\Vert_2~e^{-(1-\epsilon)\nu(t-s)}~~
\left[
\underline{\chi}_1(\delta,Q)+e^{-3(t-s)\nu}~\underline{\chi}_2(\delta,Q)
\right]~
\right\}\\
\\
\qquad\qquad\leq \sqrt{n}~\left[\overline{\sigma}(Q)+\kappa_E(Q)~ \Vert S \Vert_{2}~ \left[ \overline{\chi}_1(\delta,Q)+ e^{-3(t-s)\nu}~\overline{\chi}_2(\delta,Q)\right]\right]+\\
\\
\displaystyle\hskip3cm+~e^{-(1-\epsilon)\nu(t-s)}~\Vert X_s-x\Vert_2~~
\left[
\underline{\chi}_1(\delta,Q)+e^{-3(t-s)\nu}~\underline{\chi}_2(\delta,Q)
\right].
\end{array}
$$
This ends the proof of the theorem.\qed

\begin{theo}\label{theo-entropy}
Assume $\mbox{\rm (H)}_0$ and $\mbox{\rm (H)}_2$  are satisfied. In this situation, for any $s+\upsilon\leq t$, and any $Q\in\SS_{r}^+$ we have
$$
 \begin{array}{l}
 \mbox{\rm Ent}\left(\eta_{s,t}^{\pi}(x,Q)~|~\eta_{s,t}(x,Q)\right)\\
 \\
\qquad\qquad\qquad \leq  \displaystyle\frac{1}{2} \left(\varpi^{o}_{+}(\Ca)+ 1/\varpi^{c}_{-}\right)~
 \left[\left\Vert{\psi}^{\pi}_{s,t}(x,Q)-{\psi}_{s,t}(x,Q)\right\Vert_2^2+ \frac{5}{2}~\sqrt{r}~\Vert \phi_{s,t}(Q)-\phi^{\pi}_{s,t}(Q)\Vert_2\right].
\end{array}
$$

\end{theo}

\proof
The Boltzmann  relative entropy of $\eta_{s,t}^{\pi}(x,Q)$ w.r.t. $\eta_{s,t}(x,Q)$  is given by the formula
 $$
 \begin{array}{l}
 \mbox{\rm Ent}\left(\eta_{s,t}^{\pi}(x,Q)~|~\eta_{s,t}(x,Q)\right) \displaystyle=-\frac{1}{2}\left[
 \mbox{\rm tr}\left(I-\phi_{s,t}(Q)^{-1}\phi^{\pi}_{s,t}(Q)\right)+\log{\mbox{det}\left(\phi^{\pi}_{s,t}(Q)\phi_{s,t}(Q)^{-1}\right)}\right]\\
 \\
 \hskip4cm\displaystyle+\frac{1}{2}\left\langle\left({\psi}^{\pi}_{s,t}(x,Q)-{\psi}_{s,t}(x,Q)\right),\phi_{s,t}(Q)^{-1}\left({\psi}^{\pi}_{s,t}(x,Q)-{\psi}_{s,t}(x,Q)\right)\right\rangle.
\end{array}
$$
By Corollary~\ref{cor-H2}, for any $t\geq s+\upsilon$ we have
 $$
 \begin{array}{l}
 \mbox{\rm Ent}\left(\eta_{s,t}^{\pi}(x,Q)~|~\eta_{s,t}(x,Q)\right) \displaystyle=-\frac{1}{2}\left[
 \mbox{\rm tr}\left(I-\phi_{s,t}(Q)^{-1}\phi^{\pi}_{s,t}(Q)\right)+\log{\mbox{det}\left(\phi^{\pi}_{s,t}(Q)\phi_{s,t}(Q)^{-1}\right)}\right]\\
 \\
 \hskip6cm\displaystyle+\frac{1}{2} \left(\varpi^{o}_{+}(\Ca)+ 1/\varpi^{c}_{-}\right)~
 \left\Vert{\psi}^{\pi}_{s,t}(x,Q)-{\psi}_{s,t}(x,Q)\right\Vert_2^2.
\end{array}
$$
In addition, there exists some $\delta>0$ s.t. for any $\pi\in B(\delta)$ 
 \begin{eqnarray*}
\Vert I-\phi_{s,t}(Q)^{-1} \phi^{\pi}_{s,t}(Q)\Vert_2&=&\Vert ( \phi_{s,t}(Q)- \phi^{\pi}_{s,t}(Q)) \phi_{s,t}(Q)^{-1}\Vert_2\\
&\leq&  \left(\varpi^{o}_{+}(\Ca)+ 1/\varpi^{c}_{-}\right)~\Vert \phi_{s,t}(Q)-\phi^{\pi}_{s,t}(Q)\Vert_2 ~\leq~ \frac{1}{2\sqrt{r}}.
 \end{eqnarray*}
 This implies that
  $$
 \begin{array}{l}
 \mbox{\rm Ent}\left(\eta_{s,t}^{\pi}(x,Q)~|~\eta_{s,t}(x,Q)\right) \displaystyle=\frac{1}{2}
 \mbox{\rm tr}\left(\phi_{s,t}(Q)^{-1}\left[\phi^{\pi}_{s,t}(Q)-\phi_{s,t}(Q)\right]~\right)\\
 \\
 \qquad\qquad\qquad\quad\displaystyle+\frac{1}{2} \left(\varpi^{o}_{+}(\Ca)+ 1/\varpi^{c}_{-}\right)~
 \left[\left\Vert{\psi}^{\pi}_{s,t}(x,Q)-{\psi}_{s,t}(x,Q)\right\Vert_2^2+ \frac{3}{2}~\sqrt{r}~\Vert \phi_{s,t}(Q)-\phi^{\pi}_{s,t}(Q)\Vert_2\right].
\end{array}
$$
The last assertion is a consequence of the following lemma applied to $A=I-\phi_{s,t}(Q)^{-1} \phi^{\pi}_{s,t}(Q)$.

\begin{lem}\label{lem-tech-2}
For any $(r\times r)$-matrix $A$  we have
$$
\Vert A \Vert_2< \frac{1}{2\sqrt{r}}\Longrightarrow
\left\vert\log{\mbox{\rm det}\left(I-A\right)}\right\vert 
\leq \frac{3}{2}~\sqrt{r}~\Vert A \Vert_2.
$$

\end{lem}

\proof
For any $n\geq 1$ we have
$$
\vert\mbox{\rm tr}(A^n)\vert\leq \Vert A \Vert^n_F\leq \sqrt{r}^n~\Vert A \Vert^2_n.
$$
Using the well-known trace formulae $$
\log{\mbox{\rm det}(I-A)}=\mbox{\rm tr}(\log{(I-A)})=-\sum_{n\geq 1}~n^{-1}~\mbox{\rm tr}(A^n)
$$
 we conclude that
 $$
 \vert \log{\mbox{\rm det}(I-A)}\vert\leq -\log{\left( 1-\sqrt{r}~\Vert A \Vert_2\right)}.
 $$
The last assertion comes from the inequality
$$
0\leq -\log{(1-u)}\leq u+\frac{1}{2}~\frac{u^2}{1-u}=u\left(1+\frac{1}{2}~\frac{u}{1-u}\right)\leq 3u/2
$$
which is valid for any $u\in [0,1/2[$. This ends the proof of the lemma.\qed

To take the final step in the proof of the theorem we note that $\phi^{\pi}_{s,t}(Q)\geq \phi_{s,t}(Q)$ implies
$$
\begin{array}{l}
\begin{array}[t]{rcl}
 \mbox{\rm tr}\left(\phi_{s,t}(Q)^{-1}\left[\phi^{\pi}_{s,t}(Q)-\phi_{s,t}(Q)\right]~\right)&\leq&
  \left(\varpi^{o}_{+}(\Ca)+ 1/\varpi^{c}_{-}\right)~ \mbox{\rm tr}\left(\left[\phi^{\pi}_{s,t}(Q)-\phi_{s,t}(Q)\right]\right)\\
  &&\\
  &\leq &  \sqrt{r}~\left(\varpi^{o}_{+}(\Ca)+ 1/\varpi^{c}_{-}\right)~ \Vert\phi^{\pi}_{s,t}(Q)-\phi_{s,t}(Q)\Vert_2
\end{array}\end{array}.
$$
This ends the proof of the theorem.\qed

\begin{theo}\label{theo-wass}
Assume $\mbox{\rm (H)}_0$ and $\mbox{\rm (H)}_2$  are satisfied. For any $Q\in\SS_{r}^+$, and for $t\geq s+\upsilon$ we have the almost sure Wasserstein estimate
$$
\begin{array}{l}
\WW_2\left[\eta_{s,t}^{\pi}(x,Q),\eta_{s,t}(x,Q)\right]^2\\
\\
\qquad\qquad\leq \Vert{\psi}^{\pi}_{s,t}(x,Q)-{\psi}_{s,t}(x,Q)\Vert_2^2+\tr\left[
\phi^{\pi}_{s,t}(Q)-\phi_{s,t}(Q)\right]\\
\\
\qquad\qquad\qquad\qquad +\,4 r~\left[~\varpi^{c}_{+}(\Oa)+1/\varpi^{o}_{-}\right]\left[\varpi^{o}_{+}(\Ca)+ 1/\varpi^{c}_{-}\right] ~\Vert 
\phi_{s,t}^{\pi}(Q)-\phi_{s,t}(Q)
\Vert_2.
\end{array}
$$
In addition,
 for any $n\geq 1$ and any $t\geq s+\upsilon$ we have
 $$
\begin{array}{l}
\WW_{2n}\left[\eta_{s,t}(x_1,Q_1),\eta^{\pi}_{s,t}(x_2,Q_2)\right]\\
\\
\qquad\qquad\displaystyle\leq \Vert \psi_{s,t}(x,Q)-\psi^{\pi}_{s,t}(x,Q)\Vert\\
\\
\displaystyle\hskip3cm+\sqrt{\frac{rn}{2}}~
 \left(\varpi^{o}_{+}(\Ca)+ 1/\varpi^{c}_{-}\right)^{1/2}
\Vert {\phi}_{t-s}(x_1,Q_1)-{\phi}^{\pi}_{t-s}(x_2,Q_2)\Vert_2~
e^{1/2+\frac{3}{4n}}.
\end{array}
$$
\end{theo}

\proof
The  $\LL_2$-Wasserstein distance between the Gaussian distributions $\eta_{s,t}^{\pi}(x,Q)$, and $\eta_{s,t}(x,Q)$ is given by
$$
\begin{array}{l}
\WW_2\left[\eta_{s,t}^{\pi}(x,Q),\eta_{s,t}(x,Q)\right]^2\\
\\
\qquad\quad=\Vert{\psi}^{\pi}_{s,t}(x,Q)-{\psi}_{s,t}(x,Q)\Vert_2^2+
\tr\left[
 \phi_{s,t}(Q)+\phi^{\pi}_{s,t}(Q)-
2\left[
\phi_{s,t}(Q)^{1/2} \phi_{s,t}^{\pi}(Q) \phi_{s,t}(Q)^{1/2}\right]^{1/2}
\right].
\end{array}
$$
A proof of this formula can be found in \cite{givens,olkin}. We assume that 
$$\left[
\phi_{s,t}(Q)^{1/2} \phi_{s,t}^{\pi}(Q) \phi_{s,t}(Q)^{1/2}\right]^{1/2}\geq 0$$ is the principal square root 
of the positive definite matrix $
\phi_{s,t}(Q)^{1/2} \phi_{s,t}^{\pi}(Q) \phi_{s,t}(Q)^{1/2}\geq 0$.
Also observe that
$$
\phi_{s,t}^{\pi}(Q)\geq \phi_{s,t}(Q)~\Rightarrow \begin{array}[t]{rcl}\phi_{s,t}(Q)^{1/2} \phi_{s,t}^{\pi}(Q) \phi_{s,t}(Q)^{1/2}&\geq&
\phi_{s,t}(Q)^2\\
&\geq& \lambda_{min}(\phi_{s,t}(Q))^2~Id\\
&\geq&  \left(\varpi^{o}_{+}(\Ca)+ 1/\varpi^{c}_{-}\right)^{-2}~Id
\end{array}
$$
as soon as $t\geq s+\upsilon$.
The last estimate is a consequence of Theorem~\ref{bplem}.

Observe that
$$
 \phi_{s,t}(Q)+\phi^{\pi}_{s,t}(Q)-
2\left[
\phi_{s,t}(Q)^{1/2} \phi_{s,t}(Q) \phi_{s,t}(Q)^{1/2}\right]^{1/2}=~\phi^{\pi}_{s,t}(Q)-\phi_{s,t}(Q)~\geq ~0.
$$
This implies that
$$
\begin{array}{l}
\WW_2\left[\eta_{s,t}^{\pi}(x,Q),\eta_{s,t}^{\pi}(x,Q)\right]^2\\
\\
\qquad=~\Vert{\psi}^{\pi}_{s,t}(x,Q)-{\psi}_{s,t}(x,Q)\Vert_2^2+\tr\left[
\phi^{\pi}_{s,t}(Q)-\phi_{s,t}(Q)\right]\\
\\
\hskip2cm+2~\tr\left[
\left[\phi_{s,t}(Q)^{1/2} \phi_{s,t}(Q) \phi_{s,t}(Q)^{1/2}\right]^{1/2}
-\left[\phi_{s,t}(Q)^{1/2} \phi_{s,t}^{\pi}(Q) \phi_{s,t}(Q)^{1/2}\right]^{1/2}
\right]\\
\\
\qquad\leq ~\Vert{\psi}^{\pi}_{s,t}(x,Q)-{\psi}_{s,t}(x,Q)\Vert_2^2+\tr\left[
\phi^{\pi}_{s,t}(Q)-\phi_{s,t}(Q)\right]\\
\\
\hskip2cm +2~ r~\Vert   \left[\phi_{s,t}(Q)^{1/2} \phi_{s,t}(Q) \phi_{s,t}(Q)^{1/2}\right]^{1/2}-\left[\phi_{s,t}(Q)^{1/2} \phi_{s,t}^{\pi}(Q) \phi_{s,t}(Q)^{1/2}\right]^{1/2}\Vert_2.
\end{array}
$$
Using  (\ref{square-root-key-estimate}) we have
$$
\begin{array}{l}
\Vert   \left[\phi_{s,t}(Q)^{1/2} \phi_{s,t}(Q) \phi_{s,t}(Q)^{1/2}\right]^{1/2}-\left[\phi_{s,t}(Q)^{1/2} \phi_{s,t}^{\pi}(Q) \phi_{s,t}(Q)^{1/2}\right]^{1/2}\Vert_2\\
\\
\qquad\qquad\qquad\qquad\displaystyle\leq 2~\left(\varpi^{o}_{+}(\Ca)+ 1/\varpi^{c}_{-}\right)
\Vert   \phi_{s,t}(Q)^{1/2} \left[
\phi_{s,t}^{\pi}(Q)-\phi_{s,t}(Q)\right] \phi_{s,t}(Q)^{1/2}
\Vert_2\\
\\
\qquad\qquad\qquad\qquad\displaystyle\leq 2~\left[\varpi^{o}_{+}(\Ca)+ 1/\varpi^{c}_{-}\right] 
\Vert  \phi_{s,t}(Q)^{1/2}
\Vert_2^2 ~ \Vert 
\phi_{s,t}^{\pi}(Q)-\phi_{s,t}(Q)
\Vert_2.
\end{array}
$$
By Corollary~\ref{cor-H2} we conclude that
$$
\begin{array}{l}
\Vert   \left[\phi_{s,t}(Q)^{1/2} \phi_{s,t}(Q) \phi_{s,t}(Q)^{1/2}\right]^{1/2}-\left[\phi_{s,t}(Q)^{1/2} \phi_{s,t}^{\pi}(Q) \phi_{s,t}(Q)^{1/2}\right]^{1/2}\Vert_2\\
\\
\qquad\qquad\qquad\qquad\displaystyle\leq 2\left[~\varpi^{c}_{+}(\Oa)+1/\varpi^{o}_{-}\right]\left[\varpi^{o}_{+}(\Ca)+ 1/\varpi^{c}_{-}\right] ~\Vert 
\phi_{s,t}^{\pi}(Q)-\phi_{s,t}(Q)
\Vert_2.
\end{array}
$$
This ends the proof of the first assertion.

Observe that
$$
\overline{\psi}_{s,t}(x,Q)\stackrel{law}{=}\psi_{s,t}(x,Q)+{\phi}_{t-s}(x,Q)^{1/2}~Z
$$
and
$$
\overline{\psi}^{\pi}_{s,t}(x,Q)\stackrel{law}{=}\psi^{\pi}_{s,t}(x,Q)+{\phi}^{\pi}_{t-s}(x,Q)^{1/2}~Z
$$
where $Z$ stands for an $r$-dimensional Gaussian random variable with unit covariance matrix, and  ${\phi}_{t-s}(x,Q)^{1/2}$
stands for the principal square root of ${\phi}_{t-s}(x,Q)$. Combining (\ref{square-root-key-estimate}) with Theorem~\ref{bplem} for any $n\geq 1$ and any $t\geq s+\upsilon$ we have
$$
\begin{array}{l}
\WW_{2n}\left[\eta_{s,t}(x_1,Q_1),\eta^{\pi}_{s,t}(x_2,Q_2)\right]\\
\\
\qquad\qquad\displaystyle\leq \Vert \psi_{s,t}(x,Q)-\psi^{\pi}_{s,t}(x,Q)\Vert\\
\\
\displaystyle\hskip3cm+\sqrt{\frac{rn}{2}}~
 \left(\varpi^{o}_{+}(\Ca)+ 1/\varpi^{c}_{-}\right)^{1/2}
\Vert {\phi}_{t-s}(x_1,Q_1)-{\phi}^{\pi}_{t-s}(x_2,Q_2)\Vert_2~
e^{1/2+\frac{3}{4n}}.
\end{array}
$$
To check the last assertion, we use Stirling approximation to prove that
\begin{eqnarray*}
\EE\left[\Vert \sum_{1\leq k\leq r}Z_k^2\Vert^{n}\right]^{\frac{1}{n}}~\leq~ \sum_{1\leq k\leq r}
\EE\left[Z_1^{2n}\right]^{\frac{1}{n}} ~=~ \frac{r}{2}~\left[\frac{(2n)!}{n!}\right]^{\frac{1}{n}} ~\leq~ 2~r~n~
e^{1+\frac{3}{2n}}.
\end{eqnarray*}
This ends the proof of the theorem.
\qed

\subsection{Projection-type models}\label{projection-models-sec-2}

We consider the projection models discussed in Section~\ref{projection-models-sec}. The semigroup commutation properties (\ref{commutation-prop-ref}) already imply that
$$
\psi_{s,t}^{\pi}(x,\pi(Q))=\psi_{s,t}(x,\pi(Q))\quad\mbox{\rm and}\quad
\overline{\psi}_{s,t}^{\pi}(x,\pi(Q))=\overline{\psi}_{s,t}(x,\pi(Q)).
$$
Since $\pi(P)=P=P_{\pi}=\pi(P_{\pi})$ we the steady state Kalman-Bucy diffusions coincide; that is we have that
$$
\psi_{s,t}^{\pi}(x,P_{\pi})=\psi_{s,t}(x,P)\quad\mbox{\rm and}\quad
\overline{\psi}_{s,t}^{\pi}(x,P)=\overline{\psi}_{s,t}(x,P).
$$
By Theorem~\ref{theo-pi-1} we have
$$
\left\{
\begin{array}{rcl}
d\psi_{s,t}^{\pi}(x,Q)&=&\left[A-\pi\left(\phi^{\pi}_{s,t}(Q)\right)S\right]~\psi^{\pi}_{s,t}(x,Q)~dt+\pi\left(\phi_{s,t}^{\pi}(Q)\right)~C^{\prime}\Sigma^{-1}~dY_t\\
&&\\
\partial_{t}\pi(\phi_{s,t}^{\pi}(Q))&=&\ricc(\pi(\phi_{s,t}^{\pi}(Q))).
\end{array}
\right.
$$
This implies that
$$
\psi_{s,t}^{\pi}(x,Q)=\psi_{s,t}(x,\pi(Q))\quad\mbox{\rm and}\quad
\overline{\psi}_{s,t}^{\pi}(x,Q)=\overline{\psi}_{s,t}(x,\pi(Q)).
$$
Thus, we have the decompositions
$$
\psi_{s,t}^{\pi}(x,Q)-\psi_{s,t}(x,Q)=\psi_{s,t}(x,\pi(Q))-\psi_{s,t}(x,Q)$$
and
$$
\overline{\psi}_{s,t}^{\pi}(x,Q)-\overline{\psi}_{s,t}(x,Q)=\overline{\psi}_{s,t}(x,\pi(Q))-\overline{\psi}_{s,t}(x,Q).
$$
These formulae show that the convergence analysis of both $\psi_{s,t}^{\pi}(x,Q)-\psi_{s,t}(x,Q)$ and $\overline{\psi}_{s,t}^{\pi}(x,Q)-\overline{\psi}_{s,t}(x,Q)$ to $0$, as the time horizon $(t-s)\uparrow\infty$, reduces exactly to the stability properties of the Kalman-Bucy diffusion discussed in the article~\cite{ap-2016}. We point to this detailed study \cite{ap-2016} for the exact Kalman-Bucy convergence results.

\section{Some applications}\label{examples-and-applications}

\subsection{Variance inflation models}\label{var-inflation-sec}

We let $\Pi:=\{\pi_{\epsilon}\,:\,\epsilon\in [0,1]\}$ be the set of mappings
\begin{equation*}
	\pi_{\epsilon}(Q)=Q+\epsilon\,T\quad\Longrightarrow\quad\Gamma_{\pi_{\epsilon}}(Q)=\epsilon^2~TST 
\end{equation*}
indexed by $\epsilon\in [0,1]$ and a given reference matrix $T\geq 0$. In this situation, the $\delta$-balls around the identity
mapping are given for any $\delta\leq 1$ by the compact sets
$$
	B(\delta \Vert T\Vert_2)=\{\pi_{\epsilon}~:~\epsilon~\in [0,\delta]\}\subset\Pi.
$$
Conditions $\mbox{\rm (H)}_0$ and $\mbox{\rm (H)}_1$ are clearly met with
$$
\begin{array}{rclcrclccllclccrcl}
B_0&=&\epsilon^2~TST && B_1&=&0&& B_2&=&0&&
 \Ra(Q)&=&0
\\
R_{\pi}&=&R+\epsilon^2~TST&& A_{\pi}&=&A&&
S_{\pi}&=&S&\Longrightarrow&
 \Xi_{\pi}(Q)&=&\epsilon^2~TST.
\end{array}
$$
To check $\mbox{\rm (H)}_2$ we observe that
$$
\begin{array}{l}
  R^{-1/2} R_{\pi}R^{-1/2}-Id=\epsilon^2R^{-1/2} TSTR^{-1/2}\\
  \\
  \qquad\Longrightarrow ~  R^{-1/2} R_{\pi}R^{-1/2}\leq \left(1+\epsilon^2~\Vert
  R^{-1/2} TSTR^{-1/2}\Vert\right)\,Id\\
  \\
  \qquad\Longrightarrow~  R\leq R_{\pi}\leq R\,\left(1+\epsilon^2~\Vert
  R^{-1/2} TSTR^{-1/2}\Vert\right)\,Id \quad \Longrightarrow~\mbox{\rm (H)}_2.
  \end{array}
$$

In this situation Theorem~\ref{theo-continuity} yields the following corollary.
\begin{cor}
There exists some $\delta\in [0,1]$ such that for any $\epsilon\in [0,\delta]$
and for any time horizon $t\geq 0$ and any $Q\in \SS_{r}^+$ we have
$$
\Vert\phi^{\pi_{\epsilon}}_t(Q)-\phi_t(Q)\Vert_2\leq \epsilon^2~\left[ \chi_1(\delta)+ e^{-4t\nu}~\chi_2(\delta,\Vert Q\Vert_2)\right]
$$
for some finite constant $\chi_1(\delta)$, resp. $\chi_2(\delta,\Vert Q\Vert )$, whose values only depend on the parameter $\delta$, resp. on $(\delta,\Vert Q\Vert)$. In addition, for any $\epsilon\in [0,\delta]$ we have
$$
\Vert P_{\pi_{\epsilon}}-P\Vert_2\leq \epsilon^2~\chi_1(\delta).
$$
\end{cor}

\subsection{Block-diagonal localization} \label{block-diag-sec}

 Assume the covariance matrices associated with the Kalman filter in (\ref{nonlinear-KB-Riccati}) satisfy the property,
 $$
	\exists~\iota >0\quad:\quad \forall t\geq 0\quad \vert i-j\vert>\iota \Longrightarrow P_t(i,j)=0.
 $$
 In words, the coordinates of the signal have been arranged so that the $\iota$-long (or longer) range interactions between the state coordinates are null. The above condition is met if and only if the matrices $P_t$ are block-diagonal. Since the state variables are Gaussian, this property is equivalent to the fact that the state block components are block-two-by-two marginally independent. In this case, the signal-observation process $(X_t,Y_t)=(X_t[k],Y_t[k])_{1\leq k\leq n}$ defined in (\ref{lin-Gaussian-diffusion-filtering}) can be decomposed into $n$-independent $(r[k]\times r^{\prime}[k])$-dimensional filtering problems $(X_t[k],Y_t[k])$ of the form
 \begin{equation*}
 \left\{
\begin{array}{rcl}
dX_t[k]&=&A[k]~X_t[k]~dt~+~R^{1/2}[k]~dW_t[k]\\
dY_t[k]&=&C[k]~X_t[k]~dt~+~\Sigma^{1/2}[k]~dV_{t}[k]\quad\mbox{\rm with}\quad
1\leq k\leq n.
\end{array}
\right. 
 \end{equation*}
with $r=\sum_{1\leq i\leq n}r[i]$. In this elementary case, the resulting Kalman-Bucy filter and the associated Riccati equation collapse to $n$ independent evolution equations. In this case, the drift and the sensor matrices $(A,C)$, as well as the covariance matrices $(R,\Sigma)$ and $P_t$ are block-diagonal matrices of appropriate dimensions.  
 
Now observe that the sample covariance matrices $p_t(i,j)$ are generally non-null, even if $P_t(i,j)=0$. To mask these noisy entries, we use a localization mapping given in (\ref{ex-Ba}). It is readily checked that the mapping $\pi$ satisfies the orthogonality condition $\mbox{\rm (H)}_3$  discussed in (\ref{commutation}) with the cellular algebra $\Ba=\Ma_{r[1]}\oplus\ldots\oplus\Ma_{r[n]}$. With a little extra work, we can also check that
 $$
 n^{-1}J\leq L\leq  r^{\star}~Id\Rightarrow   n^{-1}~Q\leq \pi(Q)\leq  r^{\star}~\mbox{Diag}(Q(1,1),\ldots,Q(r,r))
 $$
 with $r^{\star}:=\vee_{1\leq k\leq n}{r(k)}$.
 
The central idea behind these mask-regularisations is to transform a given sample covariance matrix $p$ into some covariance
matrix with the same sparsity pattern as the limiting covariance $P$; or in practice, to mask spurious ``long-range'' correlations that are (almost) null in the true covariance. This idea is relevant in numerous applications of the {\tt EnKF} in which state-space interaction and signal observations are mostly local, and a kind-of ``decay-of-correlation'' effect is present; see \cite{Houtekamer2001,hamill01,sakov11}.

One difficulty is ensuring the mask-matrix $L$ is positive definite so that the projection $L\odot p$ is a positive map. In the block-diagonal model discussed above this property is clearly satisfied. In more general situations, several strategies can be underlined. The first one is to design mask-matrices as linear combinations $L=\sum_{i=1}^{n}l_i~z_i~z_i^{\prime}$ of unit rank vectors $z_i$, with $l_i\geq 0$.

\subsection{Bose-Mesner projections}\label{Bose-Mesner-sec}

We introduce the Bose-Mesner algebra and relevant projections and applications here. For a more thorough discussion on Bose-Mesner algebras and their application in statistical and quantum physics, combinatorics, coding, graph theory, and statistical covariance analysis (more particularly in experimental designs) we refer to the seminal article of Bose-Mesner~\cite{bose}, the ones of Nelder~\cite{nelder-1,nelder}, the more recent articles~\cite{delsarte,cameron,grossmann}, as well as the books~\cite{brouwer,bailey}.

\subsubsection{Association schemes}
We set $\Ia=\{1,\ldots,r\}$ the index set of the coordinates of the signal. Let $\Pa=\cup_{0\leq q\leq n}\Pa_q$ be an $n$-partition of the product set $\Ia^2$ such that
\begin{itemize}
\item The associated classes $\Pa_q$  are symmetric for any $0\leq q\leq n$, and  $\Pa_0:=\{(i,i)~:~i\in \Ia\}$.
\item For any $0\leq q_1,q_2\leq n$, there exists some integer $w^{q}_{q_1,q_2}$ 
(the parameters of the scheme; a.k.a. parameters of the first kind or the structural constants) such that
$$
\forall 0\leq q\leq n\quad
\forall (i,j)\in \Pa_{q}\quad w^{q}_{q_1,q_2}=
\mbox{\rm Card}\left\{k\in \Ia~:~(i,k)\in \Pa_{q_1}~\quad~(k,j)\in \Pa_{q_2}\right\}.
$$
\end{itemize}

These association schemes can be interpreted as a partition of the edges/arcs of a complete graph (with vertex set $\Ia$) into $n$ classes, often thought of as color classes. In this representation, there is a loop at each vertex and all the loops receive the same $0$-th color. The number of triangles with a fixed arc-base with color $q$ and the other two arcs with colors $q_1$ and $q_2$ is a number $w^{q}_{q_1,q_2}$ that doesn't depend on the choice of the arc-base. Each vertex $i$ is contained in exactly $v_q$ arcs with color $q$. The number $v_q$ is called the valency of the relation induced by $\Pa_q$. The parameters $w^{q}_{q_1,q_2}=w^{q}_{q_2, q_1}$
 are called the parameters of the scheme (a.k.a. parameters of the first kind or the structural constants).

For each $1\leq q\leq n$ we let $B_q$ be the adjacency matrix; that is
$$
B_q(k,l)=1_{(k,l)\in \Pa_q}=B_q(l,k)\quad\Longrightarrow\quad B_0=Id\quad\mbox{\rm and}\quad
\sum_{0\leq q\leq n}B_q=J.
$$
We also have 
$$
B_{q_1}B_{q_2}=B_{q_2}B_{q_1}=\sum_{0\leq q\leq n}w^{q}_{q_1,q_2}~B_{q}
\quad\mbox{\rm
and}\quad
B_{q}J=JB_{q}=v_q~J.
$$
This shows that $B_q$ has exactly $v_q$ non-zero entries in every row and every column.
Since   for any $q_1\not=q_2$ we have
$$
(B_{q_1}\circ B_{q_2})(k,l)=1_{(k,l)\in \Pa_{q_1}\cap \Pa_{q_2}}=0\Rightarrow  B_{q_1}\circ B_{q_2}=1_{q_1=q_2}~B_{q_1}
$$
the set $\Ba$ is also closed w.r.t. the Hadamard product and contains $I,J$. Thus, the set $$
\Ba:=\left\{\sum_{0\leq q\leq n}b_q~B_q~:~b=(b_q)_{0\leq q\leq n}\subset\RR^{n+1}\right\}
$$ 
is an associative commutative algebra called the Bose-Mesner algebra of the association scheme. Notice that $\Ba$ is also a matrix $\star$-algebra (i.e. closed by matrix multiplication, the transposition, addition and the scalar multiplication). These special cases of finite dimensional $\CC^{\star}$-algebra are unitarily equivalent to block-diagonal matrices. By a theorem of Von Neumann we also mention that the orthogonal projection on any matrix $\star$-algebra is a positive map.
 
An illustration when $n=2$ and $r=6$ is provided by 
 $$
B_0=Id, \qquad B_1=\left[
 \begin{array}{cccccc}
 0&1&1&0&0&0\\
 1&0&1&0&0&0\\
 1&1&0&0&0&0\\
  0&0&0&0&1&1\\
   0&0&0&1&0&1\\
    0&0&0&1&1&0\\
 \end{array}
 \right]\quad \mbox{\rm and}\quad B_2=J-\left[B_0+B_1\right].
  $$
 In this case we have $B_1^2=2Id+B_1=2B_0+B_1$ and $B_1B_2=0=B_2B_1$.

\subsubsection{Minimal orthogonal projections}

The commuting matrices $B_q$  are simultaneously diagonalizable, $\Ba$ has  a basis of minimal orthogonal idempotents $D_i$; that is, we have that
$$
D_{q_1}D_{q_2}=1_{q_1=q_2}D_{q_1}\quad\mbox{\rm and}\quad \sum_{0\leq q\leq r} D_q=Id.
$$
Without any loss of generality we can choose $D_0=r^{-1}J$. The matrices $D_q$ are called the minimal idempotents of the algebra $\Ba$.
In addition, the column vectors  $D_{i,1}\ldots,D_{i,r}$  of $D_i$ are the eigenvectors of any matrix in $\Ba$. 
The eigenvector spaces $\Da_{i}=\mbox{\rm Span}(D_{i,1}\ldots,D_{i,r})$  are mutually orthogonal and every vector $u\in \RR^r$
can be expressed uniquely as $u=\sum_{0\leq q\leq n}u_i$ with $u_i\in \Da_i$ (notice that $\Da_0$ is the $1$-dimensional space of 
constant vectors). Also notice that the dimension of  $\Da_{i}$ equals to the rank of $D_i$, which is equal to the trace of $D_i$ (since all non-zero eigenvalues of $D_i$ are equal to $1$).

In particular, we have
$$
B_qD_k=\frac{\langle B_q,D_k\rangle_F}{\langle D_k,D_k\rangle_F}~D_k ~~\Longrightarrow~~ \lambda_k(B_q)=
\frac{\langle B_q,D_k\rangle_F}{\langle D_k,D_k\rangle_F}
$$
where $ \lambda_k(B_q)$ stands for the $k$-th eigenvalue of $B_q$. Further details on these simultaneous diagonalization can be found in~\cite{barker}.

The orthogonal projection of a matrix $Q$ on $\Ba$ is given by the formulae
$$
\pi(Q)=\proj_{\Ba}(Q):=\sum_{0\leq q\leq n}~\frac{\langle Q,B_q\rangle_F}{\langle B_q,B_q\rangle_F}~B_q=
\sum_{0\leq q\leq n}~\frac{\langle Q,D_q\rangle_F}{\langle D_q,D_q\rangle_F}~D_q.
$$

To check condition $\mbox{\rm (H)}_3$ we observe that
$$
D_{q_1}D_{q_2}=1_{q_1=q_2}D_{q_1}~~\Longrightarrow~~
\proj_{\Ba}(D_q~Q)=\frac{\langle Q,D_q\rangle_F}{\langle D_q,D_q\rangle_F}~D_q=D_q~\proj_{\Ba}(Q).
$$
This yields
$$
 \forall B\in\Ba\qquad\proj_{\Ba}(B \left[Q-\proj_{\Ba}(Q)\right])=0.
$$
For any matrix $M$ we have
$$
\langle MM^{\prime},D_q\rangle_F=\tr(D_qMM^{\prime})=\tr(M^{\prime}D^2_qM)=\tr((D_qM)^{\prime}(D_qM))\geq 0.
$$
This implies that
$$
\forall Q\in\SS_{r}^+\qquad \langle Q,D_q\rangle_F\geq 0\quad\mbox{\rm and}\quad
\proj_{\Ba}(Q)=
\sum_{0\leq q\leq n}~\frac{\langle Q,D_q\rangle_F}{\langle D_q,D_q\rangle_F}~D_q~\geq 0.
$$
This shows that the orthogonal projection is a positive map from the algebra of square
matrices into itself. In addition, it is  trace-preserving and unital in the sense that
$$
\tr(\proj_{\Ba}(Q))=\tr(Q)\quad\quad\mbox{\rm and}\quad \proj_{\Ba}(Id)=Id.
$$

Last, but not least, using the decomposition
\begin{equation}\label{trace-ref-BM}
	Q=\proj_{\Ba}(Q)\stackrel{\perp}{+}\left[Q-\proj_{\Ba}(Q)\right]~~\Longrightarrow~~\Vert Q-\proj_{\Ba}(Q)\Vert_F\leq \Vert Q\Vert_F\leq \tr(Q)
\end{equation} 
as soon as $Q\in \SS_r^+$. Working a little harder, we check that
$$
	\Vert Q-\proj_{\Ba}(Q)\Vert_F ~\leq~  \Vert Q\Vert_F~\left[1-\frac{1}{n+1}\frac{1}{\wedge_{0\leq q\leq n}\tr(D_q)}~\frac{\tr(Q)^2}{\tr(Q^2)}\right]^{1/2}.
$$

\subsubsection{Distance regular graphs}
Another prototype of Bose-Mesner algebra are distance regular graphs. Given a connected graph $\Ga=(\Va,\Ea)$ with vertex set $\Va$ and arc/edges set $\Ea$, we let $\rho(i,j)$ be the path-length distance between two vertices $i,j\in \Va$. Let $$\Sa(i,q)=\{j\in \Va~:~\rho(i,j)=q\}$$ be the sphere of radius $q$. The graph $\Ga$ is distance regular if and only if we have
$$
\mbox{\rm Card}\left(\Sa(i,q_1)\cap \Sa(j,q_2)\right)=w^{\rho(i,j)}_{q_1,q_2}
$$
for some parameters $w^{q}_{q_1,q_2}$. In other words, for every two
vertices $(i,j)$ at distance $q$ there are precisely $w^{q}_{q_1,q_2}$ vertices in the graph at 
distance $q_1$ from $i$ and $q_2$ from $j$.

In these settings, the matrices $$
(B_q)(i,j)=1_{\rho(i,j)=q}\qquad\mbox{\rm with}\quad
0\leq q\leq \mbox{\rm diam}(\Ga):=\sup_{(i,j)\in \Va^2}\rho(i,j)$$
are called the distance matrices ($B_0=Id$, $B_1$ the adjacency matrix, and so on).  In this situation, the association scheme is given by
 the partition
 $$
\forall 0\leq k\leq d:=\mbox{\rm diam}(\Ga)\qquad
 \Pa_k=\{(i,j)\in \Ia^2~:~\rho(i,j)=k\}.
$$
In addition we have $w^{q_1}_{q_2,1}=0$ for any $q_1\not=0$, and
$q_2\not=\{q_1-1,q_1,q_1+1\}$. If we set
$$
a_q:=w^{0}_{q,q}\qquad b_q:=w^{q}_{q-1,1}\quad\mbox{\rm and}\quad
c_q:=w^{q}_{q+1,1}
$$
then we have
$$
B_1B_q=c_{q-1}~B_{q-1}+(a_1-b_q-c_q)~B_q+b_{q+1}B_{q+1}
$$
and
$$
B_1B_{d}=c_{d-1}B_{d-1}+(a_1-b_d)
~B_d.$$
This shows that the adjacency matrix $B_1$ generates $\Ba$ (i.e. the matrices $B_q$ 
can be written as polynomials of degree $q$ in $B_1$), so that the eigenvalues
$\left(\lambda_k(B_1)\right)_{1\leq k\leq d}$ of $B_1$ are mutually distinct.

\subsubsection{Riccati solvers}\label{riccati-solver}
A given matrix $Q$ belongs to $\Ba$ if and only if it is constant within each block. To check this claim, we observe that
\begin{eqnarray*}
Q=\proj_{\Ba}(Q)
&\displaystyle\Longleftrightarrow& \forall 0\leq q\leq n\quad Q\odot B_q=\frac{\langle Q,B_q\rangle_F}{\langle B_q,B_q\rangle_F}~B_q\\
&\displaystyle\Longleftrightarrow& \forall 0\leq q\leq n\quad \forall (i,j)\in \Pa_q\quad 
\sum_{(k,l)\in \Pa_q}Q_{k,l}=w^0_{q,q}~Q_{i,j}\\
&\displaystyle\Longleftrightarrow&\forall 0\leq q\leq n\quad\forall (i,j),(i^{\prime},j^{\prime})\in \Pa_q\qquad ~Q_{i,j}=~Q_{i^{\prime},j^{\prime}}.
\end{eqnarray*}
In other words, the matrix is constant within each block. When $r=r^{\prime}$, then $(\pi(A),\pi(S),\pi(R))=(A,R,S)$ is satisfied as soon as $(A,R,C,\Sigma^{-1})\in\Ba$.

We further assume that $(A,R,S)\in\Ba$ and we set
$$
A:=\sum_{0\leq q\leq n}~a_q~D_q\qquad R:=\sum_{0\leq q\leq n}~r_q~D_q
\quad\mbox{\rm and}\quad S:=\sum_{0\leq q\leq n}~s_q~D_q.
$$
Let $P_0=\pi(P_0)=\sum_{0\leq q\leq n}~\alpha_q(0)~D_q$ be some covariance matrix in $\Ba$. By Theorem~\ref{theo-pi-1}
we have
$$
P_t=\pi\left(P_t\right)=\sum_{0\leq q\leq n}~\alpha_q(t)~D_q.
$$
In addition, we have
\begin{eqnarray*}
\partial_tP_t&=&\sum_{0\leq q\leq n}~\partial_t\alpha_q(t)~D_q~=~\ricc(P_t)~=~
A\pi(P_t)-\pi(P_t)A^{\prime}+R-\pi(P_t)S\pi(P_t)\\
&=&
\sum_{0\leq q\leq n}~\left[2a_{q}~\alpha_{q}(t)+r_q-~\alpha_q(t)^2~s_q\right]~D_q.
\end{eqnarray*}
This implies that
$$
\left\{
\begin{array}{rcl}
\partial_t\alpha_q(t)&=&2a_{q}~\alpha_{q}(t)+r_q-\alpha_q(t)^2~s_q\\
q&=&0,\ldots,n.
\end{array}
\right.
$$
When $s_q\not=0\not=r_q$ this collection of Riccati equations
take the form
$$
\partial_t\alpha_q(t)=-s_q~\left(\alpha_q(t)-z_1(q)\right)~\left(\alpha_q(t)-z_2(q)\right)
$$
with the couple of roots
$$
z_1(q)=\frac{a_q-\sqrt{a_q^2+s_qr_q}}{s_q}<0<z_2(q)=\frac{a_q+\sqrt{a_q^2+s_q r_q}}{s_q}.
$$
The solutions of the above equations
 are given by the formulae:
\begin{equation*}
\alpha_q(t)-z_2~=~(\alpha_q(0)-z_2(q))~\frac{(z_2(q)-z_1(q))~e^{-2t\sqrt{a_q^2+s_q r_q}}}{(z_2(q)-\alpha_q(0))~e^{-2t\sqrt{a_q^2+s_q r_q}}+(\alpha_q(0)-z_1(q))}~\longrightarrow_{t\rightarrow\infty}~0.
\end{equation*}

\subsection{Stein-Shrinkage models} \label{stein-sec}

Stein-Shrinkage models are an extension of the variation inflation model to parameters
$\epsilon=\epsilon(Q)$ and target-type matrices $T=T(Q)$ that both may depend on the matrix $Q$. These models are defined by the formula
\begin{equation*}
	\pi(Q)=\epsilon(Q)~T(Q)+(1-\epsilon(Q))~Q
\end{equation*}
for some function $Q\mapsto\epsilon(Q)\in [0,1]$ and some mapping $T$ from $\SS^+_r$ into itself. It is not within scope of this article to review all the relevant covariance matrix estimators encountered in the statistics literature fitting this general model. We will just illustrate this model with three important and currently used approximations:
\begin{itemize}
\item {\em Mask matrix estimates} are associated with mappings $T$ defined by $T(Q):=L\odot Q$
with a matrix
$L$ of the form
\begin{equation}\label{mask-matrix}
L_{i,j}=1_{\vert i-j\vert< \iota}~\Longrightarrow~Q-L\odot Q=1_{\vert i-j\vert\geq  \iota}~Q_{i,j}.
\end{equation}
\item {\em Maximum likelihood} type estimates are associated with mappings $T$ defined by
$$
T(Q):=\argmax_{q\in\SS^+_r}{\left(\log{\mbox{\rm det}(q)}+\tr(q^{-1}Q)+\alpha~\Vert L\odot q\Vert\right)}
$$ 
for some $\alpha>0$, some mask matrix $L$~\cite{bien,chaudhuri,khare,lam-fan}, and some matrix norm $\Vert \point\Vert$ on $\SS^+_r$. 
\item {\em Nystr\"om estimates} are associated with mappings $T$ defined by
\begin{equation}\label{nylstrom-1}
 T(Q)=\left(J-L_{\Pa^c}\right)\odot Q+L_{\Pa^c}\odot \left[
Q_{\Pa^c,\Pa}~Q_{\Pa}^{-1}~Q_{\Pa,\Pa^c}\right]
\end{equation}
where $\{1,\ldots,r\}=\Pa\cup \Pa^c$ stands for a partition of the index coordinate set and 
$L_{\Pa^c}$ stands for the mask matrix defined by
$$
L_{\Pa^c}(i,j)=1_{\Pa^c\times\Pa^c}(i,j).
$$
At the level of the sample covariance matrices $p_0$, the matrix $T(p_0)$ is obtained by taking 
the sample covariance matrix associated with projection $\Ta_{\Va}(\zeta_l)$ of the state particle vectors 
$$
\zeta^{\prime}:=
\left[\begin{array}{c}
\zeta_1^{\prime}\\
\vdots\\
\zeta_r^{\prime}\end{array}
\right]:=\left[\xi^1_0-m_0,\ldots,\xi^N_0-m_0\right]=\left[\begin{array}{ccc}
\xi^1_0(1)-m_0(1)&\ldots&\xi^N_0(1)-m_0(1)\\
\vdots&\vdots&\vdots\\
\xi^1_0(r)-m_0(r)&\ldots&\xi^N_0(r)-m_0(r)\end{array}
\right]~
$$
 onto the vector space $\Va_{\Pa}$ of $\RR^N$ spanned by the random vectors
 $$
 V_i=\zeta_{k_i}:=\left[
 \begin{array}{c}\xi^1_0(k_i)-m_0(k_i)\\
 \vdots\\
 \xi^N_0(k_i)-m_0(k_i)
 \end{array}\right]
 \in \RR^N\quad\mbox{\rm
 with $\Pa=\{k_1,\ldots,k_s\}$ and $s:=\mbox{\rm Card}(\Pa)\leq r$}.
 $$
More precisely, if we set
\begin{eqnarray*}
N~T(p_0)&=&\left[\begin{array}{c}
\left(\Ta_{\Va}\zeta_l\right)^{\prime} \\
\vdots\\
\left(\Ta_{\Va}\zeta_r\right)^{\prime}\end{array}
\right]\left[
\Ta_{\Va}\zeta_l,
\ldots,
\Ta_{\Va}\zeta_l
\right]\\
&=&(\Ta_{\Va}\zeta)^{\prime}\Ta_{\Va}\zeta=\zeta^{\prime}\Ta_{\Va}\zeta
=\left[\begin{array}{ccc}
\langle \Ta\zeta_1,\Ta_{\Va}\zeta_1\rangle&\ldots&\langle\Ta\zeta_1,\Ta_{\Va}\zeta_r\rangle\\
\vdots\\
\langle\Ta\zeta_r,\Ta_{\Va}\zeta_1\rangle&\ldots&\langle \Ta_{\Va}\zeta_r,\Ta_{\Va}\zeta_r\rangle\end{array}
\right]
\end{eqnarray*}
then we have that
\begin{equation}\label{last-nylstrom}
\EE\left(T(p_0)\right)=\displaystyle T(P_0)+\frac{s}{N}~L_{\Pa^c}\odot \left[Q_{\Pa^c}-Q_{\Pa^c,\Pa}~Q_{\Pa}^{-1}~Q_{\Pa,\Pa^c}\right].
\end{equation} 
The proof of this bias property and related variance estimates can be found in~\cite{arcolano}. For the convenience of the reader a proof of the last assertion is provided in the Appendix.

\end{itemize}

For mask type mappings of the form (\ref{mask-matrix}), condition $\mbox{\rm (H)}_0$ is satisfied by first letting
$$
	(B_0,B_1,B_2)=(0,0,0)~\Longrightarrow~ \Gamma_{\pi}(Q)=\Ra(Q):=\epsilon(Q)^2~(L\odot Q-Q)S(L\odot Q-Q).
$$
To ensure the uniform estimate $\sup_{Q\in \SS^+_{r}}\Vert\Ra(Q)\Vert_2<\infty$ holds we use Gershgorin's theorem to show that
$$
\Vert Q-L\odot Q\Vert_2~\leq ~l_{\iota}(Q)~:=~\sup_{1\leq i\leq r}~\sum_{\vert i-j\vert\geq  \iota}\vert Q_{i,j}\vert.
$$
This yields
$$
\Vert \Ra(Q)\Vert_2\leq ~\epsilon^2(Q)~\Vert S\Vert_2~l^2_{\iota}(Q)\Rightarrow 
\Ra(Q)\leq  \overline{\Ra}(Q)~Id\quad\mbox{\rm with}\quad
 \overline{\Ra}(Q)=\epsilon^2(Q)~\Vert S\Vert_2~l^2_{\iota}(Q).
$$
When $ l_{\iota}(Q)$ is too large, the quadratic perturbation may have some destabilizing effects. To avoid these issues
we assume that $\epsilon(Q)$ is chosen so that 
$$
\epsilon(Q)=\epsilon_1~1_{l_{\iota}(Q)\leq \epsilon_2^{-1}}~\Longrightarrow~\Ra(Q)\leq  \varpi~Id\quad\mbox{\rm with}\quad
\varpi=~\Vert S\Vert_2~(\epsilon_1/\epsilon_2)^2
$$
for some  $\epsilon_1\in[0,1]$, and some threshold $\epsilon_2>0$. In this case, condition $\mbox{\rm (H)}_1$ is also met with
$$
R_{\pi}=R+(\epsilon_1/\epsilon_2)^2~\Vert S\Vert_2~Id\qquad A_{\pi}=A\quad\mbox{\rm and}\quad
S_{\pi}=S
\Longrightarrow
 \Xi_{\pi}(Q)=(\epsilon_1/\epsilon_2)^2~\Vert S\Vert_2~Id.
$$
Arguing as in the end of Section~\ref{var-inflation-sec} we have
$$
\begin{array}{l}
  R^{-1/2} R_{\pi}R^{-1/2}-Id=(\epsilon_1/\epsilon_2)^2\Vert S\Vert_2~R^{-1}\\
  \\
  \qquad \Longrightarrow  ~ R^{-1/2} R_{\pi}R^{-1/2}\leq \left(1+(\epsilon_1/\epsilon_2)^2~\Vert S\Vert_2~~\Vert
  R^{-1}\Vert\right)~Id\\
  \\
  \qquad \Longrightarrow ~ R\leq R_{\pi}\leq R~\left(1+(\epsilon_1/\epsilon_2)^2~\Vert S\Vert_2~\Vert
  R^{-1}\Vert\right)~Id \quad \Longrightarrow~\mbox{\rm (H)}_2.
  \end{array}
$$
Now we can consider the set 
$$
\Pi=\{\pi_{\epsilon_1,\epsilon_2}~:~(\epsilon_1,\epsilon_2)\in ([0,1]\times [\delta,\delta^{-1}] )\}
$$
for some given parameter $\delta$ and the just described mappings $\pi_{\epsilon_1,\epsilon_2}$ given by
$$
\pi_{\epsilon_1,\epsilon_2}(Q)=Q+\epsilon_1~1_{l_{\iota}(Q)\leq \epsilon_2^{-1}}~\left[L\odot Q-Q\right]~\Longrightarrow~
\Vert \pi_{\epsilon_1,\epsilon_2}-id\Vert_2\leq \epsilon_1/\epsilon_2.
$$ 
The associated $\delta$-balls around the identity mapping are given in this case by the compact sets
$$
B(\delta)=\{\pi_{\epsilon_1,\epsilon_2}~:~\epsilon_1/\epsilon_2\leq \delta\}
$$
for any $\delta\leq 1$.

More generally, the Stein-Shrinkage models discussed above can be extended without further work to general mappings of the following form
$$
\pi_{\epsilon_1,\epsilon_2}(Q)=Q+\epsilon_1~1_{l_T(Q)\leq \epsilon_2^{-1}}~\left[T(Q)-Q\right]~
$$
where $T$ stands for some mapping from $\SS_{r}^+$ into itself such that
$$
\Vert T(Q)-Q\Vert_2\leq l_T(Q)~
$$
for some mapping $Q\in \SS_{r}^+\mapsto l_T(Q)\in [0,\infty[$. Further examples of such mappings include the Bose-Mesner projections $T(Q)=\proj_{\Ba}(Q)$ discussed in Section~\ref{Bose-Mesner-sec} and which can be seen to fit this model via the trace operator in (\ref{trace-ref-BM}). 

In this general setting, Theorem~\ref{theo-continuity} yields the following corollary.
\begin{cor}
There exists some $\rho\in [0,1]$ such that for any  $(\epsilon_1,\epsilon_2)\in ([0,1]\times [\delta,\delta^{-1}] )$
with $\epsilon_1\leq \rho~\epsilon_2$, for any $Q\in \SS_{r}^+$ 
and any time horizon $t\geq 0$ we have
$$
\Vert\phi^{\pi_{\epsilon_1,\epsilon_2}}_t(Q)-\phi_t(Q)\Vert_2\leq (\epsilon_1/\epsilon_2)^2~\left[ \chi_1(\rho)+ e^{-4t\nu}~\chi_2(\rho,\Vert Q\Vert_2)\right]
$$
for some finite constant $\chi_1(\rho)$, resp. $\chi_2(\rho,\Vert Q\Vert )$, whose values only depend
on the parameter $\delta$, resp. on $(\rho,\Vert Q\Vert)$. In addition, for any  $\epsilon_1\leq \rho~\epsilon_2$ we have
$$
\Vert P_{\pi_{\epsilon_1,\epsilon_2}}-P\Vert_2\leq (\epsilon_1/\epsilon_2)^2~\chi_1(\rho).
$$
\end{cor}

This section illustrates how our first class of perturbation-type model captures most projection-type mappings; and consequently those results relevant to perturbation-type mappings are applicable to projection-type models (but not vice-versa).

\subsection{Mean repulsion models} \label{mean-repulsion-sec}
 
 The preceding subsections were concerned with perturbation and projecting mappings $\pi$ that directly fell within the class of models defined by (\ref{two-classes}). We also illustrated how the first class of perturbation-type model captures most projection-type mappings considered in (\ref{two-classes}). 
 
 In this subsection we illustrate that our main result (viz. Theorem~\ref{theo-continuity} and $\mbox{\rm (H)}_{2}$) on the robustness and boundedness of perturbed Riccati semigroups, captures a larger class of perturbation-type models than those simply defined by the condition $\mbox{\rm (H)}_{0}$ and (\ref{pi-Riccati-def-gamma-term}). Of course, Theorem~\ref{theo-continuity} also applies under the more constrained condition $\mbox{\rm (H)}_{0}$ as a special case, and $\mbox{\rm (H)}_{0}$ is still of specific interest in, e.g., the variance inflation and Stein-Shrinkage-type models discussed in the preceding subsections. However, $\mbox{\rm (H)}_{0}$ is not satisfied by the perturbation scheme considered in this subsection. Nevertheless, $\mbox{\rm (H)}_{2}$ is satisfied, and thus Theorem~\ref{theo-continuity} still applies.

 As their name indicates, mean repulsion models are defined by adding an extra repulsion
term around the sample averages in the nonlinear diffusion (\ref{Kalman-Bucy-filter-nonlinear-ref}).
Consider the nonlinear diffusion
\begin{eqnarray*}
d\overline{X}_t&=&\left[A~\overline{X}_t~dt-T_1(P_t)(\overline{X}_t-\widehat{X}_t)~\right]dt+~R^{1/2}~d\overline{W}_t\\
\\
&&\hskip4cm+P_tC^{\prime}\Sigma^{-1}~\left[dY_t-
\left(C\left(
\overline{X}_t+T_2(\overline{X}_t-\widehat{X}_t)\right)dt+\Sigma^{1/2}~d\overline{V}_{t}\right)\right]\\
\\
&=&[A-P_tS]\overline{X}_t-[
T_1(P_t)+P_tST_2]~(\overline{X}_t-\widehat{X}_t)~dt\\
&&\hskip6cm+~R^{1/2}~d\overline{W}_t+
P_tC^{\prime}\Sigma^{-1}~\left[dY_t-\Sigma^{1/2}~d\overline{V}_{t}\right]
\end{eqnarray*}
where $T_1\,:\,\SS_{r}^+\mapsto \MM_{r} $ stands for some mapping and $T_2$ some given matrix. 

A key feature of this class of mean repulsion models is that their $\Fa_t$-conditional projections 
coincide with the Kalman-Bucy filter, only their conditional covariance matrices are altered.

To describe the Riccati equation associated with this class of nonlinear diffusions we observe that
\begin{eqnarray*}
d(\overline{X}_t-\widehat{X}_t)
&=&\left(A-P_tS-
[T_1(P_t)+P_tST_2]~\right)~(\overline{X}_t-\widehat{X}_t)~dt+~R^{1/2}~d\overline{W}_t
-
P_tC^{\prime}R^{-1/2}_{2}d\overline{V}_{t}.
\end{eqnarray*}
Thus, the covariance evolution equation is given by the Riccati equation
\begin{eqnarray*}
\partial_tP_t&
=&\left[A-P_tS\left(Id+T_2\right)-
T_1(P_t)\right]~P_t+P_t\left[A-P_tS\left(Id+T_2\right)-
T_1(P_t)\right]^{\prime}
+R+P_tSP_t\\
\\
&=&AP_t
+P_tA^{\prime}
+R-P_tSP_t
-P_tST_2P_t-(T_1(P_t)P_t+P_tT_1(P_t)^{\prime})-P_tT_2SP_t.
\end{eqnarray*}

For instance, choosing
$$
T_1(Q)=\epsilon_1 QS\quad\mbox{\rm and}\quad T_2=\epsilon_2 Id
$$
for some $(\epsilon_1,\epsilon_2)$ such that $(\epsilon_1+\epsilon_2)>-1/2$ we find that
 \begin{eqnarray*}
\partial_tP_t&=&AP_t
+P_tA^{\prime}+R-P_tS_{\epsilon}P_t\quad\mbox{\rm with}\quad S_{\epsilon_1,\epsilon_2}:=(1+2(\epsilon_1+\epsilon_2))S.
\end{eqnarray*}

We let $\phi_{\epsilon,t}$ be the Riccati semigroup associated with the above equation, with $\epsilon=(\epsilon_1,\epsilon_2)\in \Pi=[0,1]^2$. Theorem~\ref{theo-continuity} yields the following corollary.
\begin{cor}
There exists some $\delta\in [0,1]$ such that for any $\epsilon=(\epsilon_1,\epsilon_2)\in [0,\delta]^2$
and for any time horizon $t\geq 0$ and any $Q\in \SS_{r}^+$ we have
$$
\Vert\phi_{\epsilon,t}(Q)-\phi_t(Q)\Vert_2\leq 2(\epsilon_1+\epsilon_2)~\left[ \chi_1(\delta)+ e^{-4t\nu}~\chi_2(\delta,\Vert Q\Vert_2)\right]
$$
for some finite constant $\chi_1(\delta)$, resp. $\chi_2(\delta,\Vert Q\Vert )$, whose values only depend
on the parameter $\delta$, resp. on $(\delta,\Vert Q\Vert)$. In addition, if $P_{\epsilon}=\phi_{\epsilon,t}(P_{\epsilon})$ is the fixed point of $\phi_{\epsilon,t}$, then for any $\epsilon\in [0,\delta]$ we have
$$
\Vert P_{\epsilon}-P\Vert_2\leq 2(\epsilon_1+\epsilon_2)~~\chi_1(\delta).
$$
\end{cor}

\section*{Appendix} 

\subsection*{Proof of formula (\ref{pi-Riccati-def})}

Let $\phi^{\pi}_{s,t}$ be the semigroup of equation (\ref{pi-Riccati-def}). Also let $\overline{X}_t^{\pi}$ be the time non-homogeneous diffusion given by the equation
\begin{eqnarray*}
d\overline{X}_t^{\pi}&=&A~\overline{X}_t^{\pi}~dt~+~R^{1/2}~d\overline{W}_t+\pi\left(\phi^{\pi}_{0,t}(Q)\right)~C^{\prime}\Sigma^{-1}~\left[dY_t-\left(C\overline{X}^{\pi}_tdt+\Sigma^{1/2}~d\overline{V}_{t}\right)\right]\\
&=&\left[A-\pi\left(\phi^{\pi}_{0,t}(Q)\right)S\right]~\overline{X}_t^{\pi}~dt~+\pi\left(\phi^{\pi}_{0,t}(Q)\right)~C^{\prime}\Sigma^{-1}~dY_t+dM_t^{\pi}
\end{eqnarray*}
with the $r$-valued martingale
$$
dM_t^{\pi}:=R^{1/2}~d\overline{W}_t
-\pi\left(\phi^{\pi}_{0,t}(Q)\right)C^{\prime}R^{-1/2}_{2}~d\overline{V}_{t}
$$
with covariation matrix 
$$
\partial_t\langle M^{\pi}(k),M^{\pi}(l)\rangle_t=\left[
R+\pi\left(\phi^{\pi}_{0,t}(Q)\right)S\pi\left(\phi^{\pi}_{0,t}(Q)\right)\right](k,l).
$$
We have
$$
\begin{array}{l}
\displaystyle\overline{X}_t^{\pi}
\displaystyle=\exp{\left(\oint_0^t\left[A-\pi\left(\phi^{\pi}_{0,s}(Q)\right)S\right]~ds\right)}~\overline{X}_0^{\pi}
\\
\\
\hskip3cm\displaystyle+\int_0^t
\exp{\left(\oint_s^t\left[A-\pi\left(\phi^{\pi}_{0,u}(Q)\right)S\right]~du\right)}~\pi\left(\phi^{\pi}_{0,s}(Q)\right)~C^{\prime}\Sigma^{-1}dY_s\\
\\
\hskip5cm\displaystyle+\int_0^t
\exp{\left(\oint_s^t\left[A-\pi\left(\phi^{\pi}_{0,u}(Q)\right)S\right]~du\right)}~dM_s^{\pi}.
\end{array}
$$
This implies that the conditional expectations $\widehat{X}_t^{\pi}=\EE(\overline{X}_t^{\pi}~\vert~\Fa_t)$ are given by the formula
$$
\begin{array}{l}
\widehat{X}_t^{\pi}
\displaystyle=\exp{\left(\oint_0^t\left[A-\pi\left(\phi^{\pi}_{0,s}(Q)\right)S\right]~ds\right)}~\widehat{X}_0^{\pi}
\\
\\
\hskip3cm\displaystyle+\int_0^t
\exp{\left(\oint_s^t\left[A-\pi\left(\phi^{\pi}_{0,u}(Q)\right)S\right]~du\right)}~\pi\left(\phi^{\pi}_{0,s}(Q)\right)~C^{\prime}\Sigma^{-1}dY_s.
\end{array}
$$
Equivalently, we have
$$
d\widehat{X}_t^{\pi}=A~\widehat{X}_t^{\pi}~dt~+\pi\left(\phi^{\pi}_{0,t}(Q)\right)~C^{\prime}\Sigma^{-1}~\left[dY_t-C\widehat{X}_t^{\pi}dt\right]
$$
from which we prove that
$$
d\left[\overline{X}_t^{\pi}-\widehat{X}_t^{\pi}\right]
=
\left[A-\pi\left(\phi^{\pi}_{0,t}(Q)\right)S\right]~\left[\overline{X}_t^{\pi}-\widehat{X}_t^{\pi}\right]~dt~+~dM_t^{\pi}.
$$
This implies that the covariation matrices 
\begin{eqnarray*}
Q_t^{\pi}&:=&\EE\left(\left[\overline{X}_t^{\pi}-\widehat{X}_t^{\pi}\right]\left[\overline{X}_t^{\pi}-\widehat{X}_t^{\pi}\right]^{\prime}~\vert~\Fa_t\right)=\EE\left(\left[\overline{X}_t^{\pi}-\widehat{X}_t^{\pi}\right]\left[\overline{X}_t^{\pi}-\widehat{X}_t^{\pi}\right]^{\prime}\right)
\end{eqnarray*}
don't depend on the observation process, and they satisfy the equation
$$
\partial_tQ^{\pi}_t=\left[A-\pi\left(\phi^{\pi}_{0,t}(Q)\right)S\right]Q^{\pi}_t+Q^{\pi}_t\left[A-\pi\left(\phi^{\pi}_{0,t}(Q)\right)S\right]^{\prime}+R+\pi\left(\phi^{\pi}_{0,t}(Q)\right)S\pi\left(\phi^{\pi}_{0,t}(Q)\right).
$$
Recalling that $\phi^{\pi}_{0,t}(Q)$ is the Riccati semigroup of the equation (\ref{pi-Riccati-def}) we have
$$
\partial_t\left(Q^{\pi}_t-\phi^{\pi}_{0,t}(Q)\right)=\left[A-\pi\left(\phi^{\pi}_{0,t}(Q)\right)S\right]\left(Q^{\pi}_t-\phi^{\pi}_{0,t}(Q)\right)+\left(Q^{\pi}_t-\phi^{\pi}_{0,t}(Q)\right)\left[A-\pi\left(\phi^{\pi}_{0,t}(Q)\right)S\right]^{\prime}$$
We conclude that
$$
Q^{\pi}_0=Q~~\Longrightarrow~~
Q^{\pi}_t=\phi^{\pi}_{0,t}(Q)=\Pa_{\eta^{\pi}_t}~~\Longrightarrow~~ \pi\left(Q^{\pi}_t\right)=\pi\left(\phi^{\pi}_{0,t}(Q)\right)
$$
where $\eta^{\pi}_t=\mbox{\rm Law}(\overline{X}_t^{\pi}~\vert~\Fa_t)$. This ends the proof of (\ref{pi-Riccati-def}). See also \cite[page 242]{gelb74} (among numerous other sources) for the related covariance flow of a Kalman filter with an arbitrary gain matrix.

\subsection*{Proof of Lemma~\ref{lemm-comp-flows}}\label{ref-lemm-comp-flows}
Condition (\ref{comparison-Wa-condition}) implies that
\begin{eqnarray*}
\lambda_{min}\left(\Wa_t(\pi_2)\right)\geq \varpi_{-,t}(\pi)&\Longrightarrow&
\lambda_{min}\left(\Wa_t(\pi_2)^{1/2}\right)\geq \sqrt{\varpi_{-,t}(\pi_2)}\\
&\Longrightarrow&\lambda_{max}\left(\Wa_t(\pi_2)^{-1/2}\right)\leq 1/\sqrt{\varpi_{-,t}(\pi_2)}
  \end{eqnarray*}
from which we conclude that
\begin{equation}\label{estimate-Wa-1}
\Vert \Wa_t(\pi_2)^{-1/2}\Vert_2^2\leq\varpi_{-,t}(\pi_2)^{-1}.
\end{equation}
We also have
\begin{equation}\label{comparison-Wa}
(\ref{comparison-Wa-condition})~
\Longrightarrow~
\Wa_t(\pi_1)~\geq~ \Wa_t(\pi_1,\pi_2) ~:=~\int_0^t    \Ua_s(\pi_1)  \Va_s(\pi_2)~\Ua^{\prime}_s(\pi_1)~ds.
\end{equation}
Observe that
$$
\partial_s\Wa_s(\pi_1,\pi_2)= \Ua_s(\pi_1,\pi_2)
\left[\partial_s\Wa_s(\pi_2)\right]  \Ua_s(\pi_1,\pi_2)^{\prime}
$$
with the flow of matrices
$$
\Ua_s(\pi_1,\pi_2)= \Ua_s(\pi_1) ~\Ua_s(\pi_2)^{-1}~ \Longrightarrow~\Ua_s(\pi,\pi)=Id.
$$
We set
$$
\Vert\Ua\Vert_2:=\sup_{(s,\pi)\in([0,t]\times\Pi)}{\Vert\Ua_s(\pi)\Vert_2}<\infty
\quad\mbox{\rm and}\quad
\Vert\Va\Vert_2:=\sup_{(s,\pi)\in([0,t]\times\Pi)}{\Vert\Va_s(\pi)\Vert_2}<\infty.
$$
In this notation, using the fact that
$$
 \sup_{s\in [0,t]}{\Vert\partial_s\Wa_s(\pi_2)\Vert_2}\leq t~\Vert\Ua\Vert_2^2~\Vert\Va\Vert_2
$$
we find that
$$
\sup_{s\in [0,t]}{\Vert
\partial_s\Wa_s(\pi_1,\pi_2)-\partial_s\Wa_s(\pi_2)\Vert_2}\leq 
c_{\Ua}~\Vert \pi_1-\pi_2\Vert~t~\Vert\Ua\Vert_2~\Vert\Va\Vert_2\left[2+c_{\Ua}~\Vert \pi_1-\pi_2\Vert~t~\Vert\Ua\Vert_2~\Vert\Va\Vert_2\right]
$$
from which we conclude that
\begin{equation}\label{estimate-Wa-2}
{\Vert\Wa_s(\pi_1,\pi_2)-\Wa_s(\pi_2)\Vert_2}\leq 
c_{\Ua}~\Vert \pi_1-\pi_2\Vert~t^2~\Vert\Ua\Vert_2~\Vert\Va\Vert_2\left[2+c_{\Ua}~\Vert \pi_1-\pi_2\Vert~t^2~\Vert\Ua\Vert_2~\Vert\Va\Vert_2\right].
\end{equation}

The inequality in (\ref{comparison-Wa}) implies that
\begin{eqnarray*}
\Wa_t(\pi_1)^{-1}&\leq& \Wa_t(\pi_2)^{-1/2}~\left[\Wa_t(\pi_2)^{1/2}~
\Wa_t(\pi_1,\pi_2)^{-1} \Wa_t(\pi_2)^{1/2}\right]~\Wa_t(\pi_2)^{-1/2}\\
&=& \Wa_t(\pi_2)^{-1/2}~\left[\Wa_t(\pi_2)^{-1/2}~
\Wa_t(\pi_1,\pi_2)~ \Wa_t(\pi_2)^{-1/2}\right]^{-1}~\Wa_t(\pi_2)^{-1/2}.
  \end{eqnarray*}
  On the other hand we have
  $$
  \begin{array}{l}
  \left[\Wa_t(\pi_2)^{-1/2}~
\Wa_t(\pi_1,\pi_2)~ \Wa_t(\pi_2)^{-1/2}\right]^{-1}\\
\\
\qquad\qquad=\left[Id-\Wa_t(\pi_2)^{-1/2}\left\{\Wa_t(\pi_2)-
\Wa_t(\pi_1,\pi_2)\right\} \Wa_t(\pi_2)^{-1/2}\right]^{-1}.
  \end{array}
  $$
This yields the estimate
\begin{eqnarray*}
\Wa_t(\pi_2)^{1/2}~\Wa_t(\pi_1)^{-1}~\Wa_t(\pi_2)^{1/2}
&\leq &\sum_{n\geq 0}\left[\Wa_t(\pi_2)^{-1/2}\left\{\Wa_t(\pi_2)-
\Wa_t(\pi_1,\pi_2)\right\} \Wa_t(\pi_2)^{-1/2}\right]^{n}.
  \end{eqnarray*}
Combining (\ref{estimate-Wa-1}) with 
(\ref{estimate-Wa-2}), for any $\epsilon>0$ there exists some $\delta(t,\epsilon,\pi_2)>0$ such that
$$
\Vert \pi_1-\pi_2\Vert\leq \delta(t,\epsilon)~\Longrightarrow~
\Vert \Wa_t(\pi_2)^{-1/2}\left\{\Wa_t(\pi_2)-
\Wa_t(\pi_1,\pi_2)\right\} \Wa_t(\pi_2)^{-1/2}\Vert_2\leq 1-\epsilon.
$$
This ends the proof of the lemma.
\qed

\subsection*{Proof of the bias estimate (\ref{last-nylstrom})}\label{proof-nylstrom}

Observe that if $Z\sim \mathcal{N}(0,Q)$ is Gaussian, then the conditional distribution of $Z_{\Pa^c}=(Z_k)_{k\in \Pa^c}$ given $Z_{\Pa}=(Z_k)_{k\in \Pa}$ is again a centred Gaussian with covariance matrix
$$
T_{\Pa}(Q)=Q_{\Pa^c}-Q_{\Pa^c,\Pa}~Q_{\Pa}^{-}~Q_{\Pa,\Pa^c}
$$
where $Q_{\Pa}^{-}$ stands for the Moore-Penrose pseudo-inverse of $Q_{\Pa}$.
The matrix
$T_{\Pa}(Q)$ can be seen as the Schur complement of $Q_{\Pa}$ in $Q$. This shows that
$$
Q-T(Q)=L_{\Pa^c}\odot T_{\Pa}(Q).
$$

 In this notation we have
 \begin{eqnarray*}
\zeta ~\zeta^{\prime}&=&\left[\xi^1_0-m_0,\ldots,\xi^N_0-m_0\right]\left[\begin{array}{c}
(\xi^1_0-m_0)^{\prime}\\
\vdots\\
(\xi^N_0-m_0)^{\prime}\end{array}
\right]~=~\sum_{1\leq i\leq N}~(\xi^N_0-m_0)(\xi^N_0-m_0)^{\prime}\\
&=&\left[\begin{array}{c}
\zeta_1\\
\vdots\\
\zeta_r\end{array}
\right]\left[
\zeta_1^{\prime},
\ldots,
\zeta_r^{\prime}
\right]~=~\left[\begin{array}{ccc}
\langle\zeta_1,\zeta_1\rangle&\ldots&\langle\zeta_1,\zeta_r\rangle\\
\vdots\\
\langle\zeta_r,\zeta_1\rangle&\ldots&\langle\zeta_r,\zeta_r\rangle\end{array}
\right].
 \end{eqnarray*}

 We let $g$ be the matrix
 $$
 \forall 1\leq i,j\leq s\qquad
 g_{i,j}~:=~\langle V_i,V_j\rangle~~\Longleftrightarrow~~ g ~=~\left[
 \begin{array}{c}
 V_1^{\prime}\\
 \vdots\\
 V^{\prime}_s
 \end{array}
 \right]\left[
 V_1,\ldots,V_s
 \right].
 $$
 Also let $g^-=(g^{i,j})_{1\leq i,j\leq s}$ be the pseudo-inverse of $g$. The orthogonal projection of
 a vector $\zeta_l$ with $l\not\in \Pa$ is given by
 \begin{eqnarray*}
 \proj_{\Va}(\zeta_l)&=&\sum_{1\leq i\leq s}\langle \sum_{1\leq j\leq s}~g^{i,j}~V_j,~\zeta_l\rangle~V_i~=~\left[
 V_1,\ldots,V_s
 \right]g^{-}\left[
 \begin{array}{c}
 V_1^{\prime}\\
 \vdots\\
 V^{\prime}_s
 \end{array}
 \right]\zeta_l~:=~\Ta_{\Va}\zeta_l 
 \end{eqnarray*}
 
\begin{eqnarray*}
N~T(p_0)&=&\left[\begin{array}{c}
\Ta_{\Va}\zeta_l \\
\vdots\\
\Ta_{\Va}\zeta_r\end{array}
\right]\left[
(\Ta_{\Va}\zeta_l)^{\prime},
\ldots,
(\Ta_{\Va}\zeta_l)^{\prime}
\right]~=~\left[\begin{array}{ccc}
\langle \Ta\zeta_1,\Ta_{\Va}\zeta_1\rangle&\ldots&\langle\Ta\zeta_1,\Ta_{\Va}\zeta_r\rangle\\
\vdots\\
\langle\Ta\zeta_r,\Ta_{\Va}\zeta_1\rangle&\ldots&\langle \Ta_{\Va}\zeta_r,\Ta_{\Va}\zeta_r\rangle\end{array}
\right]
\\
&=&
\left[\begin{array}{ccc}
\langle\zeta_1,\Ta_{\Va}\zeta_1\rangle&\ldots&\langle \zeta_1,\Ta_{\Va}\zeta_r\rangle\\
\vdots\\
\langle\zeta_r,\Ta _{\Va}\zeta_1\rangle&\ldots&\langle \zeta_r,\Ta_{\Va}\zeta_r\rangle\end{array}
\right]~=~\zeta \left(\Ta_{\Va}\zeta\right)^{\prime}.
\end{eqnarray*}
Given $\Va$, the $N$ random vectors $\zeta^i_{\Pa^c}=(\zeta_k^i)_{k\not\in \Pa}\in \RR^{r-s}$, with $1\leq i\leq N$ are independent random vectors in $\RR^N$
with mean
$$
\EE\left(\zeta^i_{\Pa^c}\vert\Va\right)=Q_{\Pa^c,\Pa}Q_{\Pa}^{-1}~\zeta^i_{\Pa}\quad
\mbox{\rm with}\quad ~\zeta^i_{\Pa}:=(\zeta_k^i)_{k\in \Pa}\in \RR^{s}
$$ 
and covariance matrix
$$
\EE\left(\left[\zeta^i_{\Pa^c}-\EE\left(\zeta^i_{\Pa^c}\vert\Va\right)\right]\left[\zeta^i_{\Pa^c}-\EE\left(\zeta^i_{\Pa^c}\vert\Va\right)\right]^{\prime}\vert\Va\right)=
Q_{\Pa^c}-Q_{\Pa^c,\Pa}Q_{\Pa}^{-1}Q_{\Pa,\Pa^c}.
$$
This implies that for any $k,l\not\in \Pa$ we have
\begin{eqnarray*}
\EE\left(\langle\zeta_k,\Ta _{\Va}\zeta_l\rangle\vert\Va\right)&=&\sum_{1\leq i,j\leq N}~\EE\left(
\zeta_k^i~\Ta _{\Va}(i,j)\zeta_l^j\vert\Va\right)\\
&=&\sum_{1\leq i,j\leq N}~\EE\left(
\left[\zeta_k^i-\EE\left(\zeta^i_{k}\vert\Va\right)\right]~\Ta _{\Va}(i,j)~\left[\zeta_l^j-\EE\left(\zeta^j_{l}\vert\Va\right)
\right]~\vert\Va\right)\\
&&\hskip3cm+\sum_{1\leq i,j\leq N}\EE\left(\zeta^i_{k}\vert\Va\right)~\Ta _{\Va}(i,j)~\EE\left(\zeta^j_{l}\vert\Va\right)\\
&=& \tr(\Ta _{\Va})~ (Q_{\Pa^c}-Q_{\Pa^c,\Pa}Q_{\Pa}^{-1}Q_{\Pa,\Pa^c})(k,l)\\
&&\\
&&\hskip.3cm+\sum_{1\leq i,j\leq N}\left(Q_{\Pa^c,\Pa}Q_{\Pa}^{-1}~\zeta^i_{\Pa}\right)(k)~\Ta _{\Va}(i,j)~
\left(Q_{\Pa^c,\Pa}Q_{\Pa}^{-1}~\zeta^j_{\Pa}\right)(l).
\end{eqnarray*}
On the other hand, we have
$$
\begin{array}{l}
\sum_{1\leq i,j\leq N}\left(Q_{\Pa^c,\Pa}Q_{\Pa}^{-1}~\zeta^i_{\Pa}\right)(k)~\Ta _{\Va}(i,j)~
\left(Q_{\Pa^c,\Pa}Q_{\Pa}^{-1}~\zeta^j_{\Pa}\right)(l)\\
\\
\qquad\qquad=\sum_{u,v\in \Pa}
\left(Q_{\Pa^c,\Pa}Q_{\Pa}^{-1}\right)(k,u)~\sum_{1\leq i,j\leq N}~\left[\zeta^i_{u}~\Ta _{\Va}(i,j)~
~\zeta^j_{v}\right]~\left(Q_{\Pa^c,\Pa}Q_{\Pa}^{-1}\right)(l,v)\\
\\
\qquad\qquad=\sum_{u,v\in \Pa}
\left(Q_{\Pa^c,\Pa}Q_{\Pa}^{-1}\right)(k,u)~\left\langle\zeta_{u},\Ta _{\Va}
~\zeta_{v}\right\rangle~\left(Q_{\Pa^c,\Pa}Q_{\Pa}^{-1}\right)(l,v)\\
\\
\qquad\qquad=\sum_{u,v\in \Pa}
\left(Q_{\Pa^c,\Pa}Q_{\Pa}^{-1}\right)(k,u)~\left\langle\zeta_{u},
~\zeta_{v}\right\rangle~\left(Q_{\Pa^c,\Pa}Q_{\Pa}^{-1}\right)(l,v).
\end{array}
$$
Taking the expectation we find that
\begin{eqnarray*}
\EE\left(\langle\zeta_k,\Ta _{\Va}\zeta_l\rangle\right)
&=& s~ (Q_{\Pa^c}-Q_{\Pa^c,\Pa}Q_{\Pa}^{-1}Q_{\Pa,\Pa^c})(k,l)\\
&&\\
&&\hskip.3cm+N~\sum_{u,v\in \Pa}
\left(Q_{\Pa^c,\Pa}Q_{\Pa}^{-1}\right)(k,u)~Q_{\Pa}(u,v)~\left(Q_{\Pa^c,\Pa}Q_{\Pa}^{-1}\right)(l,v)\\
&=&s~ (Q_{\Pa^c}-Q_{\Pa^c,\Pa}Q_{\Pa}^{-1}Q_{\Pa,\Pa^c})(k,l)+N
\left[Q_{\Pa^c,\Pa}Q_{\Pa}^{-1}Q_{\Pa,\Pa^c}\right](k,l)\\
&=&\left\{s~ Q_{\Pa^c}+(N-s)~\left[Q_{\Pa^c,\Pa}Q_{\Pa}^{-1}Q_{\Pa,\Pa^c}\right]\right\}(k,l).
\end{eqnarray*}
This shows that
$$
(\ref{nylstrom-1})\Longrightarrow~\begin{array}[t]{rcl}
\EE\left(T(p_0)\right)&=&\left(J-L_{\Pa^c}\right)\odot Q+
L_{\Pa^c}\odot \left(\frac{s}{N}~Q_{\Pa^c}+\left(1-\frac{s}{N}\right)~\left[Q_{\Pa^c,\Pa}Q_{\Pa}^{-1}Q_{\Pa,\Pa^c}\right]\right)\\
&&\\
&=&\displaystyle T(P_0)+\frac{s}{N}~L_{\Pa^c}\odot ~\Ta_{\Pa}(Q).
\end{array}
$$
This ends the proof of (\ref{last-nylstrom}).\qed


\begin{thebibliography}{99.}


\bibitem{abou-kandil03}
H. Abou-Kandil, G. Freiling, V. Ionescu, and G. Jank. Matrix Riccati Equations in Control and Systems Theory. Birkhauser Verlag (2003).

\bibitem{abou-kandil94}
H. Abou-Kandil, G. Freiling, and G. Jank. Solution and Asymptotic Behavior of Coupled Riccati Equations in Jump Linear Systems. IEEE Transactions on Automatic Control. vol. 34, no. 8. pp. 1631--1636 (1994).

\bibitem{anderson2003}
J.L. Anderson. A local least squares framework for ensemble filtering. Monthly Weather Review. vol. 131, no. 4. pp. 634--642 (2003).

\bibitem{anderson07}
J.L. Anderson. An adaptive covariance inflation error correction algorithm for ensemble filters. Tellus A. vol. 59, no. 2. pp. 210-224 (2007).

\bibitem{anderson09}
J.L. Anderson. Spatially and temporally varying adaptive covariance inflation for ensemble filters. Tellus A. vol. 61, no. 1. pp. 72--83 (2009).

\bibitem{anderson2012}
J.L. Anderson. Localization and Sampling Error Correction in Ensemble Kalman Filter Data Assimilation. Monthly Weather Review. vol. 140, no. 7. pp. 2359--2371 (2012).

\bibitem{anderson99}
J.L. Anderson and S.L. Anderson. A Monte Carlo Implementation of the Nonlinear Filtering Problem to Produce Ensemble Assimilations and Forecasts. Monthly Weather Review. vol. 127, no. 12. pp. 2741--2758 (1999).

\bibitem{Antsaklis}
P.J. Antsaklis and A.N. Michel. A Linear Systems Primer. Birkhäuser, Boston (2007).

\bibitem{arcolano}
N. Arcolano and P.J. Wolfe. Estimating principal components of covariance matrices using the Nystr\"om method. arXiv e-print, \href{https://arxiv.org/abs/1111.6926}{\tt arXiv:1111.6926} (2011).

\bibitem{bailey}
R.A. Bailey. Association Schemes: Designed Experiments, Algebra and Combinatorics. Cambridge University Press (2004).

\bibitem{barker}
G.P. Barker, L.Q. Eifler, and T.P. Kezlan. A non-commutative spectral theorem. Linear Algebra and Its Applications. vol. 20, no. 2. pp. 95--100 (1978).

\bibitem{basar2008}
T. Basar and P. Bernhard. H-infinity Optimal Control and Related Minimax Design Problems: A Dynamic Game Approach. Birkhauser Boston (2008).

\bibitem{Bernhard79}
P. Bernhard. Linear-quadratic, two-person, zero-sum differential games: Necessary and sufficient conditions. Journal of Optimization Theory and Applications. vol. 27, no. 1. pp. 51-69 (1979).

\bibitem{bickel}
P.J. Bickel and E. Levina. Regularized estimation of large covariance matrices. The Annals of Statistics. vol. 36, no. 1. pp. 199--227 (2008).

\bibitem{bien}
J. Bien and R. Tibshirani. Sparse Estimation of a Covariance Matrix. Biometrika. vol. 98, no. 4. pp. 807--820 (2010).

\bibitem{ap-2016}
A.N. Bishop and P. Del Moral. On the Stability of Kalman-Bucy Diffusion Processes. SIAM Journal on Control and Optimization. vol. 55, no. 6. pp 4015--4047 (2017); arxiv e-print \href{https://arxiv.org/abs/1610.04686}{\tt arXiv:1610.04686} updated.

\bibitem{2018arXiv180800235}
A.N. Bishop, P. Del Moral. On the Stability of Matrix-Valued Riccati Diffusions. arXiv e-print, \href{https://arxiv.org/abs/1808.00235}{\tt arXiv:1808.00235} (2018).

\bibitem{ap-2018-franklin}
A.N. Bishop and  P. Del Moral. On the robustness of Riccati flows to complete model misspecification. Journal of the Franklin Institute. vol. 355, no. 15. pp 7178--7200 (2018). 

\bibitem{apa-2017}
A.N. Bishop, P. Del Moral, and A. Niclas. A perturbation analysis of stochastic matrix Riccati diffusions. arXiv e-print, \href{https://arxiv.org/abs/1709.05071}{\tt arXiv:1709.05071} (2017).

\bibitem{bose}
R.C. Bose and D.M. Mesner. On linear associative algebras corresponding to association schemes of partially balanced designs. The Annals of Mathematical Statistics. vol. 30, no. 1. pp. 21--38 (1959).

\bibitem{brouwer}
A.E. Brouwer, A.M. Cohen, and A. Neumair. Distance Regular Graphs. Springer (1989).

\bibitem{callier81}
F.M. Callier and J.L. Willems. Criterion for the Convergence of the Solution of the Riccati Differential Equation. IEEE Transactions on Automatic Control. vol. 26, no. 6. pp. 1232--1242 (1981).

\bibitem{callier96}
F.M. Callier and J.J. Winkin. Asymptotic behaviour of the solution of the projection Riccati differential equation. IEEE Trans. on Automatic Control. vol. 41, no. 5. pp. 646--659 (1996).

\bibitem{cameron}
P. J. Cameron. Coherent Configurations, Association Schemes and Permutation Groups. In Proc. of Groups, Combinatorics \& Geometry (Durham, 2001), pp. 55--71. World Scientific Publishing, River Edge, NJ (2003).

\bibitem{chaudhuri}
S. Chaudhuri, M. Drton and T.S. Rochardson. Estimation of a covariance matrix with zeros. Biometrika, vol. 94, no. 1. pp. 199--216 (2007).

\bibitem{chen}
R.Y. Chen, A. Gittens and J.A. Tropp. The Masked Sample Covariance Estimator: An Analysis via Matrix Concentration Inequalities. arXiv e-print, \href{https://arxiv.org/abs/1109.1637}{\tt arXiv:1109.1637} (2011). 

\bibitem{dahl}
J. Dahl, V. Roychowdhury, and L. Vandenberghe. Maximum likelihood estimation of Gaussian graphical models: Numerical implementation and topology selection. \href{http://www.ee.ucla.edu/~vandenbe/covsel.html}{\tt UCLA Preprint}, September (2005).

\bibitem{desouza90}
C.E. de Souza and M.D. Fragoso. On the existence of maximal solution for generalized algebraic Riccati equations arising in stochastic control. Systems \& Control Letters. vol. 14, no. 3. pp. 233--239 (1990).

\bibitem{delmoral16b}
P. Del Moral, A. Kurtzmann, and J. Tugaut. On the stability and the uniform propagation of chaos of extended ensemble Kalman-Bucy filters. arXiv e-print, \href{http://arxiv.org/pdf/1606.08256.pdf}{\tt arXiv:1606.08256} (2016).

\bibitem{delfour07}
M.C. Delfour. Linear quadratic differential games: Saddle point and Riccati differential equation. SIAM Journal on Control and Optimization. vol. 46, no. 2. pp. 750-774 (2007).

\bibitem{delmoral16a}
P. Del Moral and J. Tugaut. On the stability and the uniform propagation of chaos properties of ensemble Kalman-Bucy filters. arXiv e-print, \href{https://arxiv.org/pdf/1605.09329.pdf}{\tt arXiv:1605.09329}, to appear in the Annals of Applied Probability (2017).

\bibitem{delsarte}
P. Delsarte and V.I. Levenshtein. Association schemes and coding theory. IEEE Transactions on Information Theory. vol.44, no. 6. pp. 2477--2504 (1998).

\bibitem{elkaroui}
N. El Karoui. Operator norm consistent estimation of large-dimensional sparse covariance matrices. The Annals of Statistics. vol. 36, no. 6. pp. 2717--2756 (2008).

\bibitem{evensen03}
G. Evensen. The Ensemble Kalman Filter: theoretical formulation and practical implementation. Ocean Dynamics. vol. 53, no. 4. pp. 343--367 (2003).

\bibitem{freiling96}
G. Freiling, G. Jank, and H. Abou-Kandil. Generalized Riccati difference and differential equations. Linear Algebra and Its Applications. vols. 241--243. pp. 291--303 (1996).

\bibitem{Gaspari1999}
G. Gaspari and S.E. Cohn. Construction of correlation functions in two and three dimensions. Quarterly Journal of the Royal Meteorological Society. vol. 125, no. 554. pp. 723-757 (1999).

\bibitem{Gaspari2006}
G. Gaspari, S.E. Cohn, J. Guo and S. Pawson. Construction and application of covariance functions with variable length-fields. Quarterly Journal of the Royal Meteorological Society. vol. 132, no. 619. pp. 1815-1838 (2006).	

\bibitem{gelb74}
A. Gelb (editor). Applied Optimal Estimation. MIT Press (1974).

\bibitem{givens}
C.R. Givens and R.M. Shortt. A class of Wasserstein metrics for probability distributions. Michigan Math. J. vol. 31, no. 2. pp. 231--240 (1984).

\bibitem{grossmann}
H. Grossmann. Automating the analysis of variance of orthogonal designs. Computational Statistics and Data Analysis. vol. 70. pp. 1--18  (2014).

\bibitem{Haff}
L.R. Haff. Empirical Bayes estimation of the multivariate normal covariance matrix. The Annals of Statistics. vol. 8, no. 3. pp. 586--597 (1980).

\bibitem{hamill01}
T.M. Hamill, J.S. Whitaker, and C. Snyder. Distance-Dependent Filtering of Background Error Covariance Estimates in an Ensemble Kalman Filter. Monthly Weather Review. vol. 129, no. 11. pp. 2776--2790 (2001).

\bibitem{Heemink2001}
A.W. Heemink, M. Verlaan, and A.J. Segers. Variance reduced ensemble Kalman filtering. Monthly Weather Review. vol.129, no. 7. pp. 1718--1728 (2001).

\bibitem{higham}
N.J. Higham. Functions of Matrices: Theory and Computation. SIAM, Philadelphia (2008).

\bibitem{hou16}
E. Hou, E. Lawrence, and A.O. Hero. Penalized Ensemble Kalman Filters for High Dimensional Non-linear Systems. arXiv e-print, \href{https://arxiv.org/abs/1610.00195}{\tt arXiv:1610.00195} (2016).

\bibitem{Houtekamer2001}
P.L. Houtekamer and H.L. Mitchell. A Sequential Ensemble Kalman Filter for Atmospheric Data Assimilation. Monthly Weather Review. vol. 129, no. 1. pp. 123--137 (2001).

\bibitem{johns2008}
C.J. Johns and J. Mandel. A two-stage ensemble Kalman filter for smooth data assimilation. Environmental and Ecological Statistics. vol. 15, no. 1. pp. 101-110 (2008).

\bibitem{johnson}
C.R. Johnson. Partitioned and Hadamard Product Matrix Inequalities. Journal of Research of the National Bureau of Standards. vol. B3, no. 6. pp. 585--591 (1978).

\bibitem{kalman61}
R.E. Kalman and R.S. Bucy. New Results in Linear Filtering and Prediction Theory. Journal of Basic Engineering. vol. 83, no. 1. pp. 95--108 (1961).

\bibitem{khare}
K. Khare and B. Rajaratnam. Wishart distributions for decomposable covariance graph models. The Annals of Statistics. vol. 39, no. 1. pp. 514--555 (2011).

\bibitem{kucera72}
V. Kucera. A Contribution to Matrix Quadratic Equations. IEEE Transactions on Automatic Control. vol. 17, no. 3. pp. 344--347 (1972).

\bibitem{Kwakernaak72}
H. Kwakernaak and R. Sivan. Linear Optimal Control Systems. Wiley-Interscience (1972).

\bibitem{lam-fan}
C. Lam and J. Fan. Sparsistency and rates of convergence in large covariance matrix estimation. The Annals of Statistics, vol. 37, no. 6B. pp. 4254--4278 (2009).

\bibitem{Lancaster1995}
P. Lancaster and L. Rodman. Algebraic Riccati Equations. Oxford University Press (1995).

\bibitem{law2016}
K.J.H. Law, H. Tembine and R. Tempone. Deterministic Mean-Field Ensemble Kalman Filtering. SIAM Journal on Scientific Computing. vol. 38, no. 3. pp. A1251-A1279 (2016).

\bibitem{ledoit2}
O. Ledoit and M. Wolf. A well-conditioned estimator for large-dimensional covariance matrices. Journal of Multivariate Analysis. vol. 88, no. 2. pp. 365--411 (2004).

\bibitem{legland09}
F. Le Gland, V. Monbet, V.-D. Tran. Large sample asymptotics for the ensemble Kalman filter. Research Report: RR-7014, INRIA. $<$inria-00409060$>$ (2009).

\bibitem{levina}
E. Levina and R. Vershynin. Partial estimation of covariance matrices. Probability Theory and Related Fields. vol. 153, no. 3. pp. 405--419 (2012). Available on arXiv in 2010.

\bibitem{li2009}
H. Li, E. Kalnay, and T. Miyoshi. Simultaneous estimation of covariance inflation and observation errors within an ensemble Kalman filter. Quarterly Journal of the Royal Meteorological Society. vol. 135, no. 639. pp. 523--533 (2009).

\bibitem{kelly14}
D.T.B. Kelly, K.J. Law, and A.M. Stuart. Well-posedness and accuracy of the ensemble Kalman filter in discrete and continuous time. Nonlinearity. vol. 27, no. 10. pp. 2579--2603 (2014).

\bibitem{Majda1606.09087}
A.J. Majda and X.T. Tong. Rigorous accuracy and robustness analysis for two-scale reduced random Kalman filters in high dimensions. \href{https://arxiv.org/abs/1606.09087}{\tt arXiv:1606.09087} (2016).

\bibitem{Majda16}
A.J. Majda and X.T. Tong. Performance of Ensemble Kalman filters in large dimensions. \href{https://arxiv.org/abs/1606.09321}{\tt arXiv:1606.09321} (2016).

\bibitem{mandel2011}
J. Mandel, L. Cobb, and J.D. Beezley. On the convergence of the ensemble Kalman filter. Applications of Mathematics. vol. 56, no. 6. pp. 533--541 (2011).

\bibitem{McAsey06}
M. McAsey and L. Mou. Generalized Riccati equations arising in stochastic games. Linear Algebra and Its Applications. vol. 416, no. 2-3. pp. 710--723 (2006).

\bibitem{Mitchell2002}
H.L. Mitchell, P.L. Houtekamer and G. Pellerin. Ensemble size, balance, and model-error representation in an ensemble Kalman filter. Monthly Weather Review. vol. 130, no. 11. pp. 2791--2808 (2002).

\bibitem{Molinari77}
B.P. Molinari. The time-invariant linear-quadratic optimal control problem. Automatica. vol. 13, no. 4. pp. 347--357 (1977).

\bibitem{nelder-1}
J.A. Nelder. The analysis of randomized experiments with orthogonal block structure. I. Block structure and the null analysis of variance. Proceedings of the Royal Society of London, Series A. vol. 283, no. 1393. pp. 147--162 (1965).

\bibitem{nelder}
 J.A. Nelder. The analysis of randomized experiments with orthogonal block structure. II. Treatment structure and the general analysis of variance. Proceedings of the Royal Society of London, Series A. vol. 283, no. 1393. pp. 163--178 (1965).

\bibitem{olkin}
I. Olkin and F. Pukelsheim. The distance between two random vectors with given dispersion matrices. Linear Algebra and Its Applications. vol. 48. pp. 257--263 (1982).

\bibitem{pham2001}
D.T. Pham. Stochastic methods for sequential data assimilation in strongly nonlinear systems. Monthly Weather Review. vol. 129, no. 5. pp. 1194--1207 (2001).

\bibitem{Rebeschini2015}
P. Rebeschini and R. Van Handel. Can local particle filters beat the curse of dimensionality?. The Annals of Applied Probability. vol. 25, no. 5. pp. 2809-2866 (2015).

\bibitem{Reich2013}
S. Reich and C.J. Cotter. Ensemble filter techniques for intermittent data assimilation. In Large Scale Inverse Problems: Computational Methods and Applications in the Earth Sciences (eds: M. Cullen, M.A. Freitag, S. Kindermann, R. Scheichl). pp. 91--134. De Gruyter Publishers (2013). See also: arXiv e-print, \href{https://arxiv.org/abs/1208.6572}{\tt arXiv:1208.6572} (2012).

\bibitem{sakov11}
P. Sakov and L. Bertino. Relation between two common localisation methods for the {\tt EnKF}. Computational Geosciences. vol. 15, no. 2. pp. 225-237 (2011). 

\bibitem{sakov2008a}
P. Sakov and P.R. Oke. A deterministic formulation of the ensemble Kalman filter: an alternative to ensemble square root filters. Tellus A. vol. 60, no. 2. pp. 361-371 (2008).

\bibitem{Saetrom11}
J. Sætrom and H. Omre. Ensemble Kalman filtering with shrinkage regression techniques. Computational Geosciences. vol. 15, no. 2. pp. 271--292 (2011).

\bibitem{Taghvaei2016ACC}
A. Taghvaei and P.G. Mehta. An optimal transport formulation of the linear feedback particle filter. In Proc. of the 2016 American Control Conference (ACC), Boston, USA (July, 2016). 

\bibitem{Tenenbaum2000}
J.B. Tenenbaum. V. De Silva, and J.C. Langford. A Global Geometric Framework for Nonlinear Dimensionality Reduction. Science. vol. 290, no. 5500. pp. 2319--2323 (2000).

\bibitem{tippett03}
M.K. Tippett, J.L. Anderson, C.H. Bishop, T.M. Hamill, and J.S. Whitaker. Ensemble square root filters. Monthly Weather Review. vol. 131, no. 7. pp. 1485--1490 (2003).

\bibitem{tong2016a}
X.T. Tong, A.J. Majda, and D. Kelly. Nonlinear stability and ergodicity of ensemble based Kalman filters. Nonlinearity. vol. 29, no. 2. pp 657--691 (2016).

\bibitem{tong16b}
X.T. Tong, A.J. Majda, and D. Kelly. Nonlinear stability of the ensemble Kalman filter with adaptive covariance inflation. Communications in Mathematical Sciences. vol. 14, no. 5. pp. 1283--1313 (2016).

\bibitem{hemmen}
J.L. van Hemmen and T. Ando. An inequality for trace ideals. Communications in Mathematical Physics. vol. 76, no. 143. pp. 143--148 (1980).

\bibitem{wagaman}
A.S. Wagaman and E. Levina. Discovering Sparse Covariance Structures with the Isomap. Journal of Computational and Graphical Statistics. vol. 18, no. 3. pp. 551--572 (2009).

\bibitem{wonham68}
W.M. Wonham. On a Matrix Riccati Equation of Stochastic Control. SIAM Journal of Control. vol. 6, no. 4. pp 681--697 (1968).



\end{thebibliography}
\end{document}